\documentclass[12pt]{article}
\usepackage[left=2cm,top=2cm,right=2cm,bottom=2cm]{geometry}

\usepackage[cp1251]{inputenc}
\usepackage[british]{babel}
\usepackage{wrapfig}    % for wrapfigure environment

\usepackage{graphicx,epsfig,epstopdf}
\usepackage{amsfonts,amsmath,amssymb,latexsym,xcolor,mathtools}
\usepackage{hyperref,enumitem, float}
\usepackage{amsthm,verbatim}
\usepackage[scr=boondox]{mathalfa}
 \usepackage{cite}
 \usepackage[colorinlistoftodos]{todonotes}

\numberwithin{equation}{section}
\pdfoutput=1

\theoremstyle{plain}

\theoremstyle{definition}

\newenvironment{exa}
{%
	\pushQED{\qed}\begin{exa/}}
	{\popQED\end{exa/}}
\newenvironment{rema}
{%
	\pushQED{\qed}\begin{rema/}}
	{\popQED\end{rema/}}
\newtheorem{defi/}[theorem]{Definition}
\newtheorem{rema/}[theorem]{Remark}
\newtheorem{exa/}[theorem]{Example}

\newcommand{\conj}{\dag}

\newcommand{\ptl}{\partial}

\newcommand{\bda}{\boldsymbol{a}}
\newcommand{\bdeta}{\boldsymbol{\eta}}
\newcommand{\bdzeta}{\boldsymbol{\zeta}}
\newcommand{\bdalpha}{\boldsymbol{\alpha}}

\newcommand{\bdG}{\boldsymbol{\Gamma}}

\newcommand{\bdxi}{\boldsymbol{\xi}}

\newcommand{\tva}{\text{\tmverbatim{a}}}
\newcommand{\tvb}{\text{\tmverbatim{b}}}
\newcommand{\tvc}{\text{\tmverbatim{c}}}
\newcommand{\bdm}{\boldsymbol{m}}
\DeclareMathOperator{\trace}{tr}

\usepackage{tikz}

\makeatletter
\newcommand*{\transpose}{%
	{\mathpalette\@transpose{}}%
}
\newcommand*{\@transpose}[2]{%
	\raisebox{\depth}{$\m@th#1\intercal$}%
}
\makeatother

\newcommand{\assign}{:=}
\newcommand{\cdummy}{\cdot}
\newcommand{\infixand}{\text{ and }}
\newcommand{\mathd}{\mathrm{d}}
\newcommand{\nobracket}{}

\newcommand{\tmcolor}[2]{{\color{#1}{#2}}}
\newcommand{\tmem}[1]{{\em #1\/}}
\newcommand{\tmmathbf}[1]{\ensuremath{\boldsymbol{#1}}}
\newcommand{\tmop}[1]{\ensuremath{\operatorname{#1}}}

\newcommand{\tmtextit}[1]{\text{{\itshape{#1}}}}
\newcommand{\tmverbatim}[1]{\text{{\ttfamily{#1}}}}

\pdfoutput=1

\ifpdf
\DeclareGraphicsExtensions{.eps,.pdf,.png,.jpg}
\else
\DeclareGraphicsExtensions{.eps}
\fi

\definecolor{myred}{RGB}{160,0,0}
\definecolor{mygreen}{RGB}{0,160,0}
\definecolor{myblue}{RGB}{0,0,160}
\hypersetup{
	colorlinks = true,
	linkcolor = {myred},
	citecolor = {mygreen},
	urlcolor  = {myblue},
}
\newcommand{\RED}{} % Without red comments

\newcommand{\GREEN}{} % Without blue comments

\title{A note on double Floquet-Bloch transforms\\
	and the far-field asymptotics of Green's functions\\
	tailored to periodic structures}
\author{ Andrey V. Shanin$^{\dagger}$, Rapha\"{e}l C. Assier$^{*}$, Andrey I. Korolkov$^{*}$, Oleg I. Makarov$^{\dagger}$\\
	\footnotesize{$^{\dagger}$ Department of Physics (Acoustics Division), Moscow State University, Leninskie Gory, {\rm 119992}, Moscow, Russia}\\
	\footnotesize{$^{*}$ Department of Mathematics, University of Manchester, Oxford Road, Manchester, {\rm M13 9PL}, UK}	
}

\begin{document}	

%\title[{On double Fourier integrals}]{A contribution to the mathematical theory
%	of diffraction:\\ a note on double Fourier integrals}
%
%\author[R.~C.~Assier, A.~V.~Shanin \& A.~I.~Korolkov]{R.~C.~ASSIER}
%\address{Department of Mathematics, University of Manchester, \\ Oxford Road, Manchester, {\rm M13 9PL}, UK}
%\extraauthor{A.~V.~SHANIN and A.~I.~KOROLKOV }
%\extraaddress{Department of Physics (Acoustics Division), Moscow State University, \\ Leninskie Gory, {\rm 119992}, Moscow, Russia}
%
%
%\received{\recd ??. \revd ?? }

\maketitle

%\eqnobysec

\begin{abstract}
	We propose a general procedure to study double integrals arising when considering wave propagation in periodic structures. This method, based on a complex deformation of the integration surface to bypass the integrands' singularities, is particularly efficient to estimate the Green's functions of such structures in the far field. We provide several illustrative examples and explicit asymptotic formulae. Special attention is devoted to the pathological cases of degeneracies, such as Dirac conical points for instance.
\end{abstract}

%\tableofcontents

\section{Introduction}

The study of wave propagation within periodic structures has interested
mathematicians and physicists for a long time {\cite{brillouin1946wave}}. The
emergence of metamaterials has contributed to its recent surge in popularity.
Initially focused on electro-magnetic applications (e.g.\ photonics
\cite{Shalaev2007,Soukoulis2007,Zheludev2010} or cloaking \cite{Cummer2006,Cai2007}), this field
has grown significantly in the past two decades with a myriad of applications
in elasticity and acoustics where metamaterials are now being used and
manufactured routinely, the term phononics has now gained in popularity. An
excellent, recent and comprehensive review on those so-called mechanical
metamaterials can be found in {\cite{MechMetaReview2023}}.

With this note, we wish to make a connection between our recent work on the
application of multidimensional complex analysis to diffraction theory
{\cite{Assier2018a,Shanin2019SommerfeldTI,Assier2019b,2dcontinuation2021,Assier2019a,Assier2019c,Ice2021,kunz2021diffraction,Part6A,Kunz2023,kunz2023diffractionasymp}}
and the topic of wave propagation in periodic structures. The main similarity
between these two subfields of wave study comes from the fact that, in both
cases the wave field can be written as a multiple integral of the Fourier
type. In diffraction theory, this comes from applying a multiple Fourier
transform to the BVP at hand, while for periodic structure it comes from
applying a multiple Floquet-Bloch transform. The difference lies in the
surface of integration which is $\mathbb{R}^2$ in diffraction theory and $[-
\pi, \pi]^2$ for periodic structures.

Though commonly named Floquet-Bloch transform in the research community
working on waves in periodic structures, this transform is sometimes just
called Bloch transform, generalised Fourier series, discrete Fourier transform, or $\mathcal{Z}$-transform
with a change of variables. For more about it and its properties, we refer the
reader to {\cite{Kuchment1982,Fli-2009,Radosz-2010,LECHLEITER2017}}.

In both diffraction theory and the present topic, a direct link can be made
between the singularities of the integrand and the far-field behaviour of the
physical field. Moreover they share the feature that for a wave number with a
non-zero imaginary part, the integrand is singularity-free on the integration
surface, while, as the wavenumber becomes real, singularities hit the
integration surface. In order to make sense of the physical field, one can
deform the surface of integration into $\mathbb{C}^2$, while avoiding the
singularities of the integrand. We have described this deformation process in
detail in {\cite{Part6A}} in the context of diffraction theory, and we will
show here that it can be applied in the present context. As for diffraction
theory, we show here that only special points on the singularity set of the
integrand are responsible for wave-like components in the far-field.

The topic of asymptotic estimation of multi-dimensional Fourier integrals were addressed before in \cite{Lighthill1960,Lighthill78,Jones1982-xs}. However, it was done in the context of real analysis, and for the case when  studied integrals can be treated as nested ones. Particularly, the case of polar singularity for which 2D Fourier integral can be reduced using residue theorem to a sum of 1D oscillating integrals was studied there. On the other side, the approach described in \cite{Part6A} does not have such restrictions, and is based on the local deformation of integration surfaces in $\mathbb{C}^2$.  It is worth mentioning the recent paper \cite{Levitt2023}, where such deformation is used for the numerical evaluation of Floquet-Bloch integrals, and referred to as a ``complex deformation of the Brillouin zone''.

One important phenomena in periodic structures is the occurrence of
degeneracies. Roughly speaking a degeneracy occurs when, at a given frequency,
the singularity structure of the integrand (also known as the dispersion
diagram) experiences a change in its topology. This can for example be a loop
shrinking to a point or two branches of a hyperbola that start touching each
other. These points are important as they have been shown to lead to very
interesting physical behaviour \cite{MechMetaReview2023}, such as the existence of bandgaps, and can
sometimes be exploited to give rise to so-called topologically protected modes
\cite{Laude2021}. A well known example of a topological mode is a Valley mode which arises in hegaxonal lattices due to the symmetry properties of the so-called Dirac cones corresponding to a degeneracy of the dispersion diagram in a conical point. Such cones are the main feature of unique properties in graphene \cite{CastroNeto2009}, and were recently studied in a mechanical context \cite{Miniaci2021}; they can be exploited practically to achieve acoustic cloaking and collimation (see e.g.\ \cite{Yu-Dirac}). Another example of a degeneracy is a parabolic maxima or minima of the dispersion diagram that can lead to negative refraction at material interfaces, known as a maxon or roton behaviour \cite{Santos2003,Chen2021}. Local minimums can also be exploited practically to achieve collimation (see e.g. \cite{Vanel2023}).

In this article, we will focus on the study of Green's functions specifically
tailored to the periodic structure of the media and solution to the Helmholtz
equation (continuous or discrete) forced at a specific source location and
potentially subjected to boundary conditions on periodically arranged
scatterers. Our specific aim will be to extract the far-field behaviour of
such Green's functions. \RED{The Helmholtz equation is obtained by assuming time-harmonic solutions to the linear wave equation with angular frequency $\omega$. Throughout this article, we use the $e^{-i\omega t}$ convention.} Related work in the discrete setting include
{\cite{Martin2006}} and {\cite{Vanel2016}}. In {\cite{Martin2006}} the
degeneracies are not considered however. In {\cite{Vanel2016}}, an interesting
approach inspired by the technique of high-frequency homogenisation (see, e.g.\
{\cite{Craster2010,Guzina2019,Assier2020HFH,touboul2023highfrequency}}) allows them to derive some
results at specific degeneracies. We will endeavour to show that our approach
can lead to similar results and can also be employed in the continuous
setting.

The rest of the article is organised as follows. In section
\ref{sec:general-framework}, we set the general framework within which we will
work, describing the physical fields we are interested in as double integrals
whose integrand have specific properties. In section
\ref{sec:examples-illustrative}, we show that this framework is indeed
relevant by providing two illustrative examples, one in the discrete setting
of a square lattice (section \ref{sec:discrete-example}) and one in the
continuous setting of a phononic crystal (section
\ref{sec:continuous-example}). In section \ref{sec:deformation-process}, we
explain how one can and should deform the integration surface for real
wavenumbers, before giving a brief but general description on how this can be
used to extract far-field asymptotics in section~\ref{sec:summary-procedure}. Section \ref{sec:asymptotics-special-points} is dedicated to the asymptotics
obtained away from degeneracies, while section \ref{sec:vanishing-loops}
considers a set of specific degeneracies: local extrema of the dispersion
diagram in section \ref{sec:local-max},
hyperbolic degeneracies (crossing with rebuilding) in section
\ref{sec:hyp-degeneracy}, as well a conical (Dirac) points in section~\ref{sec:Dirac}.

\section{General framework}\label{sec:general-framework}
Let us consider wave fields in doubly periodic media excited by a point source that can be written $u (\tmmathbf{r}; \tmmathbf{r}_s ; k,
\varkappa)$, where $\tmmathbf{r}= (x_1, x_2) \in \mathbb{R}^2$ represents the
space variables, $\tmmathbf{r}_s$ represents the source location, and $k + i
\varkappa$ is a scalar parameter chosen such that $k > 0$ and $\varkappa
\geqslant 0$. This parameter is generally interpreted as the wavenumber of the
problem with an added small positive imaginary part.

Let us assume that the media has double periodicity given by two vectors
$\tmmathbf{d}$ and $\tmmathbf{\ell}$, leading to this media being represented
as a repetition of some identical parallelograms. We will refer to the
parallelogram containing the origin $\tmmathbf{r}=\tmmathbf{0}$ as the unit
cell and denote it $S_0$.  Introduce a $2 \times 2$
invertible matrix $\Lambda$ whose columns are the vectors $\tmmathbf{d}$ and
$\tmmathbf{\ell}$ respectively. So that any integer combination of those
vectors $m_1 \tmmathbf{d}+ m_2 \tmmathbf{\ell}$ for some $\tmmathbf{m}= (m_1,
m_2) \in \mathbb{Z}^2$ can be written as $\Lambda \tmmathbf{m}$.

Due to the double periodicity of the media, the wave field $u$ can formally
be written for all $\tmmathbf{m} \in \mathbb{Z}^2$, as the following double
integral
\begin{align}
	u (\tmmathbf{r}+ \Lambda \tmmathbf{m}; \tmmathbf{r}_s ; k, \varkappa) & =
	\frac{1}{4 \pi^2} \iint_{\tmmathbf{\mathcal{B}}} F (\tmmathbf{r};
	\tmmathbf{r}_s ; \tmmathbf{\xi}; k, \varkappa) e^{- i\tmmathbf{m} \cdot
		\tmmathbf{\xi}} \mathd \tmmathbf{\xi},  \label{eq:inverse-BF-transform}
\end{align}
where the surface of integration $\tmmathbf{\mathcal{B}}$ is the square
$\tmmathbf{\mathcal{B}}= [- \pi, \pi]^2 \subset \mathbb{R}^2$ and is known as
the Brillouin zone. Expression (\ref{eq:inverse-BF-transform}) is known as an inverse Floquet-Bloch transform. The integration variable is $\tmmathbf{\xi}= (\xi_1,
\xi_2)$, $\mathd \tmmathbf{\xi}$ is $\mathd \xi_1 \mathd \xi_2$, and
$F$ is the direct Floquet-Bloch transform of $u$ and is defined as
\begin{align}
	F (\tmmathbf{r}; \tmmathbf{r}_s ; \tmmathbf{\xi}; k, \varkappa) & = 
	\sum_{\tmmathbf{m} \in \mathbb{Z}^2} u (\tmmathbf{r}+ \Lambda \tmmathbf{m};
	\tmmathbf{r}_s ; k, \varkappa) e^{i\tmmathbf{m} \cdot \tmmathbf{\xi}}, 
	\label{eq:BF-transform}
\end{align}
and sometimes referred to as the resolvent.  
%As for the Fourier transform,
%there are various conventions for the multiplicative constant to be used in
%the forward (\ref{eq:BF-transform}) and inverse
%(\ref{eq:inverse-BF-transform}) Floquet-Bloch transforms, as well as for the
%sign inside the exponential. Here we choose them so that it matches with our
%choices for the Fourier Transform in {\cite{Part6A}}.
The variable $\tmmathbf{\xi}$ is sometimes called the Bloch wave vector, the
Bloch quasi-momentum or the Floquet exponents. In what follows, we will think
of it as a pair of complex variables, that is we consider that $\tmmathbf{\xi}
\in \mathbb{C}^2$, and we will assume that the transform $F (\tmmathbf{r};
\tmmathbf{r}_s ; \tmmathbf{\xi}; k, \varkappa)$ is holomorphic as a function
of $\tmmathbf{\xi} \in \mathbb{C}^2$ away from its singularity set $H (k,
\varkappa)$. The singularity set is known as the dispersion diagram of the media, which is usually interpreted as the set of all eigenwaves  of the media. Indeed, each point $(\tmmathbf{\xi}, k, \varkappa)$ that belongs to $H(k,\varkappa)$  corresponds to an eigenwave $V(\tmmathbf{r} ; \tmmathbf{\xi}; k, \varkappa)$  obeying the Floquet periodicity conditions:
\begin{align*}
	V(\tmmathbf{r}+\tmmathbf{d} ; \tmmathbf{\xi}; k, \varkappa) &= e^{-i\xi_1}V(\tmmathbf{r} ; \tmmathbf{\xi}; k, \varkappa)\\
	V(\tmmathbf{r}+\tmmathbf{l}; \tmmathbf{\xi}; k, \varkappa) &= e^{-i\xi_2}V(\tmmathbf{r}; \tmmathbf{\xi}; k, \varkappa).
\end{align*}
We will further assume that $F (\tmmathbf{r};
\tmmathbf{r}_s ; \tmmathbf{\xi}; k, \varkappa)$ has the so-called real
property that is that:
\begin{itemize}
	\item \tmcolor{blue}{P1.} If $\varkappa > 0$, then $H (k, \varkappa) \cap
	\tmmathbf{\mathcal{B}}= \emptyset$, i.e. $F$ is singularity-free on
	$\tmmathbf{\mathcal{B}}$, and hence the integral
	(\ref{eq:inverse-BF-transform}) is well-defined in this case. The
	singularity set $H (k, \varkappa)$ is a union of irreducible components
	$\sigma_j (k, \varkappa)$ for $j = 1, 2, \ldots$, each described as the zero
	set of a holomorphic function of $\tmmathbf{\xi}$, depending on the parameters
	$k$ and $\varkappa$, denoted $g_j (\tmmathbf{\xi}; k, \varkappa)$.
	
	\item \tmcolor{blue}{P2.} If $\varkappa = 0$, then the singularities of $F$
	hit $\tmmathbf{\mathcal{B}}$, and $H' (k) \equiv H (k, 0) \cap
	\tmmathbf{\mathcal{B}}$ is a union of smooth one dimensional curves
	$\sigma_1' (k), \ldots, \sigma_N' (k)$ that are called the real traces of
	the irreducible singular components $\sigma_1 (k, 0), \ldots, \sigma_N (k,
	0)$ of $H (k, 0)$. It also means that $g_j(\tmmathbf{\xi}; k, 0)$ is real whenever its arguments are real. 
	\RED{In the context of periodic structures, these real traces are also sometimes called equi- or iso-frequency contours, or slowness curves.}
	%The irreducible components are sometimes referred to as
	%the {\tmem{slowness surfaces}} of the problem.
	
	\item \tmcolor{blue}{P3.} The real traces $\sigma_j'$ are {\tmem{regular}}
	in the sense that if $\tmmathbf{\xi}^{\star} \in \sigma_j'$, then
	$\nabla_{\!\! \tmmathbf{\xi}} g_j (\tmmathbf{\xi}^{\star} ; k, 0)
	\neq \tmmathbf{0}$ and $\frac{\partial g_j}{\partial \varkappa}
	(\tmmathbf{\xi}^{\star}, k, 0) \neq 0$.
\end{itemize}
Because of point \tmcolor{blue}{P2}, the integral
(\ref{eq:inverse-BF-transform}) is not well-defined when $\varkappa = 0$, and
it needs to be rewritten as
\begin{align}
	u (\tmmathbf{r}+ \Lambda \tmmathbf{m}; \tmmathbf{r}_s ; k, 0) & = 
	\frac{1}{4 \pi^2} \iint_{\tmmathbf{\Gamma}} F (\tmmathbf{r}; \tmmathbf{r}_s
	; \tmmathbf{\xi}; k, 0) e^{- i\tmmathbf{m} \cdot \tmmathbf{\xi}} \mathd
	\tmmathbf{\xi},  \label{eq:inverse-BF-transform-deformed}
\end{align}
where $\tmmathbf{\Gamma}$ is an admissible deformation/indentation of
$\tmmathbf{\mathcal{B}}$ into $\mathbb{C}^2$, in the sense that $H (k, 0) \cap
\tmmathbf{\Gamma}= \emptyset$ and that one can deform $\tmmathbf{\mathcal{B}}$
to $\tmmathbf{\Gamma}$ smoothly without hitting any singularities of $F$ for
any small enough $\varkappa$. Due to the two-dimensional version of Cauchy's
integral theorem {\cite{Shabat2}}, such deformation does not change the value
of the integral. However, because $\tmmathbf{\Gamma}$ is not real anymore,
$\mathd \tmmathbf{\xi}$ should be thought of in terms of differential forms in
(\ref{eq:inverse-BF-transform-deformed}), and we write $\mathd \tmmathbf{\xi}=
\mathd \xi_1 \wedge \mathd \xi_2$. The deformation process
$\tmmathbf{\mathcal{B}} \rightarrow \tmmathbf{\Gamma}$ will be explained in
more details in section \ref{sec:deformation-process}.

The aim of the present work is to study the far-field behaviour of $u
(\tmmathbf{r}+ \Lambda \tmmathbf{m}; \tmmathbf{r}_s ; k, 0)$ as $|
\tmmathbf{m} | \rightarrow \infty$. To do this, we will fix an
{\tmem{observation direction}} $\tilde{\tmmathbf{m}}\RED{ \in\mathbb{Z}^2}$ and study $u$ for
$\tmmathbf{m}= N \tilde{\tmmathbf{m}}$ as the integer $N$ tends to infinity.
For each chosen observation $\tilde{\tmmathbf{m}}$ we will show that
\begin{eqnarray*}
	u (\tmmathbf{r}+ \Lambda \tmmathbf{m}; \tmmathbf{r}_s ; k, 0) & \underset{N
		\rightarrow \infty}{\approx} & \sum_q u_{\tmop{loc}}^{(q)} (\tmmathbf{m};
	\tmmathbf{r}; \tmmathbf{r}_s ; k),
\end{eqnarray*}
where the far-field wave components $u_{\tmop{loc}}^{(q)}$ can be written down
explicitly and results from local consideration around a set of {\tmem{special
		points}} $\tmmathbf{\xi}_q$ belonging to $\tmmathbf{\mathcal{B}}$. This is
similar to what we found in {\cite{Part6A}} for double Fourier transforms in
the context of diffraction theory and called the {\tmem{locality principle}}.
The importance of the real traces of the integrand of multiple integrals in diffraction problems was previously highlighted in \cite{chapman-wavenumber-surface} in the context of blade-gust interactions in which the real traces are called the \textit{wavenumber surfaces}.

Additional results will be provided for values of $k$ associated with
{\tmem{degeneracies}} corresponding to the shrinking of a real trace to a
point. But before \RED{delving} any further into this, we wish to show that the
analyticity of $F$ and its real property are not unreasonable assumptions, by
considering two illustrating examples: a discrete lattice (Figure
\ref{fig:disc-lattice}, left), and a simple phononic crystal (Figure
\ref{fig:decomposition-of-boundary}, left). A common understanding is that all
singularities of $F$ are poles if there is no energy leakage from the medium
to some other subsystem. This will be the case in the examples considered
below. However, one can expect singularities of the branching type, say, if
the 2D medium is a thin layer placed atop some acoustic half-space. The method
developed in {\cite{Part6A}} and in what follows is able to deal with such
non-polar singularities.

\section{Two illustrating examples}\label{sec:examples-illustrative}

\subsection{Green's function of a discrete lattice}\label{sec:discrete-example}

Let us start by studying the discrete Helmholtz Green's function of a square
lattice. Without loss of generality, we can assume that the lattice points are
lying on the integer points $\tmmathbf{m}= (m_1, m_2) \in \mathbb{Z}^2$ of a
Cartesian plane and that the point source is at the origin. The matrix
$\Lambda$ is simply the identity matrix, $\tmmathbf{d} \equiv \tmmathbf{e}_1 =
(1, 0)^{\transpose}$ and $\tmmathbf{\ell} \equiv \tmmathbf{e}_2 = (0, 1)^{\transpose}$. Because of
the discrete nature of this media, the periodic cells are just points and
hence there are no explicit dependency on $\tmmathbf{r}$ and $\tmmathbf{r}_s$.
For these reasons, the quantity $u (\tmmathbf{r}+ \Lambda \tmmathbf{m};
\tmmathbf{r}_s ; k, \varkappa)$ of the general framework can just be written
$u (\tmmathbf{m}; k, \varkappa)$. It satisfies the discrete Helmholtz
equation, forced at the origin, given by
\begin{align}
	&u (\tmmathbf{m}+\tmmathbf{e}_1 ; k, \varkappa) + u
	(\tmmathbf{m}-\tmmathbf{e}_1 ; k, \varkappa) + u
	(\tmmathbf{m}+\tmmathbf{e}_2 ; k, \varkappa)\nonumber \\
	&+ u (\tmmathbf{m}-\tmmathbf{e}_2 ; k, \varkappa) + ((k + i \varkappa)^2 - 4)
	u (\tmmathbf{m}; k, \varkappa) =  \delta_{m_1 0} \delta_{m_2 0}, 
	\label{eq:discrete-governing}
\end{align}
where $u (\tmmathbf{m}; k, \varkappa)$ represents the out-of-plane displacement of
the node $\tmmathbf{m}= (m_1, m_2)$ for instance, $\delta$ is the Kronecker
delta, $k + i \varkappa$ is the wavenumber, with real part $k \geqslant 0$ and
a small imaginary part $\varkappa \geqslant 0$.

For $\tmmathbf{\xi}= (\xi_1, \xi_2)$, we can formally introduce the double
Floquet-Bloch transform %which is just a discrete Fourier transform
\begin{eqnarray*}
	F (\tmmathbf{\xi}; k, \varkappa) & = & \sum_{\tmmathbf{m} \in \mathbb{Z}^2}
	u (\tmmathbf{m}; k, \varkappa) e^{i\tmmathbf{m} \cdot \tmmathbf{\xi}},
\end{eqnarray*}
where, as for $u$, we have removed the dependency on $\tmmathbf{r}$ and
$\tmmathbf{r}_s$ from $F$. Applying this transform to the governing equation
(\ref{eq:discrete-governing}), we obtain an explicit expression for $F$ given
by
\begin{align}
	F (\tmmathbf{\xi}; k, \varkappa) & =  \frac{1}{2 \cos (\xi_1) + 2 \cos
		(\xi_2) + ((k + i \varkappa)^2 - 4)} \cdot 
	\label{eq:discrete-explicit-expr}
\end{align}
It is clear that, as a function of $\tmmathbf{\xi} \in \mathbb{C}^2$, $F
(\tmmathbf{\xi}; k, \varkappa)$ is holomorphic everywhere away from its
singularity set $H (k, \varkappa)$, which happens to be the polar set defined by
\begin{eqnarray*}
	H (k, \varkappa) & = & \{ \tmmathbf{\xi} \in \mathbb{C}^2, 2 \cos (\xi_1) +
	2 \cos (\xi_2) + ((k + i \varkappa)^2 - 4) = 0 \} .
\end{eqnarray*}
In fact here, for each $(k, \varkappa)$, there is only one irreducible
component $\sigma_1 (k, \varkappa)$ with defining function $g_1
(\tmmathbf{\xi}; k, \varkappa) = 2 \cos (\xi_1) + 2 \cos (\xi_2) + ((k + i
\varkappa)^2 - 4)$.

\begin{figure}[h]
	\centering{
		\includegraphics[height=0.25\textwidth]{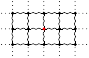} \qquad \includegraphics[height=0.25\textwidth]{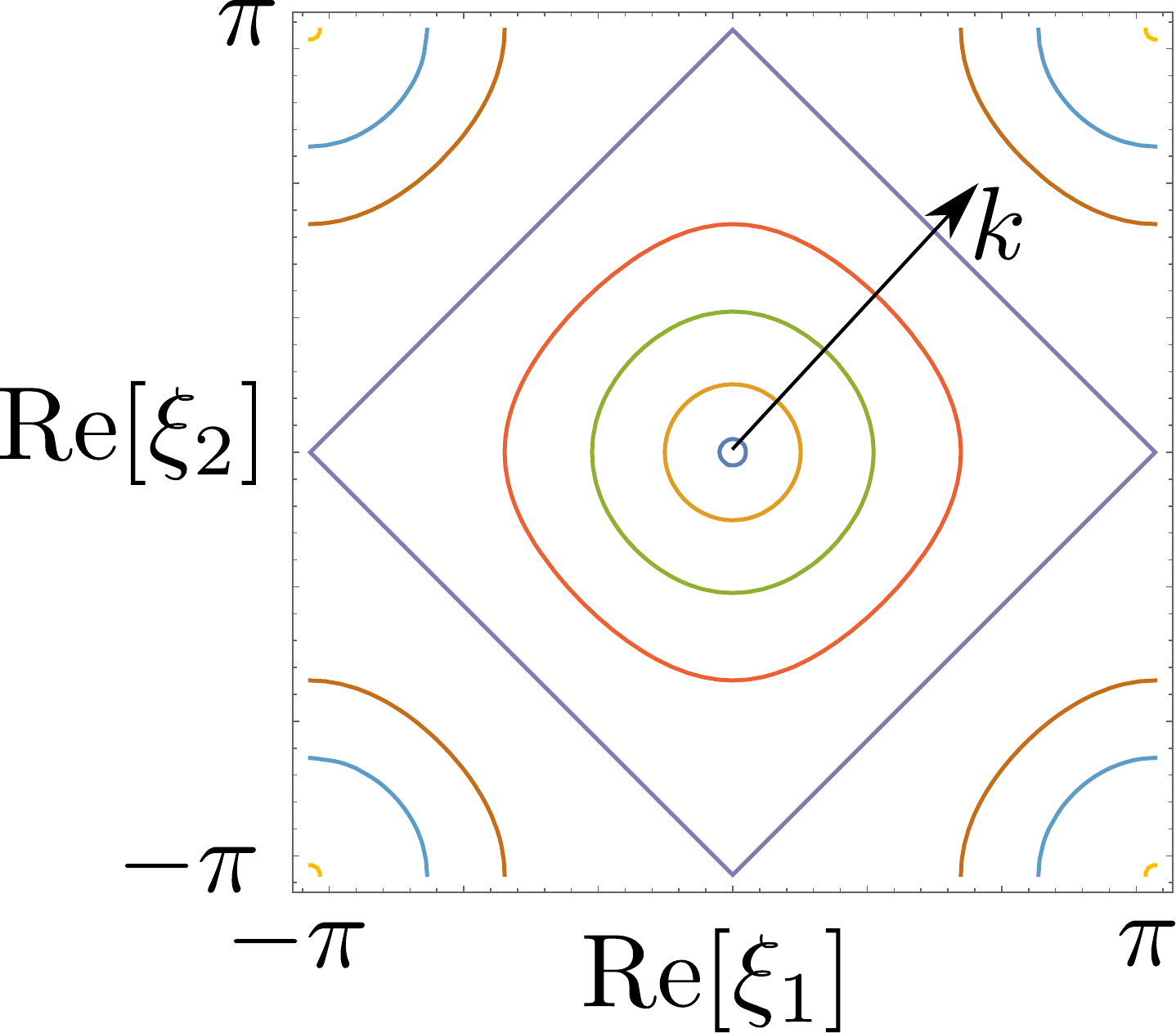}
	}
	\caption{A discrete square lattice (left) and its $H' (k)$ for different
		values of $k$, the arrow shows an increasing value of $k$ (right)}
	\label{fig:disc-lattice}
\end{figure}

It is clear that if $\varkappa > 0$, then $H (k, \varkappa) \cap
\tmmathbf{\mathcal{B}}= \emptyset$ (point \tmcolor{blue}{P1}) so that in this
case, $u$ is simply given by
\begin{eqnarray}
	u (\tmmathbf{m}, k, \varkappa) & = & \frac{1}{4 \pi^2}
	{\iint_{\tmmathbf{\mathcal{B}}}}  F (\tmmathbf{\xi}; k, \varkappa) e^{-
		i\tmmathbf{m} \cdot \tmmathbf{\xi}} \mathd \tmmathbf{\xi}. 
	\label{eq:integral-discrete}
\end{eqnarray}
We now want to study the limit $\varkappa \rightarrow 0$. For $k > 2
\sqrt{2}$, $H' (k) \equiv H (k, 0) \cap \tmmathbf{\mathcal{B}}= \emptyset$,
and the integral (\ref{eq:integral-discrete}) remains valid. However, if $0 <
k < 2$ or $2 < k < 2 \sqrt{2}$, the singularities of $F (\tmmathbf{\xi}; k,
0)$ intersect the real plane as a set of smooth one-dimensional curves
$\sigma_1' (k)$ (point \tmcolor{blue}{P2} of the real property), as illustrated in Figure
\ref{fig:disc-lattice} (right), and one can check that they are regular by
direct computation (point \tmcolor{blue}{P3} of the real property). So $F$ has the real property
for those $k$. In this case, for the integral (\ref{eq:integral-discrete}) to
make sense, the surface $\tmmathbf{\mathcal{B}}$ should be indented to a
surface $\tmmathbf{\Gamma}$ that is close to $\tmmathbf{\mathcal{B}}$
everywhere but bypasses the singularities to give
\begin{eqnarray}
	u (\tmmathbf{m}; k, 0) & = & \frac{1}{4 \pi^2} \iint_{\tmmathbf{\Gamma}} F
	(\tmmathbf{\xi}; k, 0) e^{- i\tmmathbf{m} \cdot \tmmathbf{\xi}} \mathd
	\tmmathbf{\xi},  \label{eq:integral-indented-discrete}
\end{eqnarray}
and we are hence exactly within the general framework described in section
\ref{sec:general-framework}.

Note further that, when $k = 0$ or $k = 2 \sqrt{2}$, the real traces shrink
into points, while as $k \rightarrow 2$ the curves become straight
intersecting lines. At these specific values of $k$, we lose the point
\tmcolor{blue}{P3} of the real property. A careful study of these degeneracies will be provided in
section \ref{sec:vanishing-loops}.

\subsection{Simple phononic
	crystal}\label{sec:continuous-example}

\subsubsection{Formulation}

Let us consider a 2D continuous periodic medium with double period given by
the vectors $\tmmathbf{d}$ and $\tmmathbf{\ell}$ lying in a Cartesian plane
$\tmmathbf{r}= (x_1, x_2) \in \mathbb{R}^2$. The periodically repeated
parallelogram cells $S$ all contain an identical scatterer $\Omega$ say subject
to some boundary condition (BC), as illustrated in Figure
\ref{fig:decomposition-of-boundary} (left). The discussion below is general in
nature and holds for many boundary conditions: Dirichlet ($u = 0$), Neumann
($\tmmathbf{n} \cdot \nabla u = 0$) or real positive Robin ($\tmmathbf{n}
\cdot \nabla u + \beta u = 0$, with $\beta > 0$). With very minor
modifications, it can also be applied for penetrable scatterers or
continuously varying material properties. We are interested in finding the
tailored Green's function of this medium. Let us hence assume without loss of
generality that a point source $\tmmathbf{r}_s$ is located away from $\Omega$
within the cell $S_0$ that contains the origin, so $\tmmathbf{r}_s \in S_0
\setminus \Omega$. We are interested in solving the Helmholtz equation
\begin{align}
	\Delta_{\tmmathbf{r}} u (\tmmathbf{r}; \tmmathbf{r}_s ; k, \varkappa) + (k +
	i \varkappa)^2 u (\tmmathbf{r}; \tmmathbf{r}_s ; k, \varkappa) &= \delta
	(\tmmathbf{r}-\tmmathbf{r}_s),\label{eq:continuous-governing-eq} %& \tmop{with} & \tmop{BC} \text{on each
	%	$\partial \Omega$},  \label{eq:continuous-governing-eq}
\end{align}
with BC on each $\partial \Omega$,
where $k + i \varkappa$ is the wavenumber, with positive real part $k$ and a
small positive imaginary part $\varkappa$. Let $\tmmathbf{r} \in S_0 \setminus
\Omega$ and formally introduce the double Floquet-Bloch transform
\begin{eqnarray*}
F (\tmmathbf{r}; \tmmathbf{r}_s ; \tmmathbf{\xi}; k, \varkappa) & = &
\sum_{\tmmathbf{m} \in \mathbb{Z}^2} u (\tmmathbf{r}+ \Lambda \tmmathbf{m};
\tmmathbf{r}_s ; k, \kappa) e^{i\tmmathbf{m} \cdot \tmmathbf{\xi}},
\end{eqnarray*}
By applying the Floquet-Bloch transform to the governing equation
(\ref{eq:continuous-governing-eq}), and noting that differential operators in
the $\tmmathbf{r}$ variable commute with the transform, it can be seen that,
for $\tmmathbf{r} \in S_0 \setminus \Omega$,
\begin{align}
\Delta_{\tmmathbf{r}} F (\tmmathbf{r}; \tmmathbf{r}_s ; \tmmathbf{\xi}; k,
\varkappa) + (k + i \varkappa)^2 \nobracket F (\tmmathbf{r}; \tmmathbf{r}_s
; \tmmathbf{\xi}; k, \varkappa) & =  \delta (\tmmathbf{r}-\tmmathbf{r}_s), 
\label{eq:equationpourL}
\end{align}
and that $F$ also satisfies the same BC on $\partial \Omega$ as $u$ does.
Moreover, direct calculations show that it has the Floquet periodicity
property (also called $\tmmathbf{\xi}$-quasiperiodicity w.r.t.\ $\Lambda$):
\begin{align}
\forall \tmmathbf{j} \in \mathbb{Z}^2, \quad F (\tmmathbf{r}+ \Lambda
\tmmathbf{j}; \tmmathbf{r}_s ; \tmmathbf{\xi}; k, \varkappa) & =  e^{-
	i\tmmathbf{\xi} \cdot \tmmathbf{j}} F (\tmmathbf{r}; \tmmathbf{r}_s ;
\tmmathbf{\xi}; k, \varkappa) .  \label{eq:Floquet-periodicity}
\end{align}
Note that (\ref{eq:Floquet-periodicity}) implies its gradient counterpart:
\begin{align}
\forall \tmmathbf{j} \in \mathbb{Z}^2, \quad \nabla_{\!\!
	\tmmathbf{r}} F (\tmmathbf{r}+ \Lambda \tmmathbf{j}; \tmmathbf{r}_s ;
\tmmathbf{\xi}; k, \varkappa) & =  e^{- i\tmmathbf{\xi} \cdot \tmmathbf{j}}
\nabla_{\!\! \tmmathbf{r}} F (\tmmathbf{r}; \tmmathbf{r}_s ;
\tmmathbf{\xi}; k, \varkappa) .  \label{eq:Floquet-grad-periodicity}
\end{align}

Hence, to find $F$, one needs to solve (\ref{eq:equationpourL}) on the unit
cell $S_0$ subject to BC on $\partial \Omega$ and Floquet periodicity
conditions (\ref{eq:Floquet-periodicity})-(\ref{eq:Floquet-grad-periodicity}).
For these reasons, $F$ is often referred to as the tailored quasi-periodic
Green's function of the problem associated with $\tmmathbf{\xi}$.

In what follows, we want to show that $F$ can be understood as a holomorphic
function of $\tmmathbf{\xi}$ with the real property. To get there, we first
need to study $F$ when $\tmmathbf{\xi}$ is real.

\subsubsection{Self-adjointness and an explicit formula for $F$ when $\bdxi$ is real}
\label{sec:self-adjoint}

\begin{figure}[h]
\centering{
	\includegraphics[height=0.22\textwidth]{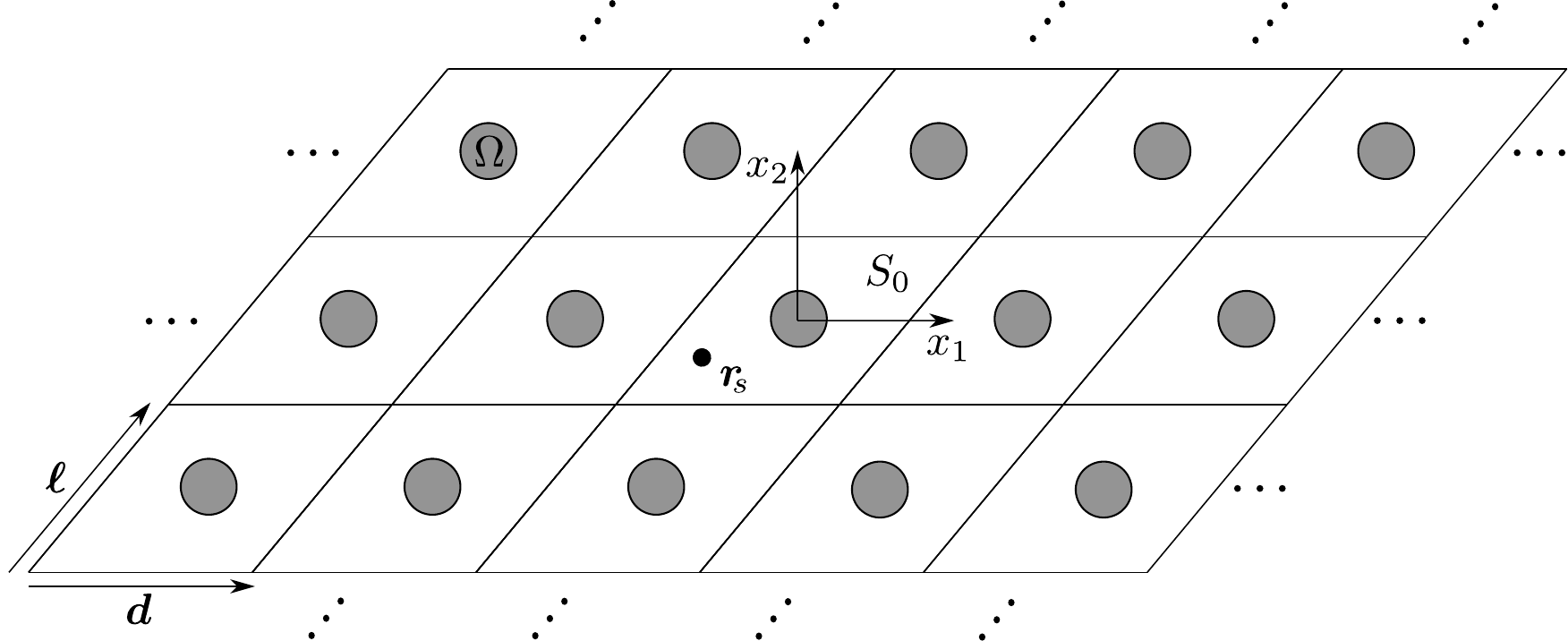}\includegraphics[height=0.22\textwidth]{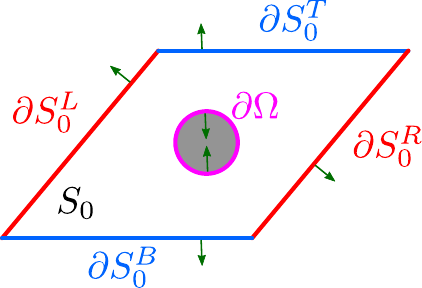}
}
\caption{A phononic crystal (left) and decomposition of the boundary
	$\partial (S_0 \setminus \Omega)$ with illustration of the orientation of
	the normal $\tmmathbf{n}$ (right).}
\label{fig:decomposition-of-boundary}
\end{figure}

Let us consider $\tmmathbf{\xi} \in \tmmathbf{\mathcal{B}} \subset
\mathbb{R}^2$ and consider the operator\footnote{In this section there are no
ambiguity about which variable the Laplacian is defined with respect to, so
there is no need for the subscript $_{\tmmathbf{r}}$.} $\mathcal{L}= - \Delta$
with a domain $\tmop{dom}_{\tmmathbf{\xi}}$ consisting of functions defined on
$S_0 \setminus \Omega$, that satisfy BC on $\partial \Omega$ and are
$\tmmathbf{\xi}$-quasiperiodic w.r.t.\ $\Lambda$. Let us consider the usual
$L_2$ inner product defined for two functions $f (\tmmathbf{r})$ and $g
(\tmmathbf{r})$ in $\tmop{dom}_{\tmmathbf{\xi}}$ by
\begin{eqnarray*}
(f, g) & = & \iint_{S_0 \setminus \Omega} fg^{\conj} \mathd S (\tmmathbf{r}),
\end{eqnarray*}
where $\conj$ means complex conjugate. By Green's identity, we have
\begin{align}
(\mathcal{L}f, g) & = \iint_{S_0 \setminus \Omega}
\nabla f \cdot \nabla g^{\conj} \mathd S
- \int_{\partial (S_0 \setminus \Omega)} (\tmmathbf{n} \cdot \nabla f)
g^{\conj} \mathd \ell,  \label{eq:operator-symmetry-part1}
\end{align}
where $\partial (S_0 \setminus \Omega)$ can be decomposed into $\partial (S_0
\setminus \Omega) = \partial S^T_0 \cup \partial S^B_0 \cup \partial S^R_0
\cup \partial S^L_0 \cup \partial \Omega$ and the normal is outgoing as
described in Figure \ref{fig:decomposition-of-boundary} (right).

Due to the $\tmmathbf{\xi}$-quasiperiodicity, and the fact that $(e^{- i
\xi_{1, 2}})^{\conj} = e^{+ i \xi_{1, 2}}$ when $\tmmathbf{\xi}$ is real, it is
direct to show that
\begin{eqnarray*}
\int_{\partial S^R \cup \partial S^L} (\tmmathbf{n} \cdot \nabla f) g^{\conj}
\mathd \ell = 0 & \text{ and } & \int_{\partial S^T \cup \partial S^B}
(\tmmathbf{n} \cdot \nabla f) g^{\conj} \mathd \ell = 0,
\end{eqnarray*}
leading to the second term in the RHS of (\ref{eq:operator-symmetry-part1}) to
simplify to $- \int_{\partial \Omega} (\tmmathbf{n} \cdot \nabla f) g^{\conj}
\mathd \ell$. Following the same approach, we can show that for $f, g \in
\tmop{dom}_{\tmmathbf{\xi}}$, we get
\begin{align}
(f, \mathcal{L}g) & =  \iint_{S_0 \setminus \Omega} \nabla f \cdot \nabla
g^{\conj} \mathd S - \int_{\partial \Omega} f (\tmmathbf{n} \cdot \nabla
g^{\conj}) \mathd \ell .  \label{eq:operator-symmetry-part2-raph}
\end{align}
Using the BC, one can show directly that $\int_{\partial \Omega} (\tmmathbf{n}
\cdot \nabla f) g^{\conj} \mathd \ell = \int_{\partial \Omega} f (\tmmathbf{n}
\cdot \nabla g^{\conj}) \mathd \ell$ (for Dirichlet and Neumann these terms are
just zero, and equal to $\beta \int_{\partial \Omega} f g^{\conj} \mathd \ell$
for Robin for some $\beta > 0$). Therefore, by comparing
(\ref{eq:operator-symmetry-part1}) and (\ref{eq:operator-symmetry-part2-raph})
we find that $(\mathcal{L}f, g) = (f, \mathcal{L}g)$ for all $f, g \in
\tmop{dom}_{\tmmathbf{\xi}}$, that is, the operator $(\mathcal{L},
\tmop{dom}_{\tmmathbf{\xi}})$ is symmetric. As a result, one automatically
gets that its eigenvalues are real. Moreover, using $g = f$ in
(\ref{eq:operator-symmetry-part1}), we find that for all $f \in
\tmop{dom}_{\tmmathbf{\xi}}$, $(\mathcal{L}f, f) \geqslant 0$ and that it is
only zero for $f = 0$. In other words, the operator $(\mathcal{L},
\tmop{dom}_{\tmmathbf{\xi}})$ is non-negative. These two properties (symmetry
and non-negativity) imply that $(\mathcal{L}, \tmop{dom}_{\tmmathbf{\xi}})$ is
self-adjoint and non-negative. So far we have not been very precise about the
regularity of the functions in $\tmop{dom}_{\tmmathbf{\xi}}$. We can follow
two approaches to show the self-adjoint property. Either we consider smooth
functions to start with, and then show that $(\mathcal{L},
\tmop{dom}_{\tmmathbf{\xi}})$ has a non-negative self-adjoint (Friedrichs)
extension (see, e.g.\ {\cite{Kato1976,Davies2007,Assier2016}}). Equivalently,
we can bypass this step and consider functions to be in an appropriate Sobolev
space from the start and in this case $(\mathcal{L},
\tmop{dom}_{\tmmathbf{\xi}})$ is directly self-adjoint.

As a consequence, the spectral theorem implies that the eigenvalues of this
operator are discrete, real, and positive and can hence be written as
\begin{eqnarray}
& 0 \leqslant \lambda_0 \leqslant \lambda_1 \leqslant \lambda_2 \leqslant
\cdots \leqslant \lambda_j \leqslant \cdots & 
\label{eq:eigenvalues-inequalities}
\end{eqnarray}
with $\lambda_j \rightarrow + \infty$ as $j \rightarrow \infty$. We are using
$\leqslant$ in (\ref{eq:eigenvalues-inequalities}) to account for multiple
eigenvalues. The eigenfunctions $V_j$ associated to the eigenvalues
$\lambda_j$ span the entire function space $\tmop{dom}_{\tmmathbf{\xi}}$ and
are orthogonal, that is $(V_j, V_{\ell}) = 0$ if $j \neq \ell$. Because
everything here depends on $\tmmathbf{\xi}$, we will from now on write
\begin{eqnarray*}
\lambda_j \equiv \lambda_j (\tmmathbf{\xi}) & \tmop{and} & V_j \equiv V_j
(\tmmathbf{r}; \tmmathbf{\xi}) .
\end{eqnarray*}
Note now that $F (\tmmathbf{r}; \tmmathbf{r}_s ; \tmmathbf{\xi}; k,
\varkappa)$ is in $\tmop{dom}_{\tmmathbf{\xi}}$, a domain spanned by the
eigenfunctions of $(\mathcal{L}, \tmop{dom}_{\tmmathbf{\xi}})$. Therefore, $F$
can be written as
\begin{eqnarray*}
F (\tmmathbf{r}; \tmmathbf{r}_s ; \tmmathbf{\xi}; k, \varkappa) & = &
\sum_{j = 0}^{\infty} A_j (\tmmathbf{r}_s ; k, \varkappa) V_j (\tmmathbf{r};
\tmmathbf{\xi}),
\end{eqnarray*}
for some coefficients $A_j$ to be determined. Using the orthogonality of the
eigenfunctions, the symmetry of the operator $(\mathcal{L},
\tmop{dom}_{\tmmathbf{\xi}})$, the fact that the eigenvalues are real, the
equation (\ref{eq:equationpourL}) and the properties of the Dirac delta
function, we find that
\begin{eqnarray*}
A_j (\tmmathbf{r}_s ; k, \varkappa) & = & \frac{1}{((k + i \varkappa)^2 -
	\lambda_j (\tmmathbf{\xi}))} \frac{V_j^{\conj} (\tmmathbf{r}_s ;
	\tmmathbf{\xi})}{(V_j (\cdot, \tmmathbf{\xi}), V_j (\cdot, \tmmathbf{\xi}))}
\cdot
\end{eqnarray*}
This leads to an explicit expression for $F$:
%\begin{widetext}
\begin{align}
	F (\tmmathbf{r}; \tmmathbf{r}_s ; \tmmathbf{\xi}; k, \varkappa) & = 
	\sum_{j = 0}^{\infty} \frac{1}{((k + i \varkappa)^2 - \lambda_j
		(\tmmathbf{\xi}))} \times \frac{V_j^{\conj} (\tmmathbf{r}_s ;
		\tmmathbf{\xi})}{\| V_j^{\conj} (. ; \tmmathbf{\xi}) \|} \times \frac{V_j
		(\tmmathbf{r}; \tmmathbf{\xi})}{\| V_j (. ; \tmmathbf{\xi}) \|}, 
	\label{eq:continuous-explicit-F}
\end{align}
where $\| . \| = \sqrt{(., .)}$ and we have used the fact that $\| V_j^{\conj}
\| = \| V_j \|$. Given the form of (\ref{eq:continuous-explicit-F}), it is
clear that, as expected, $F$ does not depend on the choice of normalisation of
the eigenfunctions. Using the fact that if $(\lambda_j (\tmmathbf{\xi}), V_j
(\tmmathbf{r}; \tmmathbf{\xi}))$ is an eigenpair of $(\mathcal{L},
\tmop{dom}_{\tmmathbf{\xi}})$ then $(\lambda_j (\tmmathbf{\xi}), V_j^{\conj}
(\tmmathbf{r}; \tmmathbf{\xi}))$ is an eigenpair of $(\mathcal{L},
\tmop{dom}_{-\tmmathbf{\xi}})$, the formula (\ref{eq:continuous-explicit-F})
can be rewritten as
\begin{align}
	F (\tmmathbf{r}; \tmmathbf{r}_s ; \tmmathbf{\xi}; k, \varkappa) & = 
	\sum_{j = 0}^{\infty} \frac{1}{((k + i \varkappa)^2 - \lambda_j
		(\tmmathbf{\xi}))} \times \frac{V_j (\tmmathbf{r}_s ; -\tmmathbf{\xi}) V_j
		(\tmmathbf{r}; \tmmathbf{\xi})}{\langle V_j (\cdot ; -\tmmathbf{\xi}), V_j
		(. ; \tmmathbf{\xi}) \rangle}, 
	\label{eq:continuous-explicit-F-no-conjugate}
\end{align}
%\end{widetext}
where the bracket $\langle \cdot, \cdot \rangle$ is defined for two functions
$f$ and $g$ as
\begin{align*}
\langle f, g \rangle & =  \iint_{S_0 \setminus \Omega} f (\tmmathbf{r}) g
(\tmmathbf{r}) \mathd S.
\end{align*}
The advantage of (\ref{eq:continuous-explicit-F-no-conjugate}) over
(\ref{eq:continuous-explicit-F}) is that it does not contain any complex
conjugate symbol, which will prove helpful when we will attempt to
analytically continue this formula for complex $\tmmathbf{\xi}$.

\begin{rema} \label{rem:discrete-case-contiuous}
The discrete example can also be written as a special case of
(\ref{eq:continuous-explicit-F}). In that case, there is only one
eigenvalue for the eigenvalue problem $- \tilde{\Delta} f = \lambda f$,
where $\tilde{\Delta}$ is the discrete Laplacian defined by $\tilde{\Delta}
f = f (\tmmathbf{m}+\tmmathbf{e}_1) + f (\tmmathbf{m}-\tmmathbf{e}_1) + f
(\tmmathbf{m}+\tmmathbf{e}_2) - f (\tmmathbf{m}-\tmmathbf{e}_2) - 4 f
(\tmmathbf{m})$, subject to Floquet quasi-periodic conditions $f
(\tmmathbf{m}+\tmmathbf{j}) = e^{- i\tmmathbf{\xi} \cdot \tmmathbf{j}} f
(\tmmathbf{m})$, with the inner product $(f, g) = f (\tmmathbf{0}) g^{\conj}
(\tmmathbf{0})$. This problem has only one eigenvalue $\lambda_0
(\tmmathbf{\xi}) = 4 - 2 \cos (\xi_1) - 2 \cos (\xi_2)$, the eigenfunction
is given by any constant $V_0 (\tmmathbf{0}; \tmmathbf{\xi}) \neq 0$.
Moreover we have $\tmmathbf{r}=\tmmathbf{r}_s =\tmmathbf{0}$ and, given the
definition of the inner product, we have that $V_0^{\conj} (\tmmathbf{0};
\tmmathbf{\xi}) V_0 (\tmmathbf{0}; \tmmathbf{\xi}) / (\| V_0^{\conj} \| \|
V_0 \|) = 1$. Hence, (\ref{eq:discrete-explicit-expr}) can indeed be thought
of as a special case of (\ref{eq:continuous-explicit-F}).
\end{rema}

\begin{rema}
Since the double integral (\ref{eq:inverse-BF-transform-deformed})  with (\ref{eq:continuous-explicit-F-no-conjugate}) has only
polar singularities, one can convert it into a
single integral by taking the integral either with respect to $\xi_1$ or
$\xi_2$. For example, let us assume that $m_2 < 0$. Fix some value of $\xi_1 \in [- \pi,
\pi]$. The contour of integration can be transformed from the segment $[- \pi,
\pi]$ into the contour shown in Figure
\ref{fig07}. The added passes along the vertical half-lines compensate each
other due to periodicity. The integrand is exponentially small at infinity, so
one can close the contour there.

\begin{figure}[h]
	\centering{
		\includegraphics[width=0.3\textwidth]{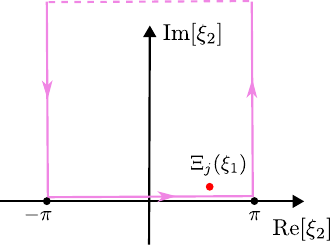}}
	\caption{Deformation of the integration contour with respect to $\xi_2$}
	\label{fig07}
\end{figure}
Finally, the integral with respect to $\xi_2$ becomes expressed a sum of poles
located inside the contour (i.e. in the upper half-plane when $\varkappa>0$). Denote the position of these poles by $\Xi_j
(\xi_1)$. The expression for the integral is as follows		
\begin{align}
	&u (\tmmathbf{r};\tmmathbf{r}_s;k,\kappa) = \label{e:02002} \\ &- \frac{i}{2\pi}
	\int_{- \pi}^{\pi} d \xi_1  \sum_{j=0}^{\infty} 
	\frac{e^{-i \tmmathbf{m}\cdot\tmmathbf{\xi}_j}}{\partial \lambda_j / \partial
		\xi_2}  \times \frac{V_{j} (\tmmathbf{r}_s ; \tmmathbf{\xi}_j)V_{j} (\tmmathbf{r} ; \tmmathbf{\xi}_j)}{\langle 			V_{j}(\cdummy ;  -\tmmathbf{\xi}_j), V_{j} (\cdummy ;  \tmmathbf{\xi}_j)\rangle}, \nonumber	
\end{align}
where $\tmmathbf{\xi}_j \equiv (\xi_1,\Xi_j (\xi_1))$, and the derivative $\partial\lambda_j / \partial
\xi_2$ is evaluated at $\tmmathbf{\xi}_j $.	
%	{\color{red} \textbf{Note for AIK}: I think the sign in the exponential was wrong in the eq above, I've also added missing indices on the $V$. Also following our discussion, I removed the discussion around $-\Xi_j$ and adapted the figure accordingly.} {\color{green} I agree}
The formula (\ref{e:02002}) is a series-integral representation of the field.
Instead of taking the integral with respect to $\xi_2$, one can use $\xi_1$
and get an alternative series-integral representation. In {\cite{Shanin2022}}
we consider a pair of representations in a similar situation and demonstrate
an equivalence between them.

Note that the representation (\ref{e:02002}) is an integral of a differential
form defined on a complex manifold (possibly with some points removed).
The integral is taken along some contours on this manifold. A special technique
{\cite{Shanin2022a,Kniazeva2022}} can be applied to estimate such an integral.
Below, however, we prefer to investigate the initial double integral
representation (\ref{eq:inverse-BF-transform}) using the methods introduced in
{\cite{Ice2021,Part6A}}. Indeed, the two methods are equivalent, but the
second one is more demonstrative and can handle singularities that are more
complicated than simple poles.
\end{rema}
\RED{To show that such phononic crystal fits the general framework of section \ref{sec:general-framework}, it remains to show that $F$ can be analytically continued when $\bdxi$ becomes complex. It can indeed be shown to be true, but the proof is a bit technical. Hence, in order not to break the flow of the paper, we have decided to give it in Appendix \ref{sec:intro-psi-xi}.}

\subsubsection{Singularity set and integral formulations}	
Given the explicit formula (\ref{eq:continuous-explicit-F}) and its analytic
property, we therefore have access to the singular set of $F$, which happens
to be a set of polar singularities. It is given by
\begin{align*}
&H (k, \varkappa) = \bigcup_{j \in \mathbb{N}} \sigma_j (k, \varkappa), \text{ where }\\
&  \sigma_j (k, \varkappa) = \{ \tmmathbf{\xi} \in \mathbb{C}^2,
g_j (\tmmathbf{\xi}; k, \varkappa) = (k + i \varkappa)^2 - \lambda_j
(\tmmathbf{\xi}) = 0 \}.
\end{align*}
Since $\lambda_j (\tmmathbf{\xi})$ is real for real $\tmmathbf{\xi}$ , it is clear that $H (k, \varkappa) \cup
\tmmathbf{\mathcal{B}}= \emptyset$ if $\varkappa > 0$ so the integral
representation (\ref{eq:inverse-BF-transform}) is valid (point
\tmcolor{blue}{P1}). Moreover, because $\lambda_j (\tmmathbf{\xi})$ are real-valued smooth
functions on $\tmmathbf{\mathcal{B}}$, its level sets are in general smooth
curves (point \tmcolor{blue}{P2}), and so, apart from some potential
pathological points where \tmcolor{blue}{P3} might be violated, $F$ does
indeed have the real property and the modified integral representation
(\ref{eq:inverse-BF-transform-deformed}) would normally be required. Therefore, this
phononic crystal example fits perfectly within the general framework discussed in
section \ref{sec:general-framework}.

Motivated by the fact that the general framework of section
\ref{sec:general-framework} is indeed relevant, we will now discuss the
deformation process $\tmmathbf{\mathcal{B}} \rightarrow \tmmathbf{\Gamma}$ in
more details.

\section{More on the indented surface of integration $\bdG$
}\label{sec:deformation-process}

\subsection{Type of deformed surfaces}

Throughout this article, we will only consider deformed surfaces of a
certain class, namely those that can be parametrised by their real parts. In
other words, they can be viewed as a small complex perturbation of
$\tmmathbf{\mathcal{B}}$. More specifically, we will consider surfaces that
can be written as follows:
\begin{eqnarray*}
\tmmathbf{\Gamma} & = & \left\{ \tmmathbf{\xi}_{\tmmathbf{\Gamma}}
=\tmmathbf{\xi}+ i\tmmathbf{\eta} (\tmmathbf{\xi}) \in \mathbb{C}^2,
\tmop{where} \tmmathbf{\xi} \in \tmmathbf{\mathcal{B}} \infixand
\tmmathbf{\eta}= (\eta_1, \eta_2)^{\transpose} \in \mathbb{R}^2 \right\} .
\end{eqnarray*}
We suppose that $\tmmathbf{\eta}$ is small $|\tmmathbf{\eta}|\ll 1$ and  impose periodicity: $\tmmathbf{\eta} (-
\pi, \xi_2) =\tmmathbf{\eta} (\pi, \xi_2)$ and $\tmmathbf{\eta} (\xi_1, - \pi)
=\tmmathbf{\eta} (\xi_1, \pi)$. Such a surface $\tmmathbf{\Gamma}$ is hence
completely described by a real vector field $\tmmathbf{\eta}$ defined over
$\tmmathbf{\mathcal{B}}$. Moreover, for any $\tmmathbf{\xi} \in
\tmmathbf{\mathcal{B}}$, we have $\tmop{Re}
[\tmmathbf{\xi}_{\tmmathbf{\Gamma}} (\tmmathbf{\xi})] =\tmmathbf{\xi}$ and
$\tmop{Im} [\tmmathbf{\xi}_{\tmmathbf{\Gamma}} (\tmmathbf{\xi})]
=\tmmathbf{\eta} (\tmmathbf{\xi})$.

\subsection{Bypassing the singularities}

We have shown in \cite{Part6A}(Theorem 3.2) that for such a surface not to
intersect the singular set $H (k, 0)$ of a function with the real property, it
is sufficient for the vector field $\tmmathbf{\eta}$ to satisfy two
conditions:
\begin{enumerate}
\item For any $\tmmathbf{\xi} \in H' (k) \equiv H (k, 0) \cap
\tmmathbf{\mathcal{B}}$, $\tmmathbf{\eta} (\tmmathbf{\xi}) \neq 0$.	
\item The vector field $\tmmathbf{\eta}$ should not be tangent to any
irreducible component of $H' (k)$.
\end{enumerate}
If these conditions are satisfied, we can always find a smooth
positive real factor function $\varepsilon (\tmmathbf{\xi})$ say such that the
surface parametrised by the vector field $\varepsilon (\tmmathbf{\xi})
\tmmathbf{\eta} (\tmmathbf{\xi})$ does not intersect $H (k, 0)$ and
$\tmmathbf{\mathcal{B}}$ can be deformed continuously to this surface without
intersection $H (k, \varkappa)$ for small enough $\varkappa > 0$. From now on,
we will omit this factor and just use $\tmmathbf{\eta}$ when describing such
surface.

\begin{figure}[h]
\centering{
	\includegraphics[width=0.25\textwidth]{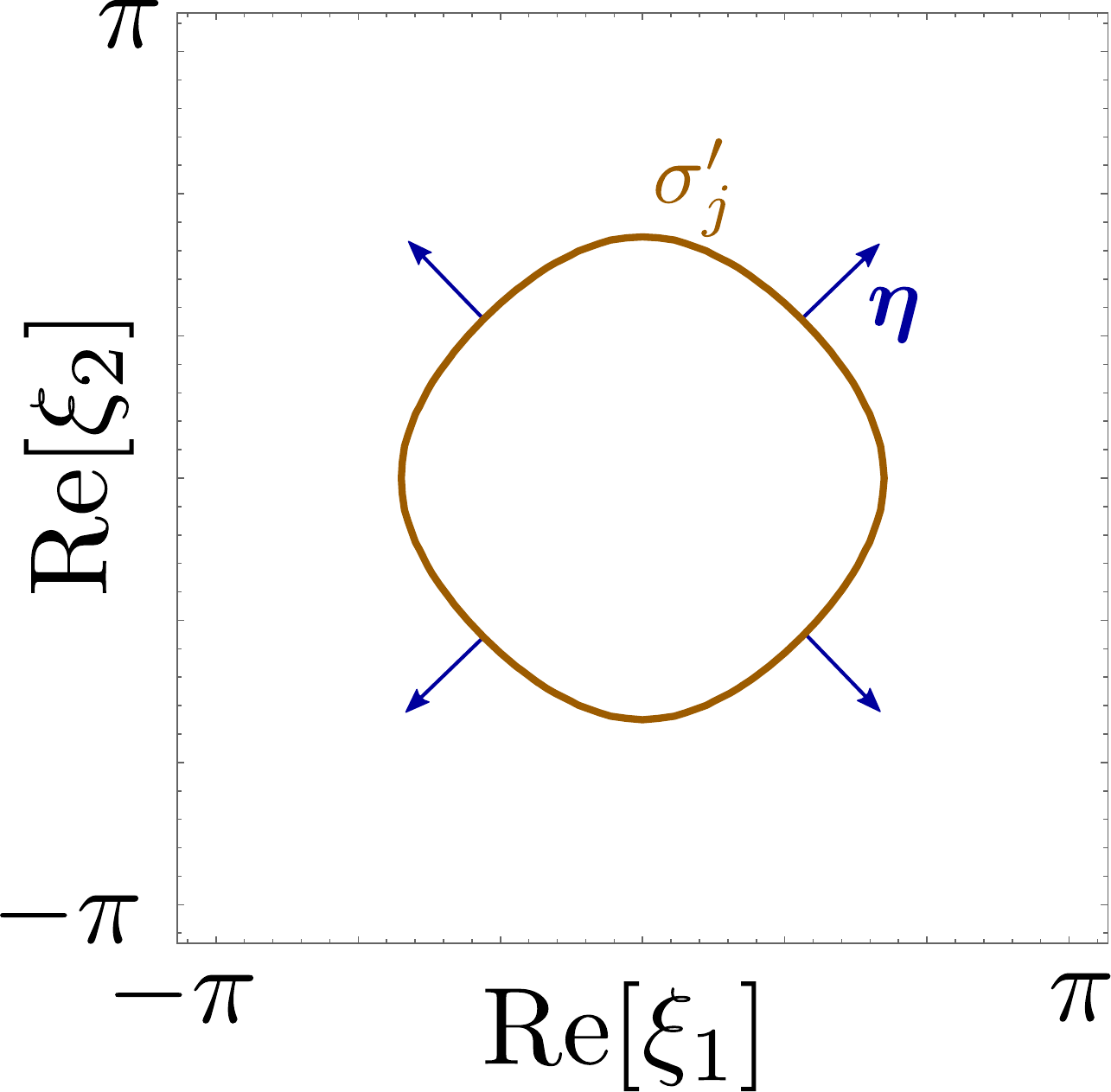}
	\quad \includegraphics[width=0.25\textwidth]{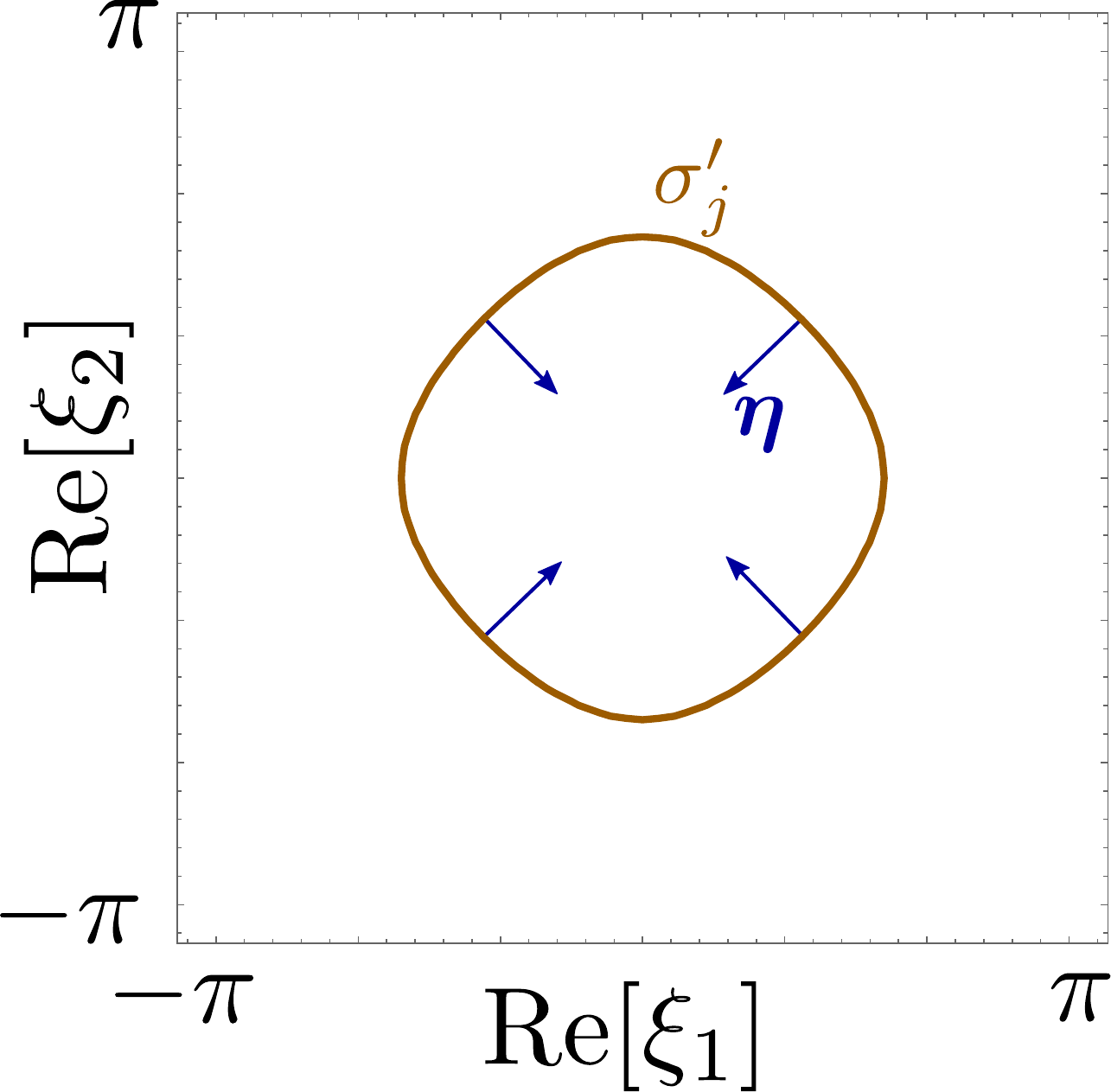}
}
\caption{The arrow indicating two different bypasses of $H' (k)$ in the case
	of the discrete lattice for some $k$ chosen such that $0 < k < 2$}
\label{fig:two-different-bypasses-closed-real-trace}
\end{figure}

\begin{figure}[h]
\centering{
	\includegraphics[width=0.35\textwidth]{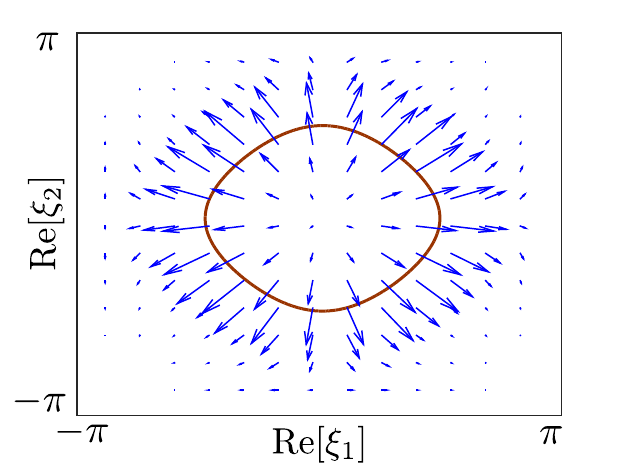}
}
\caption{The vector field $\tmmathbf{\eta}$ that provides admissible  deformation of $\tmmathbf{\mathcal{B}}$ in the case of the discrete lattice for some $k$ chosen such that $0 < k < 2$ as given by (\ref{eq:explicit-param-bypass}) using the + sign.}
\label{fig:vector-field}
\end{figure}

We have also shown in {\cite{Part6A}} that there are only two ways that such a
deformed surface can bypass an irreducible component $\sigma_j' (k)$ of $H'
(k)$. The difference between these two possibilities is the side of $\sigma_j'
(k)$ that $\tmmathbf{\eta}$ is pointing to. Moreover, knowing how
$\tmmathbf{\Gamma}$ bypasses $\sigma_j' (k)$ at one point is enough to know
the type of bypass everywhere. Indeed, because of the two properties above, it
is not possible for $\tmmathbf{\eta}$ to point to two different sides of
$\sigma_j' (k)$. A full understanding of such bypass can hence be visualised
by a simple arrow atop the real traces that indicates which side
$\tmmathbf{\eta}$ is pointing to. This arrow is the arrow part of the
{\tmem{bridge and arrow}} notation developed in {\cite{Assier2019c,Part6A}}.
For the purpose of this article, we will not need to introduce the bridge. For
instance, in Figure \ref{fig:two-different-bypasses-closed-real-trace} we
illustrate the two possible ways in which a surface $\tmmathbf{\Gamma}$ can
bypass the closed loop that is the singularity trace in the case of the
discrete lattice for $0 < k < 2$. Such  bypasses can be realized using the following parametrizations:
\begin{align}
\bdeta(\bdxi)&=(\eta_1(\bdxi),\eta_2(\bdxi)) \nonumber \\ &= \pm \exp\left(-(F(\bdxi;k,0))^{-2}\right)(\sin(\xi_1),\sin(\xi_2)), 
%\quad \eta_2 = \pm \exp\left(-F^{-2}\right)\sin(\xi_2), 
\label{eq:explicit-param-bypass}
\end{align} 
where $F(\tmmathbf{\xi};k,\varkappa)$ is given by (\ref{eq:discrete-explicit-expr}).   The vector field $\tmmathbf{\eta}$ for the bypass in Figure~\ref{fig:two-different-bypasses-closed-real-trace}, left is shown in Figure~\ref{fig:vector-field}.  

\subsection{Choosing the correct bypass\RED{, link with radiation condition}}

For a given $k$, for each irreducible component $\sigma_j' (k)$ of $H' (k)$,
we have therefore the choice between two types of bypasses. 

The procedure \RED{to make the right choice is easy to implement and can simply be summarised} as follows. Draw $\sigma_j' (k)$ and
pick $\tmmathbf{\xi}^{\star}$ on it, then for some small $\delta k > 0$, draw
$\sigma_j' (k + \delta k)$. The \RED{group velocity vector $\nabla_{\!\!\bdxi}k$ at $\bdxi^\star$} %$\nabla f(\tmmathbf{\xi}^{\star})$ 
should hence point towards $\sigma_j' (k + \delta k)$\RED{. Pick} $\tmmathbf{\eta} (\tmmathbf{\xi}^{\star})$ \RED{to} point to the opposite
side, as illustrated in Figure \ref{fig:bypass-choice}. \RED{The proof of why this procedure is valid is given in Appendix \ref{app:procedure-bypass-choice}}.

\begin{rema}
Another, equivalent, way of choosing the bypass is described in {\cite{Part6A}}.
\end{rema}

\begin{exa}
Let us consider the discrete lattice example, with $0 < k < 2$, for which
$H' (k)$ is a smooth closed curve surrounding $0$. 
Draw $\sigma_1' (k)$ and pick
$\tmmathbf{\xi}^{\star}$, draw $\sigma_1' (k + \delta k)$, deduce the
direction of \RED{the group velocity vector}, and hence
deduce the direction of $\tmmathbf{\eta} (\tmmathbf{\xi}^{\star})$.
\end{exa}

\begin{figure}[h]
\centering{
	\includegraphics[height=0.25\textwidth]{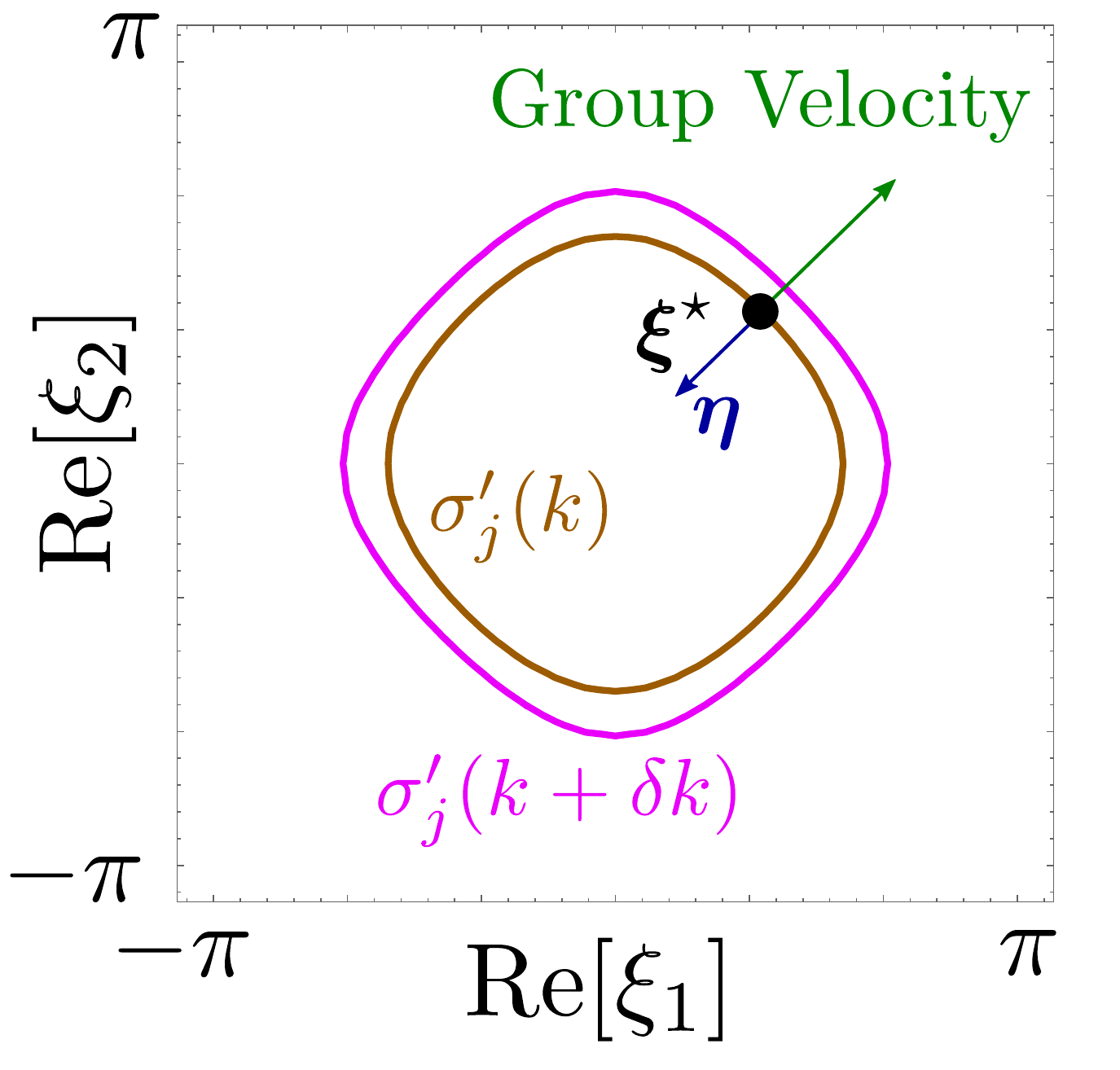}}
\caption{Correct choice of bypass for the discrete lattice example}
\label{fig:bypass-choice}
\end{figure}

\begin{exa}
Here we provide an example for which the above method fails. Consider the
function
\begin{eqnarray*}
	g_j (\tmmathbf{\xi}; k, \varkappa) & = & (\xi_1 - 1) (\xi_1 + 1) + \xi_2^2
	+ 2 (k + i \varkappa) \xi_2
\end{eqnarray*}
The point $\tmmathbf{\xi}^{\star} = (\pm 1, 0) \in \sigma_j' (k)$ for any
real $k$. Here  $\frac{\partial
	G_j}{\partial \tilde{k}} (\tmmathbf{\xi}^{\star} ; k^{\star}) = 0$. What happens in that case is that $\sigma_j' (k)$ and
$\sigma_j' (k + \delta k)$ are circles that intersect at
$\tmmathbf{\xi}^{\star}$. This prevents one to
choose a consistent orientation of the indentation vector $\tmmathbf{\eta}$
and this is due to a violation of the point \tmcolor{blue}{P3}.
\end{exa}

This choice is not
just a mathematical trick that allows us to write down an integral as
$\varkappa \rightarrow 0$, but it is intimately related to the definition of a
radiation condition at infinity through a limiting absorption principle (see,
e.g.\ {\cite{Radosz-2010}}). In fact we find that thinking of radiation
condition in this way might be more natural: the radiation condition is the
correct choice of bypasses.

\RED{By radiation condition here we mean a condition that ensures the problem is well-posed and has a unique, physical, solution. In classical scattering theory, for the problem of an incident plane wave impinging on a compact obstacle surrounded by an infinite acoustic media for instance, this is generally done in two, equivalent ways. The first way is to ensure that the Sommerfeld radiation condition is satisfied (this can be done for real wavenumber $k$). The second way is to use the limiting absorption principle (let the wavenumber $k$ have a small positive real part, and ensure that the scattered field is exponentially decaying at infinity). For periodic media, there is no such thing as a Sommerfeld radiation condition. Hence the approach usually taken is that of limiting absorption principle. In-depth consideration on this topic can for example be found in \cite{fliss2016solutions} and \cite{Radosz-2010}. The main difficulty with the limiting absorption principle is to show that the admissible solution found for $\varkappa>0$ has a limit as $\varkappa\to0$ and to write down this limit. This is exactly what we are doing here. We effectively make sense of this limit by writing it explicitly as a well-behaved integral. The only information we need to write this integral is the correct choice of bypass. Once this is done, we then have made sense of the limit $\varkappa\to0$ and hence effectively provided a radiation condition. A corollary to this is that if one cannot find an appropriate surface deformation, then there are no possible radiation conditions. Special frequencies (corresponding to degeneracies, see section \ref{sec:vanishing-loops}) at which no radiation conditions can exist have been exhibited for discrete settings in e.g.\ \cite{shaban2001radiation} and \cite{Kapanadze2021}.

Understanding the radiation condition for periodic structures, and the asymptotic behaviour of waves propagating in such structures, is essential. This could for instance pave the way to effective conditions (akin to perfectly matched layers in classical scattering theory) to be used for numerical computations on the boundary of a truncated domain to mimic an infinite periodic domain. Efforts in that direction for periodic waveguides have already taken place (see e.g.\ \cite{fliss2016solutions}).}

\section{A brief summary of the general
procedure}\label{sec:summary-procedure}

\begin{figure*}%[h]
\centering{
\includegraphics[width=0.2\textwidth]{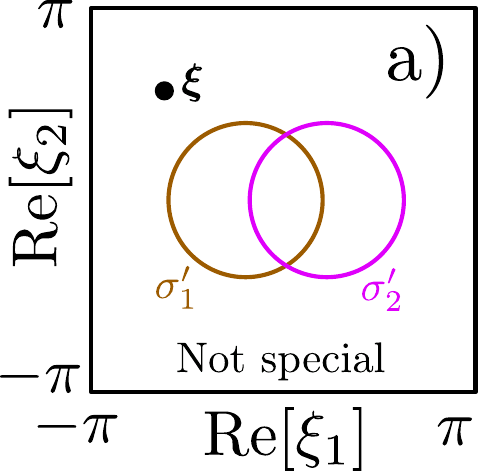}\quad\includegraphics[width=0.2\textwidth]{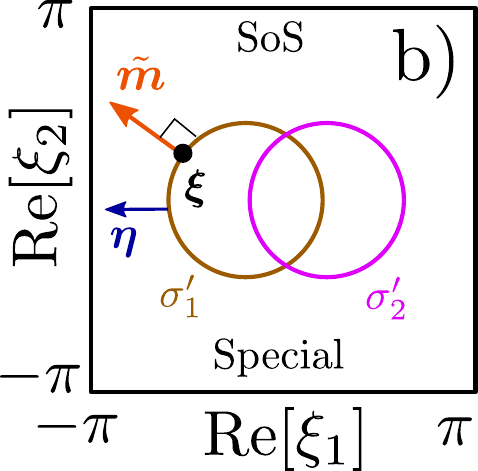}\quad\includegraphics[width=0.2\textwidth]{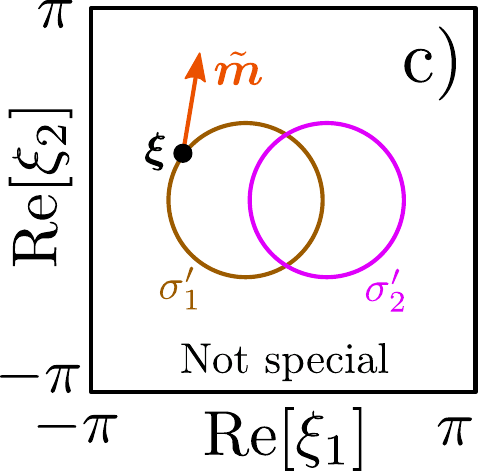}\quad\includegraphics[width=0.2\textwidth]{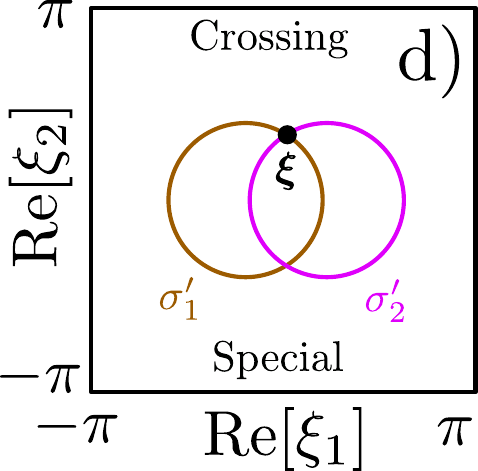}
}
\caption{Schematic illustration of some special and not special points.}
\label{fig:specialornot}
\end{figure*}

Let us consider a wave field $u$ associated to a spectral function $F$ that
satisfy the general framework of section \ref{sec:general-framework}. Let us assume a correct choice of bypass (orientation of the vector $\tmmathbf{\eta}$) for each irreducible components $\sigma_j' (k)$ of $H' (k)$, as explained in section \ref{sec:deformation-process}. We assume hereafter that $\varkappa = 0$ since this parameter was only necessary to define the correct indentation
$\tmmathbf{\Gamma}$. We will therefore simplify the notations by removing the
$\varkappa$ dependency: $H (k, 0) \equiv H (k)$, $u (\tmmathbf{r}+ \Lambda
\tmmathbf{m}; \tmmathbf{r}_s ; k, 0) \equiv u (\tmmathbf{r}+ \Lambda
\tmmathbf{m}; \tmmathbf{r}_s ; k)$ and $F (\tmmathbf{r}; \tmmathbf{r}_s ;
\tmmathbf{\xi}; k, 0) \equiv F (\tmmathbf{r}; \tmmathbf{r}_s ; \tmmathbf{\xi};
k)$. Recall that we are interested in obtaining a far-field asymptotic
expansion for $u$. To do this, we fix an observation direction
$\tilde{\tmmathbf{m}}$ and study $u (\tmmathbf{r}+ \Lambda \tmmathbf{m};
\tmmathbf{r}_s ; k, 0)$ for $\tmmathbf{m}= N \tilde{\tmmathbf{m}}$ as the
integer $N \rightarrow \infty$. For each observation direction, our aim is to
provide an expansion of the type
\begin{eqnarray}
u (\tmmathbf{r}+ \Lambda \tmmathbf{m}; \tmmathbf{r}_s ; k) & \underset{N
	\rightarrow \infty}{\approx} & \sum_q u_{\tmop{loc}}^{(q)} (\tmmathbf{m};
\tmmathbf{r}; \tmmathbf{r}_s ; k),  \label{e:03002}
\end{eqnarray}
where each of the wave components $u_{\tmop{loc}}^{(q)}$ are known explicitly.
According to the locality principle (see {\cite{Ice2021,Part6A}}), we expect
each $u_{\tmop{loc}}^{(q)}$ to result from local considerations around a
special point $\tmmathbf{\xi}_q \in \tmmathbf{\mathcal{B}}$. These special
points can be found through a second stage of integration surface deformation.
We deform $\tmmathbf{\Gamma}$ into another surface $\tmmathbf{\Gamma}'$
according to the following properties: \\ \ \\
%\begin{itemize}
\indent $\bullet$ $\tmmathbf{\Gamma}'$ is also of the type described in section
\ref{sec:deformation-process}, and can be described by a real vector field
$\tmmathbf{\eta}'$ \\
\indent $\bullet$ The deformation from $\tmmathbf{\Gamma}$ to $\tmmathbf{\Gamma}'$ is
\tmtextit{admissible}: this is a continuous deformation during which the
singularity set of $F$ is not hit. As a result, $\tmmathbf{\eta}$ and
$\tmmathbf{\eta}'$ points to the same side of the real traces of
$F$. \\ 
\indent $\bullet$ The deformation is \tmtextit{desired}: for $\tmmathbf{\xi} \in
\tmmathbf{\mathcal{B}}$, the exponential factor $e^{- i\tmmathbf{m} \cdot
\tmmathbf{\xi}_{\tmmathbf{\Gamma}'} (\tmmathbf{\xi})}$ in the integrand of
(\ref{eq:inverse-BF-transform-deformed}) should be exponentially decaying as
$N \rightarrow \infty$. This is true provided that $\tilde{\tmmathbf{m}}
\cdot \tmmathbf{\eta}' (\tmmathbf{\xi}) < 0$. \\	
\indent $\bullet$ Points $\tmmathbf{\xi} \in \tmmathbf{\mathcal{B}}$ in the
neighbourhood of which it is impossible to find a deformation
$\tmmathbf{\Gamma} \rightarrow \tmmathbf{\Gamma}'$ that is both admissible
and desired are called special points. \\ \ \\		
%\end{itemize}
\indent We are now in a position of summarising some of the findings of
{\cite{Part6A}} regarding which points are special or not. \\ \ \\
%\begin{itemize}
\indent $\bullet$ Non-singular points are not special. By non-singular we mean points
$\tmmathbf{\xi} \in \tmmathbf{\mathcal{B}} \setminus H' (k)$. This in
particular implies that if we are within a band gap, i.e.\ if $k$ is such
that $H' (k) \cap \tmmathbf{\mathcal{B}}= \emptyset$, then the Green's
function decays exponentially at infinity, as expected. See Figure
\ref{fig:specialornot} a).\\ 
\indent $\bullet$ Saddles on Singularities (SoS) are special. A SoS is a point
$\tmmathbf{\xi} \in \sigma_j' (k)$ such that at $\tmmathbf{\xi}$ we have
$\tilde{\tmmathbf{m}} \perp \sigma_j' (k)$ and $\tilde{\tmmathbf{m}} \cdot
\tmmathbf{\eta}> 0$. See Figure \ref{fig:specialornot} b). Note that this
case is well known in the theory of anisotropic media
{\cite{Lighthill1960}}: the wave vector of a wave propagating in a certain
direction is given by the point where the direction is orthogonal to the
slowness surface. This is a generalisation of the concept of the group
velocity.\\ 
\indent $\bullet$ Singular points belonging to only one irreducible component are not
special unless they are SoS. See Figure \ref{fig:specialornot} c).\\	
\indent $\bullet$ Tangential crossings between two real traces are not special unless
they are $\tmop{SoS}$.\\ 
\indent $\bullet$ Transverse crossings between two real traces are special. See Figure
\ref{fig:specialornot} d). They lead to some asymptotic contribution for
$\tilde{\tmmathbf{m}}$ belonging to one `active' quadrant of the real plane. \\ \ \\

In the next section, we will give the asymptotic formulae associated with the
special points described above. For isolated SoS (section \ref{sec:SoS}) and
transverse crossings (section \ref{sec:transverse-crossing}), the formulae can
be obtained directly from {\cite{Part6A}}. The only difference is that here
$\tmmathbf{\Gamma}$ and $\tmmathbf{\mathcal{B}}$ are finite, while they were
doubly infinite in {\cite{Part6A}}. However, since only local considerations
are needed to obtain the results, everything still holds.

\section{Asymptotics associated with special
points}\label{sec:asymptotics-special-points}

\subsection{Isolated Saddle on a Singularity (SoS)}\label{sec:SoS}

Let us start by fixing $\tilde{\tmmathbf{m}}$, $\tmmathbf{r}$ and
$\tmmathbf{r}_s$ and $k$ and a choice of bypass for $H' (k)$. Now consider a
SoS $\tmmathbf{\xi}^{\star}$ on an irreducible real trace component
$\sigma'_j$ of $H' (k)$ with defining function $g_j (\tmmathbf{\xi}; k)$. By
definition of a SoS, we have $\tilde{\tmmathbf{m}} \perp \sigma_j'$ and
$\tilde{\tmmathbf{m}} \cdot \tmmathbf{\eta}> 0$ at $\tmmathbf{\xi}^{\star}$,
see Figure \ref{fig:specialornot} b). Let us further assume that
$\tmmathbf{\xi}^{\star}$ belongs to no other irreducible component, and that
$\sigma_j'$ has non-zero curvature at $\tmmathbf{\xi}^{\star}$. The latter
assumption implies that the SoS is isolated. However, if $\sigma'_j$ were to
have straight fragments, non-isolated SoS might occur and non-local
considerations would be needed. Straight fragments were not studied in
{\cite{Part6A}}, but we will consider \RED{them} in the present paper, see section
\ref{sec:flat-sing}.

To write down the asymptotics associated with this point, we first need to
introduce some quantities that characterise the singularity. Define
$\text{\tmverbatim{a}}^{\star}_j = \partial_{\xi_1} g_j
(\tmmathbf{\xi}^{\star} ; k)$, $\text{\tmverbatim{b}}^{\star}_j =
\partial_{\xi_2} g_j (\tmmathbf{\xi}^{\star} ; k)$ and the normal vector
$\tmmathbf{n}^{\star} = \left( \text{\tmverbatim{a}}^{\star}_j,
\text{\tmverbatim{b}}^{\star}_j \right)^{\transpose} / \sqrt{\left(
\text{\tmverbatim{a}}^{\star}_j \right)^2 + \left(
\text{\tmverbatim{b}}^{\star}_j \right)^2}$. Because $\tilde{\tmmathbf{m}}
\perp \sigma'_j$, there exists a scalar $s = \pm 1$ defined such that
$\tilde{\tmmathbf{m}} = s | \tilde{\tmmathbf{m}} | \tmmathbf{n}^{\star}$. The
sign of $s$ depends on the choice of bypass and defining function. An
equivalent definition of $s$ is $s = \tmop{sign} (\tmmathbf{\eta}^{\star}
\cdot \tmmathbf{n}^{\star})$, where $\tmmathbf{\eta}^{\star} =\tmmathbf{\eta}
(\tmmathbf{\xi}^{\star})$ and $\tmmathbf{\eta}$ is the indentation vector
field associated to $\tmmathbf{\Gamma}$. For brevity, we also introduce
$N^{\star} = N | \tilde{\tmmathbf{m}} | / \sqrt{\left(
\text{\tmverbatim{a}}^{\star}_j \right)^2 + \left(
\text{\tmverbatim{b}}^{\star}_j \right)^2}$.

We further assume that the leading component of $F$ responsible for the
activity of $\tmmathbf{\xi}^{\star}$ has the following behaviour as
$\tmmathbf{\xi} \rightarrow \tmmathbf{\xi}^{\star} $:
\begin{eqnarray*}
F (\tmmathbf{r}; \tmmathbf{r}_s ; \tmmathbf{\xi}; k) & \approx & 4 \pi^2 A
(\tmmathbf{r}; \tmmathbf{r}_s)  (g_j (\tmmathbf{\xi}; k))^{- \mu} .
\end{eqnarray*}
We use the change of variable $\tmmathbf{\xi} \leftrightarrow
\tmmathbf{\zeta}$ where
\begin{eqnarray*}
\zeta_1 = \text{\tmverbatim{b}}^{\star}_j (\xi_1 - \xi_1^{\star}) -
\text{\tmverbatim{a}}^{\star}_j (\xi_2 - \xi_2^{\star}) & \tmop{and} & \zeta_2
= g_j (\tmmathbf{\xi}; k),
\end{eqnarray*}
ensuring that $\zeta_1$ is tangent to $\sigma'_j$ and $\zeta_2$ is transverse
to $\sigma'_j$. We therefore have
\begin{align*}
&\tmmathbf{\xi}-\tmmathbf{\xi}^{\star} = \Psi \tmmathbf{\zeta}+\mathcal{O}
(\zeta_1^2 + \zeta_1 \zeta_2 + \zeta_2^2),  \text{ where } \\ & \Psi =
\frac{1}{\left( \text{\tmverbatim{a}}^{\star}_j \right)^2 + \left(
\text{\tmverbatim{b}}^{\star}_j \right)^2} \left( \begin{array}{cc}
\text{\tmverbatim{b}}^{\star}_j & \text{\tmverbatim{a}}^{\star}_j\\
- \text{\tmverbatim{a}}^{\star}_j & \text{\tmverbatim{b}}^{\star}_j
\end{array} \right) .
\end{align*}
Moreover, the non-zero curvature property leads to the existence of a unique
$\alpha \neq 0$ such that
\begin{eqnarray*}
\zeta_2 & \approx & \text{\tmverbatim{a}}^{\star}_j (\xi_1 - \xi_1^{\star}) +
\text{\tmverbatim{b}}^{\star}_j (\xi_2 - \xi_2^{\star}) - \alpha \zeta_1^2
+\mathcal{O} (\zeta_2^2 + \zeta_1 \zeta_2) .
\end{eqnarray*}
We have now introduced all the required quantities to write down the far-field
wave component $u_{\tmop{loc}}$ resulting from the special point
$\tmmathbf{\xi}^{\star}$, as presented in {\cite{Part6A}}:
%\begin{widetext}
\begin{align}
u_{\tmop{loc}} (\tmmathbf{m}; \tmmathbf{r}; \tmmathbf{r}_s ; k) &
\underset{N \rightarrow \infty}{\approx}  \frac{2 \pi A (\tmmathbf{r};
	\tmmathbf{r}_s) e^{- i\tmmathbf{m} \cdot \tmmathbf{\xi}^{\star}} \sqrt{\pi}
	e^{- is \mu \pi / 2}}{\Gamma (\mu) \left( \left(
	\text{\tmverbatim{a}}^{\star}_j \right)^2 + \left(
	\text{\tmverbatim{b}}^{\star}_j \right)^2 \right)(N^{\star})^{3 /2 -\mu}}
\times \left\{ \begin{array}{ccc}
	e^{- i \pi / 4} / \sqrt{s \alpha} & \tmop{if} & s \alpha > 0\\
	e^{i \pi / 4} / \sqrt{- s \alpha} & \tmop{if} & s \alpha < 0
\end{array} \right.,
\label{SOS_asympt}	 
\end{align}
%\end{widetext}
where $\Gamma$ is the usual gamma function.

It is worth to mention, that a particular case of asymptotic (\ref{SOS_asympt}) for $(g_j (\tmmathbf{\xi}; k))^{\RED{\mu}}$ being a polynomial was obtained in \cite{Lighthill1960,Lighthill78} by methods of 1D complex analysis applied to the series-integral representation (\ref{e:02002}). 

\begin{exa}
Consider the Green's function of the discrete lattice. Let $k$ be small but finite. Then, near the singularity (where $\xi_1$ and $\xi_2$ are also small), 
(\ref{eq:discrete-explicit-expr}) can be approximated as follows:
\begin{equation}
\label{F_discr_small_k}
F \approx \frac{1}{k^2 - \xi_1^2 - \xi_2^2} \cdot
\end{equation}
Thus, the real trace of the irreducible singularity is a circle, and for any given $\tmmathbf{m}$ there is one active SoS given by $\bdxi^\star=-k \bdm/|\bdm|$. %:
%\begin{equation}
%\xi^\star_1 =\pm\frac{m_1k}{\sqrt{m_1^2+m_2^2}},\quad \xi^\star_2 =\pm\frac{m_2k}{\sqrt{m_1^2+m_2^2}}, 
%\end{equation}
%one passive and one active.  
Using (\ref{SOS_asympt}) we obtain
\begin{equation}
u_{\tmop{loc}} (\tmmathbf{m}; \tmmathbf{r}; \tmmathbf{r}_s ; k)	\underset{|\bdm| \rightarrow \infty}{\approx} -\frac{e^{i\pi/4}}{2\sqrt{2\pi}} \frac{e^{ik|\bdm|}}{\sqrt{k |\bdm|}},
\end{equation}
where we have used that $A=\tfrac{1}{4\pi^2}$, $(
\text{\tmverbatim{a}}^{\star}_j )^2 + (
\text{\tmverbatim{b}}^{\star}_j )^2=4k^2$, $N^\star=|\bdm|/(2k)$, $\mu=1$, $s=1$, and $\alpha=\tfrac{1}{4k^2}$. 
%\begin{equation}
%u_{\tmop{loc}} (\tmmathbf{m}; \tmmathbf{r}; \tmmathbf{r}_s ; k) \approx -\frac{e^{-ik\sqrt{m_1^2+m_2^2}-i\pi/4}}{2\sqrt{2\pi}} \cdot
%\end{equation}
The latter has a simple physical interpretation. It is the far-field asymptotics of the Green's function associated to the continuous problem, \RED{namely $-\tfrac{i}{4}H_0^{(1)}(k|\bdm|)$}. Indeed, when $k$ is small (\ref{eq:discrete-governing}) can be considered as a finite difference approximation of the Helmholtz equation on a plane with a point source at the origin.  
%{\color{red} \textbf{Note on eq above}: AIK: you will see that I corrected it slightly, including a square root term on the denominator and a change of sign in the exponential.} 
\end{exa}

\subsection{Transverse crossings}\label{sec:transverse-crossing}

Let us start by fixing $\tilde{\tmmathbf{m}}$, $\tmmathbf{r}$ and
$\tmmathbf{r}_s$ and $k$ again. Consider a point $\tmmathbf{\xi}^{\star} \in
\sigma_1' \cap \sigma_2'$, where $\sigma_1'$ and $\sigma_2'$ are two
irreducible components of $H' (k)$ with respective defining functions $g_1
(\tmmathbf{\xi}; k)$ and $g_2 (\tmmathbf{\xi}; k)$. Unlike for the SoS, the
position of such a point does not depend on $\tilde{\tmmathbf{m}}$. Define the
real quantities $\text{\tmverbatim{a}}_{1, 2}^{\star} = \partial_{\xi_1} g_{1,
2} (\tmmathbf{\xi}^{\star} ; k)$ and $\text{\tmverbatim{b}}^{\star}_{1, 2} =
\partial_{\xi_2} g_{1, 2} (\tmmathbf{\xi}^{\star} ; k)$. The crossing being
transverse, means that the quantity $\Delta^{\star} =
\text{\tmverbatim{a}}_1^{\star} \text{\tmverbatim{b}}^{\star}_2 -
\text{\tmverbatim{a}}_2^{\star} \text{\tmverbatim{b}}^{\star}_1$ is not zero.
Because it is always possible to swap notations for $\sigma_1'$ and
$\sigma_2'$, we can assume without loss of generality that $\Delta^{\star} >
0$. We can naturally define the normals $\tmmathbf{n}_{1, 2}^{\star}$ as
$\tmmathbf{n}_1^{\star} = \left( \text{\tmverbatim{a}}^{\star}_1,
\text{\tmverbatim{b}}^{\star}_1 \right)^{\transpose} / \sqrt{\left(
\text{\tmverbatim{a}}^{\star}_1 \right)^2 + \left(
\text{\tmverbatim{b}}^{\star}_1 \right)^2}$ and $\tmmathbf{n}_2^{\star} =
\left( \text{\tmverbatim{a}}^{\star}_2, \text{\tmverbatim{b}}^{\star}_2
\right)^{\transpose} / \sqrt{\left( \text{\tmverbatim{a}}^{\star}_2 \right)^2 + \left(
\text{\tmverbatim{b}}^{\star}_2 \right)^2}$ and the sign factors $s_{1, 2} =
\pm 1$ can be defined by $s_1 = \tmop{sign} (\tmmathbf{\eta}^{\star} \cdot
\tmmathbf{n}_1^{\star})$ and $s_2 = \tmop{sign} (\tmmathbf{\eta}^{\star} \cdot
\tmmathbf{n}_2^{\star})$, where $\tmmathbf{\eta}^{\star} =\tmmathbf{\eta}
(\tmmathbf{\xi}^{\star})$ and $\tmmathbf{\eta}$ is the indentation vector
field associated to $\tmmathbf{\Gamma}$.

Let us now assume that the leading component of $F$ responsible for the
activity of $\tmmathbf{\xi}^{\star}$ has the following behaviour as
$\tmmathbf{\xi} \rightarrow \tmmathbf{\xi}^{\star}$
\begin{align}
F (\tmmathbf{r}; \tmmathbf{r}_s ; \tmmathbf{\xi}; k) & \approx  4 \pi^2 A
(\tmmathbf{r}; \tmmathbf{r}_s) (g_1 (\tmmathbf{\xi}))^{- \mu_1} (g_2
(\tmmathbf{\xi}))^{- \mu_2} .  \label{eq:local-behav-crossing}
\end{align}
We have now introduced all the required quantities to write down the far-field
wave component $u_{\tmop{loc}}$ resulting from the special point
$\tmmathbf{\xi}^{\star}$, as presented in {\cite{Part6A}}:
%\begin{widetext}
\begin{align}
u_{\tmop{loc}} &
\underset{N \rightarrow \infty}{\approx} \frac{4 \pi^2 A (\tmmathbf{r};
	\tmmathbf{r}_s) e^{- i\tmmathbf{m} \cdot \tmmathbf{\xi}^{\star}} e^{- i
		\frac{\pi}{2} (s_1 \mu_1 + s_2 \mu_2)}}{\Gamma (\mu_1) \Gamma (\mu_2)
	(\Delta^{\star})^{\mu_1 + \mu_2 - 1}} \times \frac{\mathcal{H} \left( s_1
	\left( m_1 \text{\tmverbatim{b}}^{\star}_2 - m_2
	\text{\tmverbatim{a}}_2^{\star} \right) \right)}{{\left| m_1
		\text{\tmverbatim{b}}^{\star}_2 - m_2 \text{\tmverbatim{a}}_2^{\star}
		\right|^{1 - \mu_1}} } \times \frac{\mathcal{H} \left( s_2 \left( - m_1
	\text{\tmverbatim{b}}^{\star}_1 + m_2 \text{\tmverbatim{a}}_1^{\star}
	\right) \right)}{\left| - m_1 \text{\tmverbatim{b}}^{\star}_1 + m_2
	\text{\tmverbatim{a}}_1^{\star} \right|^{1 - \mu_2}} .
\end{align}
%\end{widetext}

Note that the presence of the Heaviside function $\mathcal{H}$ implies that
such point $\tmmathbf{\xi}^{\star}$ has a region of activity and a region of
inactivity. These depend on the choice of indentation associated to
$\sigma_1'$ and $\sigma_2'$. A visual way of determining those regions of
activity has been provided in {\cite{Part6A}}.

\paragraph{Additive crossings.}The local approximation
(\ref{eq:local-behav-crossing}) of $F$ plays an important part in what is
discussed above. In {\cite{Part6A}}, we have shown that if the crossing
$\tmmathbf{\xi}^{\star}$ is additive, then $\tmmathbf{\xi}^{\star}$ does not
lead to any asymptotic contribution. A crossing $\tmmathbf{\xi}^{\star}$ is
additive, if, locally around $\tmmathbf{\xi}^{\star}$, we can write $F
(\tmmathbf{\xi}) = F_1 (\tmmathbf{\xi}) + F_2 (\tmmathbf{\xi})$, where $F_1$
is not singular on $\sigma_1'$ and $F_2$ is not singular on $\sigma_2'$.
Additive crossings have been shown to occur and to be important in diffraction
theory (see e.g.\
{\cite{Assier2018a,2dcontinuation2021,Assier2019c,Kunz2023}}), but have not been observed in periodic media so far. Given the form
of the eigenexpansion (\ref{eq:continuous-explicit-F-no-conjugate}), it is
likely that crossings resulting from an intersection of two eigenvalues are
going to be additive in the examples considered and will therefore not
contribute to the far-field asymptotics of $u$.

\subsection{Flat singularities and non-isolated SoS}\label{sec:flat-sing}

\begin{figure}[h]
\centering{
\includegraphics[width=0.25\textwidth]{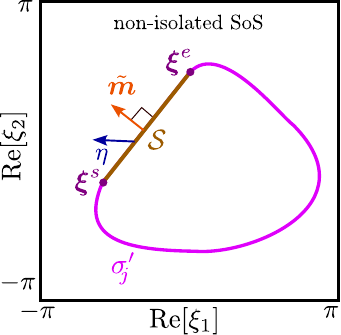}
}
\caption{Schematic illustration of non-isolated SoS}
\label{fig:flatsingularity}
\end{figure}

Let us consider the case when a portion of an irreducible real trace component $\sigma'_j$ is a straight segment $\mathcal{S}$, as shown in Figure~\ref{fig:flatsingularity} . 
Let us further assume that $\tilde{\tmmathbf{m}} \perp \mathcal{S}$ and that $\tilde{\bdm}\cdot\bdeta>0$, so that all points $\bdxi^\star\in\mathcal{S}$ are non-isolated SoS. They all need to be considered to obtain a far-field approximation in the observation direction $\tilde{\bdm}$. This will be done as follows.
%Then, each point of the segment will provide asymptotically non-small contribution to the far-field, i.e. the segment is a non-isolated SoS. Let us write down the asymptotics associated with it. 
%Introduce quantities characterising the singularity as follows. 
Because it is straight, $\mathcal{S}$ can be associated to a defining function $g(\tmmathbf{\xi};k)$ of the form 
\[
g(\tmmathbf{\xi};k) = \tva\xi_1 + \tvb\xi_2 + \tvc,
\] 
with normal vector $\tmmathbf{n} = \left( \tva,
\tvb\right)^{\transpose} / \sqrt{\tva^2 + \tvb^2}$.
Denote the endpoints of $\mathcal{S}$ as $\tmmathbf{\xi}^s$ and $\tmmathbf{\xi}^e$.  Since $\tilde{\tmmathbf{m}}\perp \mathcal{S}$ we can write $\tilde{\tmmathbf{m}} = s | \tilde{\tmmathbf{m}} | \tmmathbf{n}$, with $s=\pm 1$. Similarly to \RED{Section \ref{sec:SoS}}, we introduce $N^\star = N | \tilde{\tmmathbf{m}} | / \sqrt{\tva^2 + \tvb^2}$, as well as the new coordinates: $\zeta_1 = \tvb \xi_1 - \tva\xi_2$ and $\zeta_2
= g (\tmmathbf{\xi}; k)$.

We further assume that, in the vicinity of $\mathcal{S}$, $F$ behaves as 
\[
F (\tmmathbf{r}; \tmmathbf{r}_s ; \tmmathbf{\xi}; k)  \approx  4 \pi^2 A
(\tmmathbf{r}; \tmmathbf{r}_s)  \zeta_2^{- \mu} h(\zeta_1),
\]	 
where $h(\zeta_1)$ is some regular function of the longitudinal coordinate $\zeta_1$.

Using the latter, and the fact that the Jacobian of the change of variables $\tmmathbf{\xi}\to \tmmathbf{\zeta}$ behaves as $(\tva^2+\tvb^2)^{-1}$, the integral (\ref{eq:BF-transform}) can be approximated as follows: 
%\begin{widetext}
\begin{equation}
u_{\text{loc}} \underset{N \rightarrow \infty}\approx 	\frac{A
	(\tmmathbf{r}; \tmmathbf{r}_s)e^{-is\tvc N^\star} }{\tva^2+\tvb^2}\int_{\zeta_1^s}^{\zeta_1^e}h(\zeta_1)d\zeta_1\int_{-\infty+is\epsilon}^{\infty+is\epsilon}\zeta_2^{-\mu} e^{-isN^\star\zeta_2}d\zeta_2,
\end{equation}
%\end{widetext}
where $\zeta_1^s = \tvb \xi^s_1 - \tva\xi^s_2$, $\zeta_1^e = \tvb \xi^e_1 - \tva\xi^e_2$. Note, that only an exponentially small error has been introduced by making the limits of integration in $\zeta_2$ go to infinity. 
The integral in $\zeta_2$ can be evaluated explicitly in terms of the gamma function, which leads to
\begin{equation}
u_{\text{loc}} \underset{N \rightarrow \infty} \approx 	\frac{2\pi A
(\tmmathbf{r}; \tmmathbf{r}_s)e^{-is\tvc N^\star} e^{-is\mu\pi/2} }{\Gamma(\mu)(\tva^2+\tvb^2)( N^\star)^{1-\mu}}\int_{\zeta_1^s}^{\zeta_1^e}h(\zeta_1)d\zeta_1.
\end{equation}
%{\color{green} Agree}
A similar consideration for the case of a polynomial singularity has been given in \cite{Lighthill1960,Lighthill78}. It is an interesting formula, as one can see that it leads to some kind of beam for this specific value of $\bdm$ on which the wave will decay less fast (or not at all if $\mu=1$) tha\RED{n} in the other observation directions associated with a usual isolated SoS behaviour.
% {\color{red}\textbf{Note for AIK}: I think there was a missing $\tvc$ in your original formula as well as minus sign pb (see correction in red). The rest I agree with. It is interesting to see that the power of $N^\star$ is different ($1-\mu$) than it was for an isolated SoS ($3/2-\mu$). In particular, with this new power, it looks like if the singularity is a simple pole (i.e. $\mu=1$), there's no decay. }

Finally, the situation where special points are located very close to each
other or are merging requires special care and will provide intermediate
asymptotics, often referred to as penumbra zone in diffraction theory. In the
next section we will concentrate on {\tmem{degeneracies}}, a phenomenon that
did not appear in our studies on diffraction theory, but that is fundamental
to the study of wave propagation in periodic structures.

\section{Degeneracies}\label{sec:vanishing-loops}

Without losing too much in generality, below we will focus on a field $u$ described
by a spectral function $F$ that can be written as an eigenprojection of the
type (\ref{eq:continuous-explicit-F-no-conjugate}). As mentioned earlier this
encompasses at least the two examples of sections \ref{sec:discrete-example}
and \ref{sec:continuous-example}. We remain interested in the case of real
$k$. For simplicity, we will rename the terms containing the eigenfunctions
as $E_j$ so that we can write
\begin{eqnarray*}
F (\tmmathbf{r}; \tmmathbf{r}_s ; \tmmathbf{\xi}; k) & = & \sum_{j =
0}^{\infty} \frac{E_j (\tmmathbf{r}; \tmmathbf{r}_s ; \tmmathbf{\xi})}{(k^2
- \lambda_j (\tmmathbf{\xi}))}
\end{eqnarray*}
for real $\tmmathbf{\xi} \in \tmmathbf{\mathcal{B}}$. In this section, we will
consider another type of special points that we call {\tmem{degeneracies}}. 
These can be loosely defined as a change of topology in the real traces of $F$
at a specific wavenumber. These occur for example when $H (k)$ contains a loop
that becomes smaller and smaller as $k$ changes, until ultimately shrinking at
some $k = k^{\star}$ or when $H (k)$ contains two branches of a hyperbola that
eventually touch at some $k = k^{\star}$. Indeed, in both cases the surface of integration becomes {\it pinched} between singularities, and thus the resulting integral becomes singular.
Each type of degeneracy needs to be
treated differently. Our aim here is to provide an approximation for $u
(\tmmathbf{r}+ \Lambda \tmmathbf{m}; \tmmathbf{r}_s ; k)$ that is valid when
both $| \tmmathbf{m} | \rightarrow \infty$ and $k \rightarrow k^{\star}$. A
similar aim, though with a different approach base on the high-frequency
homogenization technique was successfully followed in {\cite{Vanel2016}}. We
will start by considering degeneracies that remain a simple eigenvalue,
typically represented by local extrema of the dispersion diagram or hyperbolic
crossing, and then consider degeneracies involving two eigenvalues
simultaneously such as Dirac conical points.

\subsection{General procedure for a simple non-zero eigenvalue}

Let us choose a pair $(\tmmathbf{\xi}^{\star}, k^{\star}) \in
\tmmathbf{\mathcal{B}} \times \mathbb{R}^+$ that is a degeneracy associated to
a simple eigenvalue. In other words, there is a unique $\ell$ such that
$\lambda_{\ell} (\tmmathbf{\xi}^{\star}) = (k^{\star})^2 > 0$ and the topology
of the real trace $\sigma_{\ell}'$ associated to $\lambda_{\ell}$ changes as
$k \rightarrow k^{\star}$. Our aim is to provide a far-field approximation for
$u$ that is valid as $k$ approaches $k^{\star}$. Let us change the variables
from $(\tmmathbf{\xi}, k)$ to $(\tmmathbf{\zeta}, \hat{k})$ by introducing
$\hat{k} = k - k^{\star}$ and defining $\tmmathbf{\zeta}$ as local variables
centred at $\tmmathbf{\xi}^{\star}$, with a non-singular Jacobian matrix
$\Psi^{\star}$ such that as $\tmmathbf{\xi} \rightarrow
\tmmathbf{\xi}^{\star}$ (and hence $\tmmathbf{\zeta} \rightarrow
\tmmathbf{0}$), we have $\tmmathbf{\xi}-\tmmathbf{\xi}^{\star} \approx
\Psi^{\star} \tmmathbf{\zeta}$.

\begin{figure*}%[h]
\centering{
\includegraphics[width=0.45\textwidth]{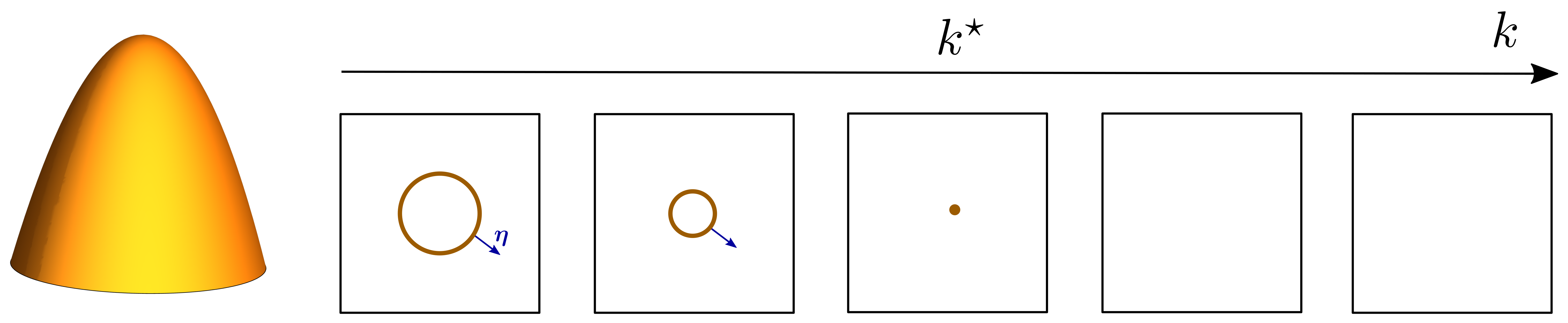} \qquad \includegraphics[width=0.45\textwidth]{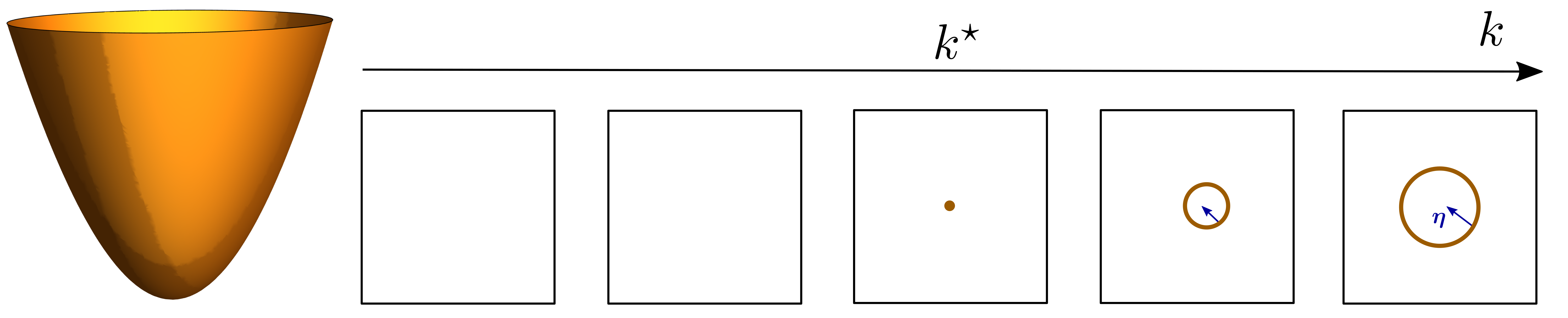}}
\caption{Typical simple eigenvalue degeneracies: local maximum (left) and
minimum (right)}
\label{fig:local-extrememums}
\end{figure*}

Because $\lambda_{\ell} (\tmmathbf{\xi}^{\star}) > 0$, we can, at least
locally, introduce the function $\Omega (\tmmathbf{\zeta})$ defined by
\begin{eqnarray*}
\Omega (\tmmathbf{\zeta}) & = & \sqrt{\lambda_{\ell} (\tmmathbf{\xi})} -
k^{\star},
\end{eqnarray*}
chosen such that, for $(\tmmathbf{\xi}, k)$ in some neighbourhood of
$(\tmmathbf{\xi}^{\star}, k^{\star})$, we have the singularity of $F$
described by $k = \sqrt{\lambda_{\ell} (\tmmathbf{\xi})}$ or in our new
variables, $\hat{k} = \Omega (\tmmathbf{\zeta})$. In a neighbourhood of
$(\tmmathbf{\xi}^{\star}, k^{\star})$, we can therefore write $F$ as follows
\begin{eqnarray*}
F (\tmmathbf{r}; \tmmathbf{r}_s ; \tmmathbf{\xi}; k) & = & \frac{E_{\ell}
(\tmmathbf{r}; \tmmathbf{r}_s ; \tmmathbf{\xi})}{(k^2 - \lambda_{\ell}
(\tmmathbf{\xi}))} + \tmop{regular} \tmop{terms}\\
& \approx & \frac{E_{\ell} (\tmmathbf{r}; \tmmathbf{r}_s ;
\tmmathbf{\xi}^{\star})}{{(k^{\star} + \hat{k})^2}  - (k^{\star} + \Omega
(\tmmathbf{\zeta}))^2} + \tmop{regular} \tmop{terms},
\end{eqnarray*}
where we have used the definition of $\Omega$ and the fact that $E_{\RED{\ell}}
(\tmmathbf{r}; \tmmathbf{r}_s ; \tmmathbf{\xi})$ is analytic around
$\tmmathbf{\xi}^{\star}$. According to {\cite{Part6A}}, the regular terms do
not contribute to the asymptotics of $u$, and we are just left with the
contribution
%\begin{widetext}
\begin{eqnarray*}
u_{\tmop{loc}} (\tmmathbf{m}; \tmmathbf{r}; \tmmathbf{r}_s ; k) = \frac{E_{\RED{\ell}}
	(\tmmathbf{r}; \tmmathbf{r}_s ; \tmmathbf{\xi}^{\star}) \det (\Psi^{\star})
	e^{- i\tmmathbf{m} \cdot \tmmathbf{\xi}^{\star}}}{4 \pi^2} \mathcal{I}, &
\quad \text{where} \quad   & \mathcal{I}= \iint_{\tmmathbf{\Pi}} \frac{e^{-
		i\tmmathbf{\alpha} \cdot \tmmathbf{\zeta}}}{{(k^{\star} + \hat{k})^2}  -
	(k^{\star} + \Omega (\tmmathbf{\zeta}))^2} \mathd \tmmathbf{\zeta},
\end{eqnarray*}
%\end{widetext}
is the {\tmem{canonical integral}} associated to the degeneracy point
$(\tmmathbf{\xi}^{\star}, k^{\star})$. Here $\tmmathbf{\alpha}=
(\Psi^{\star})^{\transpose} \tmmathbf{m}$, $\det (\Psi^{\star})$ is the determinant of
the Jacobian arising from the change of variable, and $\tmmathbf{\Pi}$ is the
image of $\tmmathbf{\Gamma}$ in the new variable $\tmmathbf{\zeta}$. Because
$k^{\star} > 0$, this local integral can be further simplified to
\begin{eqnarray*}
\mathcal{I} & \approx & \frac{1}{2 k^{\star}} \iint_{\tmmathbf{\Pi}}
\frac{e^{- i\tmmathbf{\alpha} \cdot \tmmathbf{\zeta}}}{\hat{k} - \Omega_0
(\tmmathbf{\zeta})} \mathd \tmmathbf{\zeta},
\end{eqnarray*}
where $\Omega_0 (\tmmathbf{\zeta})$ is a local (Taylor-like) approximation to
$\Omega (\tmmathbf{\zeta})$ as $\tmmathbf{\zeta} \rightarrow \tmmathbf{0}$.
The integrand is not periodic anymore, and is only singular in the
neighbourhood of $\tmmathbf{\zeta}=\tmmathbf{0}$. We can then expand the
finite surface of integration $\tmmathbf{\Pi}$ to a doubly infinite one,
denoted $\tmmathbf{\Pi}^{\infty}$, that is almost $\mathbb{R}^2$ everywhere
apart from the neighbourhood of $\tmmathbf{0}$ and on which the integrand is
exponentially decaying as $| \tmmathbf{\alpha} | \rightarrow \infty$ (or,
equivalently, $| \tmmathbf{m} | \rightarrow \infty$). Therefore, the canonical
integral can be approximated by
\begin{eqnarray}
\mathcal{I} & \approx & \frac{1}{2 k^{\star}}
\iint_{\tmmathbf{\Pi}^{\infty}} \frac{e^{- i\tmmathbf{\alpha} \cdot
	\tmmathbf{\zeta}}}{\hat{k} - \Omega_0 (\tmmathbf{\zeta})} \mathd
\tmmathbf{\zeta},  \label{eq:approx-Omega-kstar-non-zero}
\end{eqnarray}
where $\tmmathbf{\Pi}^{\infty}$ can be seen as a small deformation of
$\mathbb{R}^2$ that bypasses the singularity at
$\tmmathbf{\zeta}=\tmmathbf{0}$. We will consider three cases specifically:
local minimums and local maximums (shrinking loops, see Figure
\ref{fig:local-extrememums}), as well as hyperbolic crossings (see Figure
\ref{fig:hyperbolic-crossing}).
%\begin{widetext}
%\end{widetext}
%\tmcolor{red}{{\tmstrong{Questions:}} Is it possible to have a quadratic
%	behaviour at $k^{\star} = 0$? Does it still count as a simple eigenvalue? Not
%	with Helmholtz, but yes with bi-harmonic plate equation., not sure if we
%	should mention it. Also I have a feeling that whenever the dispersion diagram
%	goes through $k^{\star} = 0$, it is always as a multiple eigenvalue of some
%	sort so not sure how to talk about $k^{\star} = 0$ here. }

\subsubsection{Local maximums and minimums }\label{sec:local-max}

Because the eigenvalues are smooth away from double eigenvalues, a local
maximum is required to be smooth and is usually quadratic. The function
$\Omega$ can then be typically assumed to behave like
\begin{eqnarray*}
\Omega (\tmmathbf{\zeta}) & \underset{\tmmathbf{\zeta} \rightarrow
0}{\approx} & \Omega_0 (\tmmathbf{\zeta}) = - \GREEN{\tfrac{1}{\Lambda}}(\zeta_1^2 + \zeta_2^2) .
\end{eqnarray*}
\GREEN{for some $\Lambda>0$.} Using (\ref{eq:approx-Omega-kstar-non-zero}) we have
\begin{eqnarray*}
\mathcal{I} & \approx & \mathcal{I}_{\max} = \frac{\GREEN{\Lambda}}{2 k^{\star}}
\iint_{\tmmathbf{\Pi}^{\infty}} \frac{e^{- i\tmmathbf{\alpha} \cdot
	\tmmathbf{\zeta}}}{\GREEN{\Lambda}\hat{k} + \zeta_1^2 + \zeta_2^2} \mathd \tmmathbf{\zeta}.
\end{eqnarray*}
%As shown in Appendix \ref{app:canonical-integrals}, this integral can be taken
%exactly and given in (\ref{eqapp:Imax}). 
%Note that by definition of
%$\tmmathbf{\alpha}$, $N$ and $\tmmathbf{m}$, there exists a constant
%$\mathfrak{C}> 0$ such that
%\begin{eqnarray*}
%	\sqrt{\alpha_1^2 + \alpha_2^2} & = & \mathfrak{C}N,
%\end{eqnarray*}
%where the constant $\mathfrak{C}> 0$ is defined by $\mathfrak{C}= |
%\tilde{\tmmathbf{\alpha}} | = | (\Psi^{\star})^{\transpose} \tilde{\tmmathbf{m}} |$.
\RED{Using the bypass of Figure~\ref{fig:local-extrememums}, left,} this integral can be evaluated:
\begin{eqnarray*}
\mathcal{I}_{\max} & = & \left\{ \begin{array}{ccc}
\frac{\GREEN{\Lambda}\pi}{k^{\star}} K_0 \left( |\tilde \bdalpha|N \sqrt{\GREEN{\Lambda}\hat{k}} \right) &
\tmop{if} & \hat{k} > 0\\
\RED{-}\frac{i\GREEN{\Lambda} \pi^2}{2 k^{\star}} H_0^{(2)} \left(|\tilde \bdalpha|N \sqrt{-
	\GREEN{\Lambda}\hat{k}} \right) & \tmop{if} & \hat{k} < 0
\end{array} \right. ,
\end{eqnarray*}
where
\begin{equation}
\tilde{\tmmathbf{\alpha}} = \tmmathbf{\alpha}/N , \quad |\tilde \bdalpha| = \sqrt{\tilde{\alpha}_1^2 + \tilde{\alpha}_1^2}.
\label{eq:tilde_alpha}
\end{equation} 

Using the usual far-field approximations for the modified Bessel function and
the Hankel functions, we recover that, as $N \sqrt{| \hat{k} |} \rightarrow
\infty$, the integral $\mathcal{I}_{\max}$ is exponentially decaying if
$\hat{k} > 0$ (no asymptotic contribution), and oscillatory and algebraically
decaying if $\hat{k} < 0$ (wave component asymptotic contribution), as one
might expect from a local maximum of a dispersion diagram. This far-field
behaviour is consistent with the fact that for $\hat{k} < 0$ we have a
$\tmop{SoS}$ for all $\tilde{\tmmathbf{m}}$, while for $\hat{k} > 0$, no SoS
are present. \RED{Note that if $k$ were to have a small imaginary part $\varkappa>0$, then the imaginary part of  $\sqrt{-\GREEN{\Lambda}\hat{k}}$ would be negative and $H_0^{(2)} \left(|\tilde \bdalpha|N \sqrt{-
\GREEN{\Lambda}\hat{k}} \right)$ would indeed decay exponentially to zero as $N\to\infty$, which is consistent with the limiting absorption principle.}

Using the near-field approximations for the modified Bessel function and the
Hankel functions, we find that, as $N \sqrt{| \hat{k} |} \rightarrow 0$,
$\mathcal{I}_{\max}$ blows up like $\log \left( |\tilde\bdalpha| N \sqrt{| \hat{k}
|} \right)$, highlighting a resonance phenomena, as expected.

Similarly, a local minimum can be described by 
\begin{eqnarray*}
\Omega (\tmmathbf{\zeta}) & \underset{\tmmathbf{\zeta} \rightarrow
0}{\approx} & \Omega_0 (\tmmathbf{\zeta}) = \GREEN{\tfrac{1}{\Lambda}}(\zeta_1^2 + \zeta_2^2),
\end{eqnarray*}
\GREEN{for some $\Lambda>0$}, that leads to 
\begin{eqnarray*}
\mathcal{I}_{\min} & = & \left\{ \begin{array}{ccc}
\frac{i\GREEN{\Lambda} \pi^2}{2 k^{\star}} H_0^{(\RED{1})} \left( |\tilde\alpha| N \sqrt{\GREEN{\Lambda}\hat{k}}
\right) & \tmop{if} & \hat{k} > 0\\
\frac{\GREEN{\Lambda}\pi}{k^{\star}} K_0 \left( |\tilde\alpha| N \sqrt{- \GREEN{\Lambda}\hat{k}} \right) &
\tmop{if} & \hat{k} < 0
\end{array} \right. .
\end{eqnarray*}

Acoustic waves or phonons in gases, liquids, and solids with a  parabolic maximum and minimum are referred to as maxons and rotons respectively \cite{Chen2021,Santos2003}. 

\begin{exa}
Consider the Green's function of a discrete lattice for some $k$ close to $k^\star=2\sqrt{2}$. Then, in the vicinity of $\bdxi^\star=(\pi,\pi)$ (\ref{eq:discrete-explicit-expr}) can be approximated as
\begin{equation}
F \approx \frac{1}{\RED{2 k^\star\hat{k}}+\zeta_1^2+\zeta_2^2},\quad  \zeta_n = \RED{\xi_n-\pi}, \quad \RED{n\in \{1,2\}} 
\end{equation}
i.e.\ it has a local maximum at $(\pi,\pi)$. %\RED{Should I make an explicit link with the notations of the general case, i.e\ say that here we have $\Lambda=2k^\star$, $\det(\Psi^\star)=1$, $E_\ell=1$ and refer to Remark \ref{rem:discrete-case-contiuous}. Also I prefer $\zeta_n=\xi_n-\pi$.}  
Maxons can also be observed at the points $(-\pi,\pi)$, $(\pi,-\pi)$, $(-\pi,-\pi)$. %Analogously, $F$ has a roton at the origin where it is described by $(\ref{F_discr_small_k})$.  
\end{exa}

%\subsubsection{Local minimum (roton)}\label{sec:local-min}
%Because the eigenvalues are smooth away from double eigenvalues, such a local
%minimum is also required to be smooth. The function $\Omega$ can then be
%typically assumed to behaved like
%\begin{eqnarray*}
%	\Omega (\tmmathbf{\zeta}) & \underset{\tmmathbf{\zeta} \rightarrow
%		0}{\approx} & \Omega_0 (\tmmathbf{\zeta}) = (\zeta_1^2 + \zeta_2^2) .
%\end{eqnarray*}
%Using (\ref{eq:approx-Omega-kstar-non-zero}) we have
%\begin{eqnarray*}
%	\mathcal{I} & \approx & \mathcal{I}_{\min} = \frac{1}{2 k^{\star}}
%	\iint_{\mathbb{R}^2} \frac{e^{- i\tmmathbf{\alpha} \cdot
%			\tmmathbf{\zeta}}}{\hat{k} - (\zeta_1^2 + \zeta_2^2)} \mathd
%	\tmmathbf{\zeta}.
%\end{eqnarray*}
%As shown in Appendix \ref{app:canonical-integrals}, this integral can be taken
%exactly and given in (\ref{eqapp:Imin}) and reads
%\begin{eqnarray*}
%	\mathcal{I}_{\min} & = & \left\{ \begin{array}{ccc}
%		\frac{i \pi^2}{2 k^{\star}} H_0^{(1)} \left( -\mathfrak{C}N \sqrt{\hat{k}}
%		\right) & \tmop{if} & \hat{k} > 0\\
%		\frac{\pi}{k^{\star}} K_0 \left( \mathfrak{C}N \sqrt{- \hat{k}} \right) &
%		\tmop{if} & \hat{k} < 0
%	\end{array} \right.
%\end{eqnarray*}
%With the same reasoning as above we see that as $N \sqrt{| \hat{k} |}
%\rightarrow \infty$, $\mathcal{I}_{\min}$ is exponentially decaying if
%$\hat{k} < 0$, and represent a wave component if $\hat{k} > 0$. Moreover, as
%$N \sqrt{| \hat{k} |} \rightarrow 0$, $\mathcal{I}_{\min}$ also blows up
%logarithmically, highlighting a resonance phenomena.

\subsubsection{Crossing with rebuilding: a seeming
contradiction}\label{sec:hyp-degeneracy}

\begin{figure}[h]
\centering{
\includegraphics[width=0.7\textwidth]{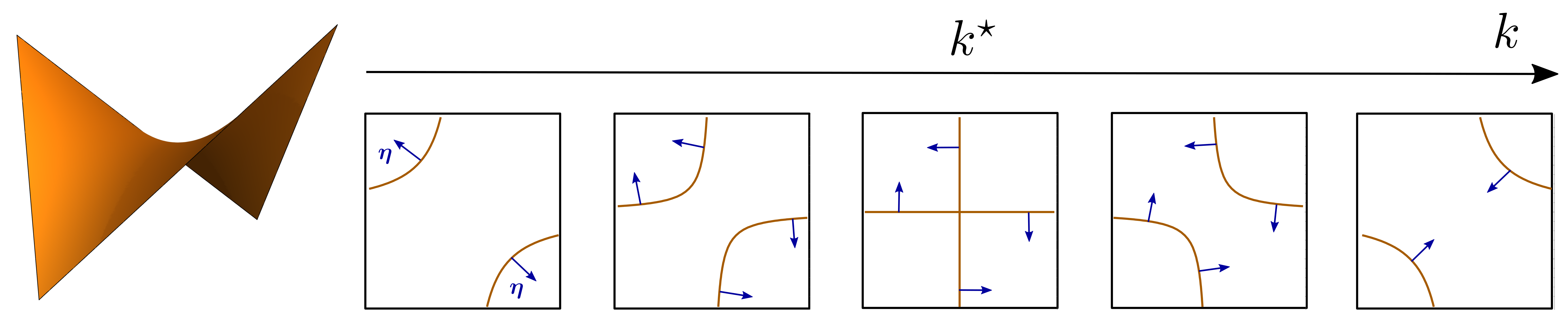}}
\caption{Typical simple eigenvalue hyperbolic degeneracy, together with the
orientation of the indentation vector $\tmmathbf{\eta}$: example of crossing
and rebuilding.}
\label{fig:hyperbolic-crossing}
\end{figure}

Let us now consider the case when the function $\Omega (\tmmathbf{\zeta})$ can be
approximated by
\begin{eqnarray}
\Omega (\tmmathbf{\zeta}) & \underset{\tmmathbf{\zeta} \rightarrow
0}{\approx} & \Omega_0 (\tmmathbf{\zeta}) = \zeta_1 \zeta_2 . 
\label{eq:Omega-hyp}
\end{eqnarray}
The evolution of the dispersion diagram is illustrated in Figure
\ref{fig:hyperbolic-crossing}. If $\hat{k} > 0$ or $\hat{k} < 0$, it is
possible to define a valid indentation of the surface of integration given by
the orientation of the vector field $\tmmathbf{\eta}$. However, when $\hat{k}
= 0$ (i.e.\ $k = k^{\star}$), we seem to reach a contradiction. Indeed if we
enforce continuity of the bypass as $\hat{k} \rightarrow 0$, as it should be,
we obtain the middle configuration of Figure \ref{fig:hyperbolic-crossing}. In
this case, we have two clear singularity components: the vertical line
$\zeta_1 = 0$ and the horizontal line $\zeta_2 = 0$. Following the bypass
rules of {\cite{Part6A}} that we reminded in section
\ref{sec:deformation-process}, the bypass should be continuous along each of
these lines. But we see in Figure \ref{fig:hyperbolic-crossing} (middle graph)
that it changes at the origin when $k = k^{\star}$ and therefore violates the
bypass rules. This means that for this crossing, we cannot define a proper
indented surface of integration, and the methods of {\cite{Part6A}} and
{\cite{Ice2021}} cannot be used to describe the wave field resulting from this
crossing. In other words, it is impossible to make sense of the integral
(\ref{eq:inverse-BF-transform}) for this value of $k$, as no admissible
deformation of the surface of integration can be made. This is a well-known
phenomena for lattices and metamaterials. For example, in
{\cite{Kapanadze2021}} it is stated that one cannot formulate a radiation
condition for the corresponding stationary problem. %Such degeneracy occurs for
%the discrete example treated in Section \ref{sec:discrete-example} for
%$k^{\star} = 2$ and $\tmmathbf{\xi}^{\star} = (\pi, 0)$.

%\paragraph{Approximation when $\hat{k} \neq 0$}
Let $k$ be close to $k^{\star}$. In this case, the canonical integral to consider
is of the type
\begin{eqnarray}
\mathcal{I} & \approx & \mathcal{I}_{\tmop{hyp}} = \frac{1}{2 k^{\star}}
\iint_{\tmmathbf{\Pi}^{\infty}} \frac{e^{- i\tmmathbf{\alpha} \cdot
	\tmmathbf{\zeta}}}{\hat{k} - \zeta_1 \zeta_2} \mathd \tmmathbf{\zeta}.
\label{eq:Ihyp-def}
\end{eqnarray}
Let us for simplicity deal with $\hat{k} > 0$, the case $\hat{k} < 0$ can be
dealt with in a very similar manner. This integral is somewhat standard and can be
taken exactly to give
%\begin{widetext}
\begin{eqnarray}
\mathcal{I}_{\tmop{hyp}} & \RED{=} &   \left\{ 
\begin{array}{ccccc}
	\frac{\pi^2}{k^{\star}} H_0^{(\RED{1})} \left( 2 \sqrt{\tilde{\alpha}_1
		\tilde{\alpha}_2} N \sqrt{\hat{k}} \right) & \tmop{if} & \alpha_1 > 0 &
	\infixand & \alpha_2 > 0\\
	\frac{- 2 i \pi}{k^{\star}} K_0 \left( 2 \sqrt{- \tilde{\alpha}_1
		\tilde{\alpha}_2} N \sqrt{\hat{k}} \right) & \tmop{if} & \alpha_1 > 0 &
	\infixand & \alpha_2 < 0\\
	\frac{- 2 i \pi}{k^{\star}} K_0 \left( 2 \sqrt{- \tilde{\alpha}_1
		\tilde{\alpha}_2} N \sqrt{\hat{k}} \right) & \tmop{if} & \alpha_1 < 0 &
	\infixand & \alpha_2 > 0\\
	\frac{\pi^2}{k^{\star}} H_0^{(\RED{1})} \left( 2 \sqrt{\tilde{\alpha}_1
		\tilde{\alpha}_2} N \sqrt{\hat{k}} \right) & \tmop{if} & \alpha_1 < 0 &
	\infixand & \alpha_2 < 0
\end{array} \right.,  \label{eq:Ihypfinal}
\end{eqnarray}
%\end{widetext}
%where $\tilde{\tmmathbf{{\alpha}}}=(\Psi^\star)^{\transpose} \tilde{\tmmathbf{m}}$ and
%$\tmmathbf{\alpha}= N \tilde{\tmmathbf{\alpha}}$. 
This represents a
propagating wave in the quadrant for which $\alpha_1$ and $\alpha_2$ have the
same sign, and no wave component otherwise. This is expected and can be
recovered by studying the SoS of this problem. For $\hat{k} < 0$, the `propagating' and
`evanescent' quadrants are swapped. Note that only one branch of the hyperbola
contains an SoS when we are in a propagating quadrant. Using the near-field asymptotics of the Hankel and modified Bessel functions, one can see that the integral has a logarithmic singularity as $\hat k \to 0$ which is a consequence of a pinching of the integration surface at the crossing point of the dispersion diagram. %\RED{Shall I say somewhere that for $\hat{k}<0$, we'd get $H_0^{(2)}$ ?}

\begin{exa}
Consider again the Green's function of the discrete lattice. Let $k$ be close to $k^\star=2$. Then, in the vicinity of the crossing point $\bdxi^\star = \RED{(0,\pi)}$ the dispersion diagram can be approximated as
\begin{equation}
F\approx\frac{\RED{1/4}}{ \RED{\hat{k}-} \zeta_1\zeta_2}, \quad \zeta_1 = \RED{\frac{\xi_1-(\xi_2-\pi)}{2}},\quad \zeta_2 = \RED{\frac{\xi_1+\xi_2-\pi}{2}},
\end{equation}
which leads to an integral of the type (\ref{eq:Ihyp-def}). While the stationary problem for (\ref{eq:discrete-governing}) does not make sense for $k=k^\star=2$, the physical question of what will happen if the lattice is excited exactly at this resonance frequency remains valid. We study the corresponding time domain problem in Appendix~\ref{App:E} for the sake of completeness. 
\end{exa}

%\tmcolor{red}{{\tmstrong{AS:}} Why can't we do as before and just show that
%	the solution for $\hat{k} \neq 0$ blows up logarithmically as $\hat{k}
%	\rightarrow 0$ and avoid doing all the time dependent stuff?}

\subsection{General procedure for a double eigenvalue}

Let us choose a pair $(\tmmathbf{\xi}^{\star}, k^{\star}) \in
\tmmathbf{\mathcal{B}} \times \mathbb{R}^+$ and an index $\ell \in
\mathbb{N}$, such that $\lambda_{\ell} (\tmmathbf{\xi}^{\star}) \geqslant 0$
is a double eigenvalue of $(\mathcal{L}, \tmop{dom}_{\tmmathbf{\xi}^{\star}})$
and $(k^{\star})^2 = \lambda_{\ell} (\tmmathbf{\xi}^{\star})$.

Let us denote by $V_{\ell}^{(1)} (\tmmathbf{r})$ and $V_{\ell}^{(2)}
(\tmmathbf{r})$ the two associated orthogonal eigenfunctions, normalised such
that $(V_{\ell}^{(\mathscr{m})}, V_{\ell}^{(\mathscr{k})}) =
\delta_{\mathscr{m}, \mathscr{k}}$ for $\mathscr{m}, \mathscr{k} \in \{ 1, 2
\}$. Note, that dependence on $\bdxi^\star$ is not indicated for $V_{\ell}^{(m)} (\tmmathbf{r})$. 

Our aim remains to provide a far-field approximation for $u$ that is valid
as $k$ approaches $k^{\star}$. Using perturbation analysis, we show in Appendix~\ref{app:A} that as $\tmmathbf{\xi} \rightarrow \tmmathbf{\xi}^{\star}$, \RED{$F(\tmmathbf{r}; \tmmathbf{r}_s ; \tmmathbf{\xi}; k)$} can be approximated as
\begin{eqnarray}
F (\tmmathbf{r}; \tmmathbf{r}_s ; \tmmathbf{\xi}; k) & \approx & a_1
V_{\ell}^{(1)} (\tmmathbf{r}) + a_2 V_{\ell}^{(2)} (\tmmathbf{r}), 
\label{eq:simple-approx-F-double-eig}
\end{eqnarray}
where the vector $\bda=(a_1,a_2)^\transpose$ has the following form:
\begin{eqnarray}
\label{eq:a_coef_cont}
\tmmathbf{a}& = & (\mathcal{Y}- \delta_{\lambda}^{(\ell)}{\rm I}_2)^{-1} \left(
\begin{array}{c}
(V_{\ell}^{(1)} (\tmmathbf{r}_s))^{\conj}\\
(V_{\ell}^{(2)} (\tmmathbf{r}_s))^{\conj}
\end{array} \right),
\label{eq:asandbs}	
\end{eqnarray}
where $\delta_\lambda^\ell\assign \lambda_{\ell} (\tmmathbf{\xi}^{\star}) -
k^2 = (k^{\star})^2 - k^2$ and ${\rm I}_2$ is the $2 \times 2$ identity matrix. The matrix $\mathcal{Y}$ is 
%\begin{widetext}
\begin{equation}
\mathcal{Y}= \delta_{\xi_1} \mathcal{Y}^{(L)} + \delta_{\xi_2}
\mathcal{Y}^{(B)}, \quad \mathcal{Y}^{(L)} = \left( \begin{array}{cc}
	Y_{1, 1}^{(L)} & Y_{1, 2}^{(L)}\\
	Y_{2, 1}^{(L)} & Y_{2, 2}^{(L)}
\end{array} \right), \quad \mathcal{Y}^{(B)} = \left( \begin{array}{cc}
	Y_{1, 1}^{(B)} & Y_{1, 2}^{(B)}\\
	Y_{2, 1}^{(B)} & Y_{2, 2}^{(B)}
\end{array} \right),
\end{equation}

where $\delta_{\xi_1} = (\xi_1-\xi^\star_1)$, $\delta_{\xi_2} = (\xi_2-\xi^\star_2)$, and the quantities 	$Y_{\mathscr{m}, \mathscr{k}}^{(L,B)}$ are given by 
%${\rm I}_2$ is the $2 \times 2$ identity matrix,
\begin{eqnarray*}
Y_{\mathscr{m}, \mathscr{k}}^{(L)} & = & i \int_{\partial S_0^L}
V_{\ell}^{(\mathscr{k})} (\tmmathbf{r}) (\tmmathbf{n} \cdot \nabla
V_{\ell}^{(\mathscr{m})})^{\conj} - (\tmmathbf{n} \cdot \nabla
V_{\ell}^{(\mathscr{k})} (\tmmathbf{r})) (V_{\ell}^{(\mathscr{m})})^{\conj}
\mathd s, \\
Y_{\mathscr{m}, \mathscr{k}}^{(B)} & = & i \int_{\partial S_0^B}
V_{\ell}^{(\mathscr{k})} (\tmmathbf{r}) (\tmmathbf{n} \cdot \nabla
V_{\ell}^{(\mathscr{m})})^{\conj} - (\tmmathbf{n} \cdot \nabla
V_{\ell}^{(\mathscr{k})} (\tmmathbf{r})) (V_{\ell}^{(\mathscr{m})})^{\conj}
\mathd s.
\end{eqnarray*}
%\end{widetext}

Consider the  eigenfunction $V_{\ell}
(\tmmathbf{r}; \tmmathbf{\xi})$ for $\tmmathbf{\xi}$ in the neighbourhood of $\tmmathbf{\xi}^\star$, that corresponds to the eigenvalue $\lambda_{\ell}(\tmmathbf{\xi})$. Similarly, it can be approximated as:
\begin{eqnarray}
V_{\ell} (\tmmathbf{r}; \tmmathbf{\xi}) & \approx & \tilde{a}_1
V_{\ell}^{(1)} (\tmmathbf{r}) + \tilde{a}_2 V_{\ell}^{(2)} (\tmmathbf{r}) 
\label{eq:double-eigenproj-Vl}
\end{eqnarray}
for some new vector $\tilde \bda =(\tilde a_{1},\tilde a_{2})^\transpose$ that \RED{is} given as solution of the following homogeneous system of linear equations:
\begin{equation}
(\mathcal{Y}- \delta_{\lambda}^{(\ell)} {\rm I}_2) \tmmathbf{\tilde a}=\tmmathbf{0}.
\end{equation}
The latter can only have a non-trivial solution
provided that
\begin{eqnarray}
\det (\mathcal{Y}- \delta_{\lambda}^{(\ell)} {\rm I}_2) & = & 0. 
\label{eq:detzerodisp}
\end{eqnarray}
This equation gives a local approximation to the dispersion diagram.

\subsubsection{Dirac conical points}\label{sec:Dirac}

\begin{figure}[h]
\centering{
\includegraphics[width=0.7\textwidth]{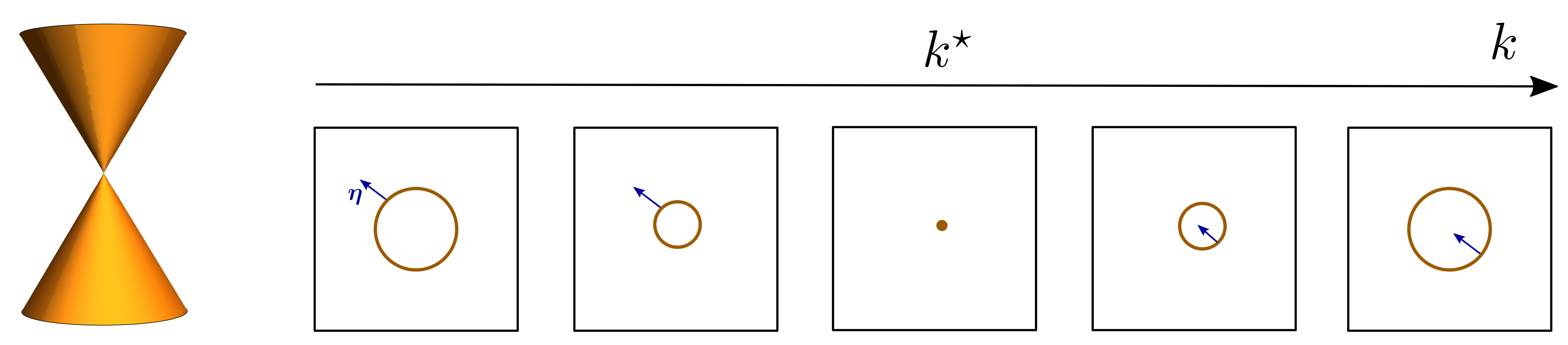}}
\caption{Typical conical Dirac point double eigenvalue degeneracy}
\label{fig:degeneracy-dirac}
\end{figure}

As we have done previously, let us consider a local change of variables
$\tmmathbf{\xi}= (\xi_1, \xi_2) \leftrightarrow \tmmathbf{\zeta}= (\zeta_1,
\zeta_2)$, such that, as $\tmmathbf{\xi} \rightarrow \tmmathbf{\xi}^{\star}$,
we have $\tmmathbf{\xi}-\tmmathbf{\xi}^{\star} \approx \Psi^{\star}
\tmmathbf{\zeta}$, where $\Psi^{\star}$ is the invertible Jacobian matrix
associated to the change of variable. We can then rewrite $\mathcal{Y}= \delta_{\xi_1} \mathcal{Y}^{(L)} + \delta_{\xi_\RED{2}}
\mathcal{Y}^{(B)}$ in the new coordinate system as
\begin{eqnarray*}
\mathcal{Y} & = & \zeta_1 \tilde{\mathcal{Y}}^{(1)} + \zeta_2
\tilde{\mathcal{Y}}^{(2)},
\end{eqnarray*}
where the matrices $\tilde{\mathcal{Y}}^{(1, 2)}$ can be written in terms of
$\mathcal{Y}^{(L, B)}$ and the entries of $\Psi^{\star}$. Let us assume that
we are in a situation where (see example below, and  Appendix~\ref{FEM_construction})
\begin{eqnarray}
\label{Yanzats}
\tilde{\mathcal{Y}}^{(1)} = \left( \begin{array}{cc}
q_1 & q^\conj_2\\
q_2 & - q_1
\end{array} \right) & \tmop{and} & \tilde{\mathcal{Y}}^{(2)} = \left(
\begin{array}{cc}
q_3 & {q}^\conj_4\\
q_4 & -q_3
\end{array} \right),
\end{eqnarray}
where $q_{j}$ are such that
\[
q_1^2 + |q_2|^2 = 1, \quad q_3^2 + |q_4|^2 = 1,\quad q_4{q}^\conj_2 + {q}^\conj_4q_2 + 2q_1q_3 = 0.
\]
%Remembering that $\delta_{\lambda}^{(\ell)} = (k^{\star})^2 - k^2$, 
Therefore, (\ref{eq:detzerodisp})
implies that the dispersion diagram is locally approximated by $(\delta_{\lambda}^{(\ell)})^2  =  \zeta_1^2 + \zeta_2^2$.
%\begin{eqnarray*}
%	(k^2 - (k^{\star})^2)^2 & = & \zeta_1^2 + \zeta_2^2 .
%\end{eqnarray*}
Then, from (\ref{eq:asandbs}) we get
%\begin{widetext}
\begin{equation}
\left( \begin{array}{c}
	a_1\\
	a_2
\end{array} \right)  =  \frac{1}{\zeta_1^2 + \zeta_2^2 - (\delta_{\lambda}^{(\ell)})^2} \left( \begin{array}{cc}
	\delta_{\lambda}^{(\ell)} + q_1 \zeta_1 + q_3 \zeta_2  & {q}^\conj_2 \zeta_1 + {q}^\conj_4\zeta_2\\
	q_2 \zeta_1 + q_4\zeta_2 & \delta_{\lambda}^{(\ell)} - q_1 \zeta_1 - q_3 \zeta_2
\end{array} \right) \left( \begin{array}{c}
	(V_{\ell}^{(1)} (\tmmathbf{r}_s))^{\conj}\\
	(V_{\ell}^{(2)} (\tmmathbf{r}_s))^{\conj}
\end{array} \right).  \label{eq:coeffa1a2zeta}
\end{equation}
%\end{widetext}
%
%\begin{equation}
%	\left( \begin{array}{c}
%		a_1\\
%		a_2
%	\end{array} \right)  =  \frac{1}{\zeta_1^2 + \zeta_2^2 - (k^2 -
%		(k^{\star})^2)^2} \left( \begin{array}{cc}
%		(k^{\star})^2 - k^2 + q_1 \zeta_1 + q_3 \zeta_2  & {q}^\conj_2 \zeta_1 + {q}^\conj_4\zeta_2\\
%		q_2 \zeta_1 + q_4\zeta_2 & (k^{\star})^2 - k^2 - q_1 \zeta_1 - q_3 \zeta_2
%	\end{array} \right) \left( \begin{array}{c}
%		(V_{\ell}^{(1)} (\tmmathbf{r}_s))^{\conj}\\
%		(V_{\ell}^{(2)} (\tmmathbf{r}_s))^{\conj}
%	\end{array} \right)  \label{eq:coeffa1a2zeta}
%\end{equation}
If $k^{\star} \neq 0$, upon introducing $\hat{k} = k - k^{\star}$, and remembering that $\delta_{\lambda}^{(\ell)} = (k^{\star})^2 - k^2$  the
dispersion diagram can then be locally approximated by $(2 k^{\star}
\hat{k})^2 = \zeta_1^2 + \zeta_2^2$, which is a circular cone (see Figure \ref{fig:degeneracy-dirac}) merging at
$\hat{k} = 0$, each part of which being defined by
\begin{eqnarray*}
\hat{k} & = & \pm \frac{1}{2 k^{\star}} \sqrt{\zeta_1^2 + \zeta_2^2}
\end{eqnarray*}

We now focus on finding the asymptotic component of $u$ resulting from such
degeneracy. Upon introducing $\rho = 2 k^{\star} \hat{k}$ \RED{and noting that $\delta_{\lambda}^{(\ell)}\approx-\rho$}, the equation
(\ref{eq:coeffa1a2zeta}) can be further simplified to 
%\RED{CAREFUL HERE: with this definition, we have $\delta_{\lambda}^{\ell}\approx-\rho$ so we should have $-$ signs in front of the $\rho$s in the eq below}
%\begin{widetext}
\begin{eqnarray}
\left( \begin{array}{c}
	a_1\\
	a_2
\end{array} \right) & \approx & \frac{1}{\zeta_1^2 + \zeta_2^2 - \rho^2}
\left( \begin{array}{cc}
	q_1 \zeta_1+q_3 \zeta_2 \RED{-} \rho & {q}^\conj_2 \zeta_1 + {q}^\conj_4\zeta_2\\
	q_2 \zeta_1 + q_4\zeta_2 & -q_1 \zeta_1-q_3 \zeta_2 \RED{-} \rho
\end{array} \right) \left( \begin{array}{c}
	(V_{\ell}^{(1)} (\tmmathbf{r}_s))^{\conj}\\
	(V_{\ell}^{(2)} (\tmmathbf{r}_s))^{\conj}
\end{array} \right). \label{eq:coeffa1a2zetasi,plenonwerokstqr}
\end{eqnarray}
Using the approximation (\ref{eq:simple-approx-F-double-eig}) we can write the
associated far-field wave component $u_{\tmop{loc}}$ as
\begin{eqnarray}\label{eq:Dirac-point-far-field}
u_{\tmop{loc}} (\tmmathbf{m}; \tmmathbf{r}; \tmmathbf{r}_s ; k) & = &
\frac{\det (\Psi^{\star}) e^{- i\tmmathbf{m} \cdot
		\tmmathbf{\xi}^{\star}}}{4 \pi^2} \times \left( \sum_{\mathscr{m},
	\mathscr{k}= 1}^2 \mathcal{I}_{\mathscr{m}\mathscr{k}}
V_{\ell}^{(\mathscr{m})} (\tmmathbf{r}) (V_{\ell}^{(\mathscr{k})}
(\tmmathbf{r}_s))^{\conj} \right)
\end{eqnarray}
where the integrals $\mathcal{I}_{\mathscr{m}\mathscr{k}}$ are given by
\begin{alignat*}{2}
\mathcal{I}_{11} &= \iint_{\tmmathbf{\Pi}^{\infty}} \frac{(q_1 \zeta_1+q_3 \zeta_2 \RED{-} \rho) e^{- i\tmmathbf{\alpha} \cdot \tmmathbf{\zeta}}}{\zeta_1^2 + \zeta_2^2
	- \rho^2} \mathd \tmmathbf{\zeta}, \qquad && \mathcal{I}_{12} =
\iint_{\tmmathbf{\Pi}^{\infty}} \frac{({q}^\conj_2 \zeta_1 + {q}^\conj_4\zeta_2) e^{-
		i\tmmathbf{\alpha} \cdot \tmmathbf{\zeta}}}{\zeta_1^2 + \zeta_2^2 - \rho^2}
\mathd \tmmathbf{\zeta},\\
\mathcal{I}_{21} &= \iint_{\tmmathbf{\Pi}^{\infty}} \frac{({q}_2 \zeta_1 + {q}_4\zeta_2) e^{- i\tmmathbf{\alpha} \cdot \tmmathbf{\zeta}}}{\zeta_1^2 +
	\zeta_2^2 - \rho^2} \mathd \tmmathbf{\zeta}, && \mathcal{I}_{22} =
-\iint_{\tmmathbf{\Pi}^{\infty}} \frac{(q_1 \zeta_1+q_3 \zeta_2 \RED{+} \rho) e^{-
		i\tmmathbf{\alpha} \cdot \tmmathbf{\zeta}}}{\zeta_1^2 + \zeta_2^2 - \rho^2}
\mathd \tmmathbf{\zeta}.
\end{alignat*}
%{\color{red} \textbf{AIK:} need to check the signs in those integrals, lots of changes from version 04 it seems! }{\color{green}Done, see the corrections in green starting from (\ref{eq:coeffa1a2zeta})} 
Computing the $\mathcal{I}_{\mathscr{m}\mathscr{k}}$ integrals can be reduced
to the computation of only three canonical integrals, $\mathcal{J}_0$,
$\mathcal{J}_1$ and $\mathcal{J}_2$ given by
\begin{equation}
\label{eq:canon_int_Dirac}
\mathcal{J}_0(\tmmathbf{\alpha},\rho)= \iint_{\tmmathbf{\Pi}^{\infty}} \frac{e^{-
		i\tmmathbf{\alpha} \cdot \tmmathbf{\zeta}}}{\zeta_1^2 + \zeta_2^2 - \rho^2}
\mathd \tmmathbf{\zeta}, \quad  \mathcal{J}_{1,2}(\tmmathbf{\alpha},\rho) = \iint_{\tmmathbf{\Pi}^{\infty}}
\frac{\zeta_{1,2} e^{- i\tmmathbf{\alpha} \cdot \tmmathbf{\zeta}}}{\zeta_1^2 +
	\zeta_2^2 - \rho^2} \mathd \tmmathbf{\zeta},  
\end{equation}
as indeed we have
\begin{eqnarray*}
\mathcal{I}_{11} = q_1 \mathcal{J}_1 +q_3 \mathcal{J}_2 \RED{-}\rho\mathcal{J}_0, \quad & \mathcal{I}_{12} =
{q}^\conj_2 \mathcal{J}_1 + {q}^\conj_4\mathcal{J}_2, \quad & \mathcal{I}_{21} = q_2 \mathcal{J}_1 +
{q}_4\mathcal{J}_2, \quad \mathcal{I}_{22} = -q_1 \mathcal{J}_1-q_3 \mathcal{J}_2\RED{-}\rho\mathcal{J}_0 .
\end{eqnarray*}

\RED{For $\hat{k}>0$,} these canonical integrals  are given by
\begin{eqnarray*}
\mathcal{J}_0 (\tmmathbf{\alpha},\rho) =i \pi^2 H_0^{(1)} (\rho |\bdalpha| ) &
\quad\text{ and }\quad & \mathcal{J}_{1, 2}(\tmmathbf{\alpha},\rho) =\frac{\pi^2 \rho{\alpha}_{1, 2}}{| {\tmmathbf{\alpha}} |}   H_1^{(1)}
(\rho |\bdalpha|),
\end{eqnarray*}
%\end{widetext}
\RED{while for $\hat{k}<0$ they will be written in terms of $H_0^{(2)}(-\rho|\bdalpha|)$ and $H_1^{(2)}
(-\rho |\bdalpha|)$.} Note that for fixed large $N$ and $\hat{k} \rightarrow 0$, $\mathcal{J}_0(\tilde\bdalpha N,2 k^{\star}
\hat{k})$
behaves like $\hat{k} \log (\hat{k} N)$, which tends to zero as $\hat{k}
\rightarrow 0$. Moreover the quantities $\mathcal{J}_{1, 2}(\tilde\bdalpha N,2 k^{\star}
\hat{k})$ behave like $1 /N$ and are therefore not zero in that limit. Note further that as $\hat{k} N
\rightarrow \infty$, all the terms, $\mathcal{J}_0(\tilde\bdalpha N,2 k^{\star}
\hat{k})$ and $\mathcal{J}_{1, 2}(\tilde\bdalpha N,2 k^{\star}
\hat{k})$
have a wave-like behaviour dictated by the far-field of the Hankel functions
and can be shown to behave like $e^{2 i k^{\star}
\hat{k} |\tilde\bdalpha|N} \sqrt{\hat{k} / N}$.
\begin{exa}
\label{ex:FEM-stuff}
Consider a finite element formulation for the phononic crystal problem with Neumann boundary conditions on circular scatterers. Let each cell be approximated by  linear triangular elements. The discrete field can be described by the set of values $u(\tmmathbf{m},j)$, where $\tmmathbf{m}$ is the number of the cell, and $j$ is the node index inside the cell, where $j$ is running through $\{1,\ldots, N_{\text{cell}}\}$, $N_{\text{cell}}$ being the number of nodes in the cell. The field satisfies a matrix equation of the form:
\begin{equation}
{\rm K U} - k^2{\rm M U} = {\rm E},
\end{equation}
where ${\rm K}$ and ${\rm M}$ are infinite stiffness and mass matrices, ${\rm U}$ is an infinite vector of the field values $u(\mathbf{m},j)$, and ${\rm E}$ is a forcing term. Due to the geometry of the problem, the matrices ${\rm K}$ and ${\rm M}$ are periodic in the sense that they are made up of a repetition of some $N_{\text{cell}}\times N_{\text{cell}}$ matrices ${\rm K_0}$ and ${\rm M_0}$. Let us take the Floquet-Bloch transform
\begin{equation}
{F}_j(\tmmathbf{\xi}) = \sum_{\tmmathbf{m} \in \mathbb{Z}^2}u(\tmmathbf{m},j)e^{i\tmmathbf{m}\cdot\tmmathbf{\xi}},
\end{equation}
and construct a $N_{\text{cell}}\times1$ vector ${\rm F}(\tmmathbf{\xi})$ whose components are denoted ${F}_j(\tmmathbf{\xi})$. As shown in Appendix~\ref{FEM_construction}, because of quasi-periodicity, all the information contained in ${\rm F}$ can be reduced to a smaller $\hat{N}\times1$ vector $\hat{F}$, that satisfies the following matrix equation:
\begin{equation}
\label{eq:discreteFloqueteq}
{ \rm \hat K}(\tmmathbf{\xi}){\rm \hat{F}}(\tmmathbf{\xi}) - k^2{ \rm \hat M}(\tmmathbf{\xi}) {\rm \hat{F}}(\tmmathbf{\xi}) = {\rm \hat E}(\tmmathbf{\xi}),
\end{equation}
where ${\rm\hat K}$  and ${\rm \hat M}$ are $\hat{N} \times \hat{N}$ mass and stiffness matrices for a unit cell with $\tmmathbf{\xi}$-quasiperiodic boundary conditions and ${\rm \hat E}$ is a forcing $\hat{N}\times1$ vector. 
%We discuss the structure of ${\rm\hat K}$, ${\rm \hat M}$ and ${\rm \hat F}$ in Appendix~\ref{FEM_construction} in details using an energy based approach. 
Since (\ref{eq:discreteFloqueteq}) is a matrix equation of finite dimension, it can be solved by inversion, leading to $\rm \hat{F}$ and then $\rm F$. The field can then be recovered using the inverse Floquet-Bloch transform:
\begin{equation}
\label{eq:discret_inverse_floquet}
u(\tmmathbf{m},j) = 	\frac{1}{4 \pi^2} \iint_{\tmmathbf{\mathcal{B}}} {F}_j (\tmmathbf{\xi}) e^{- i\tmmathbf{m} \cdot \tmmathbf{\xi}} \mathd \tmmathbf{\xi}.
\end{equation} 
%		Note also that the vector ${\rm \hat{F}}$ can be written as:
%		\begin{equation}
%			\label{eq:F_discrete_series}
%			{\rm \hat{F}}(\bdxi) = \sum_{j=1}^{\hat{N}} a_j {\hat{\rm V}}_{j}(\tmmathbf{\xi}), 
%		\end{equation}
%		where $a_j$ are some coefficients, and 

\RED{Let} ${\hat{\rm V}}_{j}(\tmmathbf{\xi})$ \RED{be} the eigenvector associated to the eigenvalue $\hat{\lambda}_j(\bdxi)$ of the generalized eigenvalue problem $\hat{\rm K}(\tmmathbf{\xi}) {\hat{\rm V}}= \lambda \hat{\rm M}(\tmmathbf{\xi}) {\hat{\rm V}} $, i.e.\ we have
\begin{equation}
\label{eq:eigenequation}
\hat{\rm K}(\tmmathbf{\xi}) {\hat{\rm V}}_{j}(\tmmathbf{\xi})= \hat{\lambda}_j(\bdxi)\hat{\rm M}(\tmmathbf{\xi}) {\hat{\rm V}}_{j}(\tmmathbf{\xi}).
\end{equation}
%	For simplicity, we assume that the eigenvalues are simple here.

Consider the mesh of equilateral triangles shown in Figure~\ref{fig:FEM_unit_cell} (left), and plot the first two eigenvalues as functions of $\tmmathbf{\xi}$. The corresponding dispersion surfaces are shown in  Figure~\ref{fig:FEM_unit_cell} (right). 
\begin{figure*}%[h]
\centering{
	\includegraphics[width=0.49\textwidth]{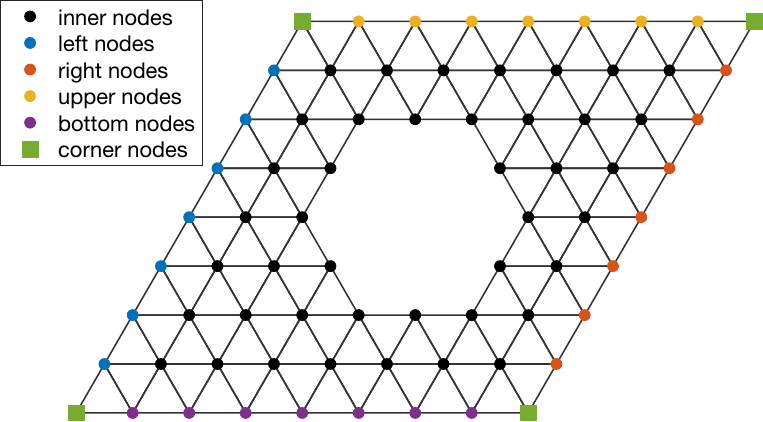}
	\includegraphics[width=0.49\textwidth]{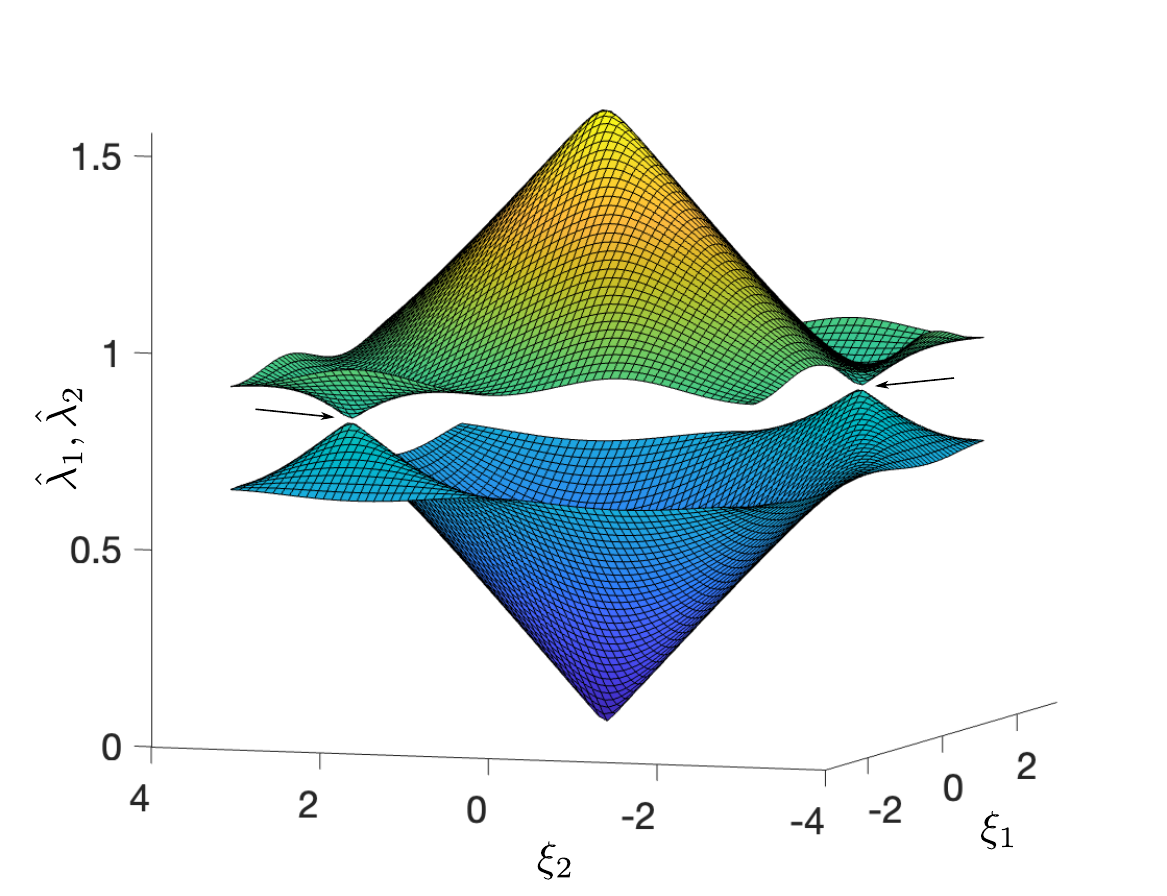}
	\caption{Geometry of the FEM problem (left), and dispersion surfaces of the first two eigenvalues~(right)}
	\label{fig:FEM_unit_cell}
}	
\end{figure*}
%\begin{figure}[h!]
%    \centering{
%	\includegraphics[width=0.7\textwidth]{figures/FEM_Disp_Surfaces.eps}
%	\caption{Dispersion surfaces of the problem (first two eigenvalues)}
%        \label{fig:FEM_Disp_Surfaces}
%        }	
%\end{figure}
As indicated by arrows on Figure~\ref{fig:FEM_unit_cell}, there are two points corresponding to a double eigenvalue. A result similar to  (\ref{eq:simple-approx-F-double-eig}) can be derived for the FEM problem, which is done in Appendix~\ref{FEM_construction}, and  it is shown there that the matrices (their analogues for the discrete problem to be precise)  $\tilde{\mathcal{Y}}^{(1,2)}$ have the same structure (\ref{Yanzats}). Hence, those double eigenvalue points are indeed Dirac conical points.   
\end{exa}

\section{Conclusion and perspectives}

We have developed a general procedure to study double integrals of Floquet-Bloch type occuring when studying wave propagation in complex structures. In particular, we gave explicit formulae for the Green's functions associated to these structures and estimated their far-field behaviour. To do this we used results of multidimensional complex analysis to efficiently deform the surface of integration in $\mathbb{C}^2$. The general framework developed was shown to be relevant for practical examples such as discrete lattices and continuous phononic crystals. We also used our approach to shed some light on degeneracies, including Dirac conical points, maxons, rotons and hyperbolic degeneracies. Here we have insisted on obtaining explicit asymptotic formulae, but we note that the deformations proposed can also be used for an efficient numerical evaluation of such highly-oscillatory double integrals. Extensions to the present work, will be to describe optimal surface deformations and also to extend our approach to three-dimensional wave propagation.

%\begin{acknowledgements}
\paragraph{Acknowledgement and Funding}
RCA and AVS would like to thank the
Isaac Newton Institute (INI) for Mathematical Sciences for support and hospitality
during the programme {\tmem{Mathematical theory and applications of multiple
	wave scattering}} when work on this paper was undertaken. RCA would particularly like to acknowledge interesting discussions on the topic of this paper with Profs Sonia Fliss and Alex Barnett. This programme was
supported by EPSRC grant EP/R014604/1. AVS thanks the UK Acoustics Network (UKAN) for financially supporting his visit to the INI and RCA acknowledges support from the EPSRC grant EP/W018381/1. AIK has been supported by the Royal Society via Dr Kisil's Dorothy Hodgkin Research Fellowship. 
%\end{acknowledgements}

\appendix
\section{\RED{How to define $F$ when $\bdxi$ is complex for the phononic crystal of section \ref{sec:continuous-example}?}}\label{sec:intro-psi-xi}

Let us consider $\tmmathbf{\xi} \in \mathbb{C}^2 \setminus \mathbb{R}^2$. We
can define the exact same operator $(\mathcal{L},
\tmop{dom}_{\tmmathbf{\xi}})$, following the same technique as in section
\ref{sec:self-adjoint}, but keeping in mind that this time we have $(e^{- i
\xi_{1, 2}})^{\conj} \neq e^{+ i \xi_{1, 2}}$, we can show that $(\mathcal{L}f,
g) \neq (f, \mathcal{L}g)$. Hence the operator is non-symmetric and cannot be
self-adjoint anymore, so it is not quite clear how to find $F$ in that case.

\begin{figure}[h]
\centering{\includegraphics[width=0.3\textwidth]{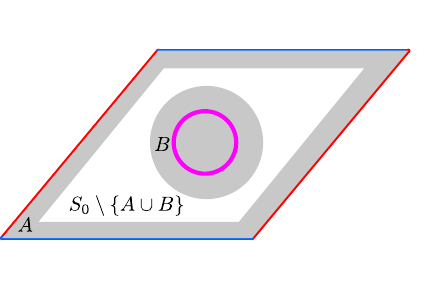}}
\caption{Illustration of the auxiliary regions $A$ and $B$ within the unit
cell.}
\label{fig:def-A-B-cell}
\end{figure}

Remember that our aim remains to show that $F$, as a function of
$\tmmathbf{\xi}$, can be analytically continued in a $\mathbb{C}^2$
neighbourhood of $\tmmathbf{\mathcal{B}}$. The easiest way to prove this is to
introduce a new function and a new operator. Start by decomposing $S_0$ into
three distinct regions $A$, $B$ and $S_0 \setminus \{ A \cup B \}$, where $A$
is a neighbourhood of $\partial S^T_0 \cup \partial S^B_0 \cup \partial S^R_0
\cup \partial S^L_0$, and $B$ is a neighbourhood of $\partial \Omega$, as
illustrated in Figure \ref{fig:def-A-B-cell}. For these regions, define the
function $\chi (\tmmathbf{r})$ that is smooth on $S_0$, strictly positive on
$S_0 \setminus \{ A \cup B \}$, and satisfies
\begin{eqnarray*}
\chi (\tmmathbf{r}) & = & \left\{ \begin{array}{ccc}
1 & \tmop{if} & \tmmathbf{r} \in A\\
0 & \tmop{if} & \tmmathbf{r} \in B
\end{array} \right. \cdot
\end{eqnarray*}

Consider now $f$ in $\tmop{dom}_{\tmmathbf{\xi}}$, that is $f$ is
$\tmmathbf{\xi}$-quasiperiodic w.r.t.\ $\Lambda$ and satisfies BC on $\partial
\Omega$, and introduce the new function $\phi$ defined by
\begin{align}
\phi (\tmmathbf{r}) = f (\tmmathbf{r}) / \psi_{\tmmathbf{\xi}}
(\tmmathbf{r}), & \text{ where }  \psi_{\tmmathbf{\xi}} (\tmmathbf{r}) = e^{+
i\tmmathbf{\xi} \cdot \Lambda^{- 1} \tmmathbf{r} \chi (\tmmathbf{r})} .  \label{eq:def-Psi-Xi}
\end{align}
This is an analytic function of $\tmmathbf{\xi}$ that is never zero on $S_0$.
%Note that $\psi_{\tmmathbf{\xi}}$ varies analytically w.r.t.\ $\tmmathbf{\xi}$
%and that as $\tmmathbf{\xi} \rightarrow \tmmathbf{0}$, $\psi_{\tmmathbf{\xi}}
%(\tmmathbf{r}) \rightarrow 1$. So we can assume that in a neighbourhood of
%$\tmmathbf{\xi}=\tmmathbf{0}$, $\psi_{\tmmathbf{\xi}} (\tmmathbf{r})$ is never
%zero on $S_0$. 
With this definition we have that
\begin{eqnarray*}
\phi (\tmmathbf{r}) & = & \left\{ \begin{array}{ccc}
e^{- i\tmmathbf{\xi} \cdot \Lambda^{- 1} \tmmathbf{r}} f (\tmmathbf{r}) &
\tmop{if} & \tmmathbf{r} \in A\\
f (\tmmathbf{r}) & \tmop{if} & \tmmathbf{r} \in B
\end{array} \right.,
\end{eqnarray*}
which implies that $\phi$ gives rise to a $\Lambda$-periodic function and
satisfies the same BC as $f$ on $\partial \Omega$. In other words, if $f \in
\tmop{dom}_{\tmmathbf{\xi}}$, then $\phi \in \tmop{dom}_{\tmmathbf{0}}$. By
direct calculations, it can be shown that
\begin{eqnarray*}
\mathcal{L}f & = & \phi \mathcal{L} \psi_{\tmmathbf{\xi}} - 2 \left(
\frac{\partial \psi_{\tmmathbf{\xi}}}{\partial x_1} \frac{\partial
\phi}{\partial x_1} + \frac{\partial \psi_{\tmmathbf{\xi}}}{\partial x_2}
\frac{\partial \phi}{\partial x_2} \right) + \psi_{\tmmathbf{\xi}}
\mathcal{L} \phi .
\end{eqnarray*}
This suggest introducing a new differential operator
$\mathcal{L}_{\tmmathbf{\xi}}$ defined by
\begin{eqnarray*}
\mathcal{L}_{\tmmathbf{\xi}} \phi & \equiv &
\frac{1}{\psi_{\tmmathbf{\xi}}} \mathcal{L}f =
\frac{\phi}{\psi_{\tmmathbf{\xi}}} \mathcal{L} \psi_{\tmmathbf{\xi}} -
\frac{2}{\psi_{\tmmathbf{\xi}}} \left( \frac{\partial
\psi_{\tmmathbf{\xi}}}{\partial x_1} \frac{\partial \phi}{\partial x_1} +
\frac{\partial \psi_{\tmmathbf{\xi}}}{\partial x_2} \frac{\partial
\phi}{\partial x_2} \right) +\mathcal{L} \phi
\end{eqnarray*}
the dependence of $\mathcal{L}_{\tmmathbf{\xi}}$ on $\tmmathbf{\xi}$ is
explicit, and given the definition of $\psi_{\tmmathbf{\xi}}$, it is analytic.
With this definition, we have $\mathcal{L}f = \psi_{\tmmathbf{\xi}}
\mathcal{L}_{\tmmathbf{\xi}} \phi$, and therefore $\mathcal{L}f = \lambda f
\Leftrightarrow \mathcal{L}_{\tmmathbf{\xi}} \phi = \lambda \phi$. In other
words, solving the eigenvalue problem associated to the operator
$(\mathcal{L}, \tmop{dom}_{\tmmathbf{\xi}})$ is equivalent to solving the
eigenvalue problem associated to the operator $(\mathcal{L}_{\tmmathbf{\xi}},
\tmop{dom}_{\tmmathbf{0}})$, where $\tmop{dom}_{\tmmathbf{0}}$ is the set of
$\Lambda$-periodic functions that satisfy the correct BC on $\partial \Omega$.
The advantage of this new formulation is that all the operators
$\mathcal{L}_{\tmmathbf{\xi}}$ have the same domain, and that they depend
analytically on $\tmmathbf{\xi}$. As such the family of operators
$(\mathcal{L}_{\tmmathbf{\xi}}, \tmop{dom}_{\tmmathbf{0}})$ is a holomorphic
family of operators or type (A) as treated in
{\cite{Kato1976}}(VII.{\textsection}2). Moreover, we know that for real
$\tmmathbf{\xi}$, the operators $(\mathcal{L}, \tmop{dom}_{\tmmathbf{\xi}})$
are self-adjoint and so this is also the case for
$(\mathcal{L}_{\tmmathbf{\xi}}, \tmop{dom}_{\tmmathbf{0}})$. Hence the family
of operators $(\mathcal{L}_{\tmmathbf{\xi}}, \tmop{dom}_{\tmmathbf{0}})$ is a
self-adjoint real holomorphic family of operators, see
{\cite{Kato1976}}(VII.{\textsection}3.2).

A direct consequence of this is that the eigenvalues $\lambda_j
(\tmmathbf{\xi})$ with multiplicity 1 are real-analytic functions of
$\tmmathbf{\xi}$, and so are the eigenfunctions. Of course, one needs to be
careful when saying that the eigenfunctions are analytic in $\tmmathbf{\xi}$
since they are defined up to an arbitrary multiplication constant that may
depend non-analytically on $\tmmathbf{\xi}$. However the quantity
$\frac{V_j (\tmmathbf{r}_s ; -\tmmathbf{\xi}) V_j (\tmmathbf{r};
\tmmathbf{\xi})}{\langle V_j (\cdot ; -\tmmathbf{\xi}), V_j (. ;
\tmmathbf{\xi}) \rangle}$ is completely independent of the choice of
normalisation, so this term is clearly analytic in $\tmmathbf{\xi}$.

This result stops being valid for an eigenvalue with multiplicity $> 1$, see
e.g.\ {\cite{Wilcox1978}} and {\cite{Kato1976}}(II.{\textsection}7.1). In
what follows, unless specified otherwise, we will assume that all eigenvalues
have multiplicity one. Though we will also discuss specific cases of multiple
eigenvalues. Because of this real-analyticity property, we automatically
deduce that eigenvalues and eigenprojections can be analytically continued
into the complex plane in a neighbourhood of $\tmmathbf{\mathcal{B}}$. As a
result, the formula (\ref{eq:continuous-explicit-F-no-conjugate}), and
therefore $F$ can also be analytically continued for $\tmmathbf{\xi} \in
\mathbb{C}^2$ that lives in a neighbourhood of $\tmmathbf{\mathcal{B}}$.

\begin{rema}
For Dirichlet BC or continuously varying properties, one can simply let $A =
S_0$, $B = \emptyset$ and $\chi (\tmmathbf{r}) \equiv 1$. The operator
$\mathcal{L}_{\tmmathbf{\xi}}$ then simplifies to $- \left( \left(
\frac{\partial}{\partial x_1} - i \beta_1 \right)^2 + \left(
\frac{\partial}{\partial x_2} - i \beta_2 \right)^2 \right)$, where
$\tmmathbf{\beta}= (\beta_1, \beta_2)$ is a linear function of
$\tmmathbf{\xi}$ defined by $\tmmathbf{\beta} \equiv (\Lambda^{- 1})^{\transpose}
\tmmathbf{\xi}$. Introducing such operator is quite standard in media with
periodic material properties, see e.g.\ {\cite{Bensoussan1978}} (chapter 4, p. 615),
{\cite{OdehKeller64}}, where it is known as the operator for the 	{\tmem{shifted cell problem}}.
\end{rema}

\section{\RED{Proof of the bypass choice procedure} } \label{app:procedure-bypass-choice}

\RED{Consider an irreducible component $\sigma_j(k,\varkappa)$ of the singularity set of $F$ with defining analytic function $g_j(\bdxi;k,\varkappa)$.} For simplicity, let us introduce the function $G_j (\tmmathbf{\xi};
\tilde{k})$, such that $G_j (\tmmathbf{\xi}; k + i \varkappa) \equiv g_j
(\tmmathbf{\xi}; k, \varkappa)$. Let us  fix $k > 0$ and consider a point
$\tmmathbf{\xi}^{\star} \in \sigma_j' (k)$, that is we have $G_j
(\tmmathbf{\xi}^{\star} ; k) \equiv g_j (\tmmathbf{\xi}^{\star} ; k, 0) = 0$.
Now, provided that $G_j$ is analytic in a $\mathbb{C}^3$ neighbourhood of
$(\tmmathbf{\xi}^{\star} ; k)$ and that $\frac{\partial G_j}{\partial
\tilde{k}} (\tmmathbf{\xi}^{\star} ; k) \neq 0$, the implicit function theorem
implies that for $(\tmmathbf{\xi}; \tilde{k})$ in a $\mathbb{C}^3$
neighbourhood of $(\tmmathbf{\xi}^{\star} ; k)$, there exists a unique
analytic function $f$ such that $G_j (\tmmathbf{\xi}; \tilde{k}) = 0
\Leftrightarrow \tilde{k} = f (\tmmathbf{\xi})$.

Let us consider $\tilde{k} = k + i \varkappa$ for some small $\varkappa > 0$
and look for points $\tmmathbf{\xi} \in \sigma_j (k, \varkappa)$ such that
$\tmmathbf{\xi}=\tmmathbf{\xi}^{\star} {+ i\tmmathbf{\xi}^{\ddagger}} $ for
some ${\tmmathbf{\xi}^{\ddagger}}  \in \mathbb{R}^2$. Using the implicit
function theorem, valid since we can choose $\varkappa$ as small as we want
and ${\tmmathbf{\xi}^{\ddagger}}  \rightarrow \tmmathbf{0}$ as $\varkappa
\rightarrow 0$, we obtain $k + i \varkappa = f \left( \tmmathbf{\xi}^{\star}
{+ i\tmmathbf{\xi}^{\ddagger}}  \right)$, which can be Taylor expanded to
obtain
\begin{eqnarray}
\varkappa & \approx & {\tmmathbf{\xi}^{\ddagger}}  \cdot \nabla
f (\tmmathbf{\xi}^{\star}).  \label{eq:secondwaymagic}
\end{eqnarray}
By definition of $f$, we know that $\nabla f
(\tmmathbf{\xi}^{\star}) \perp \sigma_j' (k)$, and therefore, since $\varkappa
> 0$, (\ref{eq:secondwaymagic}) implies that ${\tmmathbf{\xi}^{\ddagger}} $
and $\nabla f (\tmmathbf{\xi}^{\star})$ point towards the same
side of $\sigma_j' (k)$. 
Now we want to choose a surface indentation about $\bdxi^\star$ such that it does not cross the singularity, that is that in the limiting process $\bdxi^\ddagger \to \boldsymbol{0}$, the surface is not crossed. In other words, we should have $\bdxi^\star+i \bdeta \neq  \bdxi^\star+i \bdxi^\ddagger$ as $\bdxi^\ddagger \to \boldsymbol{0}$. To ensure that this is the case, it is enough to insist that $\bdeta$ and $\bdxi^\ddagger$ point to a different side of $\sigma'(k)$.
%{\color{red} Since ${\tmmathbf{\xi}^{\ddagger}}$ corresponds to as singularity, and the surface of integration cannot cross the singularity as ${\tmmathbf{\xi}^{\ddagger}} \to 0$,  $\tmmathbf{\eta} (\tmmathbf{\xi}^{\star})$  points in the opposite direction.}
Hence, if we can
find which side of $\sigma_j' (k)$ the vector $\nabla f
(\tmmathbf{\xi}^{\star})$ is pointing to, the bypass is determined by choosing $\bdeta$ to point in the other direction.
But, by definition, $\nabla f (\tmmathbf{\xi}^{\star})$,
which is also known as the {\tmem{group velocity vector}}, points towards
increasing real $k$ (see e.g. \cite{Lighthill78,chapman-wavenumber-surface}).

\section{Resonance excitation of a discrete lattice} 
\label{App:E}
Assume that the unknown $\mathfrak{u} (\tmmathbf{m}, t)$ satisfies the discrete linear wave
equation with effective wave speed $c$ subject to a harmonic forcing at the
origin starting at $t = 0$:
\begin{align}
&\mathfrak{u} (\tmmathbf{m}+\tmmathbf{e}_1 ; t) +\mathfrak{u}
(\tmmathbf{m}-\tmmathbf{e}_1 ; t) +\mathfrak{u} (\tmmathbf{m}+\tmmathbf{e}_2
; t)+\mathfrak{u} (\tmmathbf{m}-\tmmathbf{e}_2 ; t)  \nonumber\\
& - 4\mathfrak{u}
(\tmmathbf{m}; t) - \frac{1}{c^2} \frac{\mathd^2
\mathfrak{u}}{\mathd t^2} (\tmmathbf{m}; t)  =  \delta_{m_1 0} \delta_{m_2
0} e^{\RED{-}i \omega^{\star} t} \mathcal{H} (t)  , 
\label{eq:discrete-governing-time}
\end{align}
where $\omega^{\star}$ is specifically chosen to be $\omega^{\star} =
k^{\star} c$, for some $k^{\star} > 0$ and $\mathcal{H}$ is the Heaviside
function. Let us apply the double Floquet-Bloch transform in $\tmmathbf{m}$
and the Fourier transform in $t$ of (\ref{eq:discrete-governing-time}). Upon
denoting
\begin{eqnarray*}
\mathfrak{F} (\tmmathbf{\xi}; \omega) & : = & \int_{- \infty}^{\infty} \quad
\sum_{\tmmathbf{m} \in \mathbb{Z}^2} \, \mathfrak{u} (\tmmathbf{m}, t)
e^{i\tmmathbf{m} \cdot \tmmathbf{\xi}}  e^{\RED{i} \omega t}
\mathd t,
\end{eqnarray*}
direct computations lead to
\begin{eqnarray*}
\mathfrak{F} (\tmmathbf{\xi}; \omega) & = & \frac{i}{\left( 2 \cos (\xi_1) +
2 \cos (\xi_2) + \left( \frac{\omega^2}{c^2} - 4 \right) \right)
(\RED{\omega} - \RED{\omega^{\star}})}
\end{eqnarray*}
with inverse
\begin{eqnarray*}
\mathfrak{u} (\tmmathbf{m}, t) & = & \frac{1}{8 \pi^3} \int_{- \infty
{\RED{+} i \varepsilon}}^{\infty {\RED{+} i \varepsilon}}
\iint_{\tmmathbf{\mathcal{B}}} \mathfrak{F} (\tmmathbf{\xi}; \omega) e^{-
i\tmmathbf{m} \cdot \tmmathbf{\xi}}{e^{\RED{-}i \omega t}} \mathd
\tmmathbf{\xi} \mathd \omega .
\end{eqnarray*}
Let us now  change the variable $\omega
\leftrightarrow kc$ in the $\omega$ integral to get
%\begin{widetext}
\begin{eqnarray*}
\mathfrak{u} (\tmmathbf{m}, t) & = & \frac{\RED{i}}{8 \pi^3}
\iint_{\tmmathbf{\mathcal{B}}} e^{- i\tmmathbf{m} \cdot \tmmathbf{\xi}}
\int_{- \infty \RED{+} i \varepsilon}^{\infty \RED{+} i \varepsilon} \frac{e^{
		\RED{-}ikct}}{(k - \Theta (\tmmathbf{\xi})) (k + \Theta (\tmmathbf{\xi})) (k -
	k^{\star})} \mathd k \mathd \tmmathbf{\xi},
\end{eqnarray*}
where $\Theta (\tmmathbf{\xi}) = \sqrt{4 - 2 \cos (\xi_1) - 2 \cos (\xi_2)}$.
Let us consider the inner integral for some fixed $\tmmathbf{\xi} \in
\tmmathbf{\mathcal{B}}$. The integrand has three poles in the $k$-plane: $k_1
= k^{\star}$, $k_2 = \Theta (\tmmathbf{\xi})$ and $k_3 = - \Theta
(\tmmathbf{\xi})$ that happen to be real. By closing the contour of the $k$
integral in the \RED{lower}-half $k$-plane, we get
\begin{eqnarray}
\mathfrak{u} (\tmmathbf{m}, t) & = & \frac{1}{4 \pi^2}
\iint_{\tmmathbf{\mathcal{B}}} e^{- i\tmmathbf{m} \cdot \tmmathbf{\xi}} G
(\tmmathbf{\xi}; k^{\star} ; t) \mathd \tmmathbf{\xi}, 
\label{eq:doubleGintegral}
\end{eqnarray}
where
\begin{eqnarray}
G (\tmmathbf{\xi}; k^{\star} ; t) & = & \frac{e^{\RED{-} ik^{\star}
		ct}}{(k^{\star} - \Theta (\tmmathbf{\xi})) (k^{\star} + \Theta
	(\tmmathbf{\xi}))} - \frac{e^{ \RED{-}i \Theta (\tmmathbf{\xi}) ct}}{2 \Theta
	(\tmmathbf{\xi}) (k^{\star} - \Theta (\tmmathbf{\xi}))} + \frac{e^{\RED{ i}
		\Theta (\tmmathbf{\xi}) ct}}{2 \Theta (\tmmathbf{\xi}) (k^{\star} + \Theta
	(\tmmathbf{\xi}))} \cdot  \label{eq:G-3-terms}
\end{eqnarray}
%\end{widetext}
Recall that we did all this because we were interested in the singularity
defined by $k^{\star} = \Theta (\tmmathbf{\xi})$. It is remarkable that
through this change of integration order, the resulting integrand $G
(\tmmathbf{\xi}; k^{\star};t)$ is actually regular on the set $k^{\star} =
\Theta (\tmmathbf{\xi})$. Indeed, the third term is straightforwardly regular
there, while the singularities of each of the first two terms cancel each
other. Therefore, the integral (\ref{eq:doubleGintegral}) is actually
well-defined, unlike the integral (\ref{eq:inverse-BF-transform}).

In the simple case when $k^{\star} \neq 2$, we can consider that the first
term of (\ref{eq:G-3-terms}) yields the time harmonic solution, while the
other two terms are transient. Indeed, just by considering the first term we
obtain the exact same integral as (\ref{eq:integral-discrete}) multiplied by
the time harmonic factor $e^{\RED{-} i \omega^{\star} t}$. {The other
two terms can be shown to decay as $t \rightarrow \infty$.}
%\tmcolor{blue}{{\tmstrong{AS:}} I've not yet convinced myself on how to show
%	the last sentence. By taking an absorptive media, did you mean to consider $c$
%	and $k$ complex but $\omega$ real?} So everything is consistent and shows how
%the limiting absorption principle can be thought of as a radiation condition.

Let us now assume that $k^{\star} = 2$ and study the specific crossing at
$\tmmathbf{\xi}^{\star} = \RED{(0, \pi)}$. Consider the linear change of variables
$\tmmathbf{\xi} \leftrightarrow \tmmathbf{\zeta}$ defined by
\begin{eqnarray*}
\zeta_1 = \frac{\RED{\xi_1-(\xi_2-\pi)}}{2} & \infixand & \zeta_2 = \frac{\RED{\xi_1+\xi_2-\pi}}{2},
\end{eqnarray*}
leading to $\tmmathbf{\xi}-\tmmathbf{\xi}^{\star} = \Psi^{\star}
\tmmathbf{\zeta}$, where $\Psi^{\star} = \left( \begin{array}{cc}
\RED{1} & \RED{1}\\
\RED{-}1 & \RED{1}
\end{array} \right)$ and $\det (\Psi^{\star}) = 2$. Note that, as
$\tmmathbf{\xi} \rightarrow \tmmathbf{\xi}^{\star}$, $k^{\star} - \Theta
(\tmmathbf{\xi}) \approx \RED{-}\zeta_1 \zeta_2$, and that this coincides with the
ansatz (\ref{eq:Omega-hyp}). The vector $\tmmathbf{\alpha}$ defined by
$\tmmathbf{{\alpha}}=(\Psi^\star)^{\transpose} \tmmathbf{m}$ is given by
$\alpha_1 = m_{\RED{1}} - m_{\RED{2}}$ and $\alpha_2 =  m_1 \RED{+} m_2$.  Using this change of
variables, the local contribution due to $\tmmathbf{\xi}^{\star}$ arising from
(\ref{eq:doubleGintegral}) is given by
\begin{align*}
\mathfrak{u}_{\tmop{loc}} (\tmmathbf{m}, t) &= \RED{-}(-1)^{m_{\RED{2}}}\frac{e^{\RED{-}ik^{\star} ct}}{8
\pi^2} \times \mathcal{I}_{\tmop{time}}, \quad \tmop{where} \\
\mathcal{I}_{\tmop{time}} &= \iint_{\tmmathbf{\Pi}^{\infty}} \frac{(1 - e^{-i
	\zeta_1 \zeta_2 ct})}{\zeta_1 \zeta_2} e^{- i\tmmathbf{\alpha} \cdot
\tmmathbf{\zeta}} \mathd \tmmathbf{\zeta}.
\end{align*}
%{\color{red}\textbf{AIK}:The $(-1)^{m_1}$ comes from a term $e^{-i\bdm\cdot\bdxi^\star}$, which I think we forgot in previous versions. Do you agree?}{\color{green}Agree}
When $\alpha_1 > 0$ and
$\alpha_2 > 0$, the integral $\mathcal{I}_{\tmop{time}}$ can be shown to be 
\begin{align}
\mathcal{I}_{\tmop{time}} &= 2 i \pi e^{ i \frac{\alpha_1 \alpha_2}{ct}}
\Phi \left( \frac{\alpha_1 \alpha_2}{ct} \right), \text{ where } \Phi (z)
\assign \int_0^{\infty} \frac{e^{ iz \tau}}{1 + \tau} \mathd \tau . 
\label{eq:Itimeexpr}
\end{align}
Noting that $\Phi (z) \sim \log (z)$ as $z \rightarrow 0$, one can see that
$\mathcal{I}_{\tmop{time}}$ blows up like $\log (1 / t)$ as $t \rightarrow
\infty$, which is indeed representative of a resonance phenomena, and
illustrates why it is not possible to obtain a stationary solution. 
\section{Derivation of  (\ref{eq:simple-approx-F-double-eig})}\label{app:A}

%\subsection{Approximating $\psi_{\tmmathbf{\xi}} /
%	\psi_{\tmmathbf{\xi}^{\star}}$ to simplify (\ref{eq:double-eigenproj-F})}
\RED{Remember that we consider a pair $(\tmmathbf{\xi}^{\star}, k^{\star}) \in
\tmmathbf{\mathcal{B}} \times \mathbb{R}^+$ and an index $\ell \in
\mathbb{N}$, such that $\lambda_{\ell} (\tmmathbf{\xi}^{\star}) \geqslant 0$
is a double eigenvalue of $(\mathcal{L}, \tmop{dom}_{\tmmathbf{\xi}^{\star}})$
and $(k^{\star})^2 = \lambda_{\ell} (\tmmathbf{\xi}^{\star})$.

The two associated orthogonal eigenfunctions are $V_{\ell}^{(1)} (\tmmathbf{r})$ and $V_{\ell}^{(2)}
(\tmmathbf{r})$, which are normalised such
that $(V_{\ell}^{(\mathscr{m})}, V_{\ell}^{(\mathscr{k})}) =
\delta_{\mathscr{m}, \mathscr{k}}$ for $\mathscr{m}, \mathscr{k} \in \{ 1, 2
\}$.}
For simplicity, we assume that all the other eigenvalues are simple and
remain so in a neighbourhood of $\tmmathbf{\xi}^{\star}$. We also assume that
the other orthogonal eigenfunctions $V_j (\tmmathbf{r};
\tmmathbf{\xi}^{\star})$ are also normalised so that $(V_j (\cdot ;
\tmmathbf{\xi}^{\star}), V_j (\cdot ; \tmmathbf{\xi}^{\star})) = 1$. This set
of eigenfunctions spans $\tmop{dom}_{\tmmathbf{\xi}^{\star}}$, and using the
ideas of \RED{Appendix} \ref{sec:intro-psi-xi}, we know that these eigenfunctions
divided by $\psi_{\tmmathbf{\xi}^{\star}}$ \RED{(defined in \ref{eq:def-Psi-Xi})} span $\tmop{dom}_{\tmmathbf{0}}$.

So let us consider $(\tmmathbf{\xi}, k)$ in the
neighbourhood of $(\tmmathbf{\xi}^{\star}, k^{\star})$. Since the function $F
(\tmmathbf{r}; \tmmathbf{r}_s ; \tmmathbf{\xi}; k)$ is in
$\tmop{dom}_{\tmmathbf{\xi}}$, we know that $F (\tmmathbf{r}; \tmmathbf{r}_s ;
\tmmathbf{\xi}; k) / \psi_{\tmmathbf{\xi}}$ is in $\tmop{dom}_{\tmmathbf{0}}$.
This allows us to write $F$ in terms of the eigenfunctions associated to
$\tmmathbf{\xi}^{\star}$:
%\begin{widetext}
\begin{eqnarray}
F (\tmmathbf{r}; \tmmathbf{r}_s ; \tmmathbf{\xi}; k) & = & \left( a_1
V_{\ell}^{(1)} (\tmmathbf{r}) + a_2 V_{\ell}^{(2)} (\tmmathbf{r}) + \sum_{j
	\neq \ell} b_j V_j (\tmmathbf{r}; \tmmathbf{\xi}^{\star}) \right)
\frac{\psi_{\tmmathbf{\xi}} (\tmmathbf{r})}{\psi_{\tmmathbf{\xi}^{\star}}
	(\tmmathbf{r})}  \label{eq:double-eigenproj-F}
\end{eqnarray}
%\end{widetext}
for some unknown constants $a_1, a_2$ and $b_j$ to be determined.

Let us apply perturbation analysis to (\ref{eq:double-eigenproj-F}). Let $\delta_\xi = |\tmmathbf{\xi}-\tmmathbf{\xi}^\star|$ be a small parameter. Then,  
\begin{equation}
\frac{\psi_{\tmmathbf{\xi}}}{\psi_{\tmmathbf{\xi}^{\star}}}\approx 1 + i(\tmmathbf{\xi}-\tmmathbf{\xi}^\star)\cdot\Lambda^{-1}\tmmathbf{r}\chi(\tmmathbf{r}).
\end{equation}
Suppose that 
\begin{equation}
\label{eq:pert_dep}
a_{1,2} = \mathcal{O}(1/\delta_\xi), \quad b_{j} = \mathcal{O}(1).
\end{equation}
Then, using the latter, approximate  $F$ as follows:
\begin{align}
F (\tmmathbf{r}; \tmmathbf{r}_s ; \tmmathbf{\xi}; k) &\approx  a_1 V_{\ell}^{(1)} (\tmmathbf{r}) + a_2 V_{\ell}^{(2)}
(\tmmathbf{r}) \nonumber \\ &+ \sum_{j \neq \ell} b_j V_j (\tmmathbf{r};
\tmmathbf{\xi}^{\star}) + \Upsilon (\tmmathbf{r}), 
\label{eq:rewrite-Vl-part3}
\end{align}
where we defined
\begin{align}
\Upsilon (\tmmathbf{r}) & \equiv  i(\tmmathbf{\xi}-\tmmathbf{\xi}^\star)\cdot\Lambda^{-1}\tmmathbf{r}\chi(\tmmathbf{r})
(a_1 V_{\ell}^{(1)} (\tmmathbf{r}) + a_2 V_{\ell}^{(2)}
(\tmmathbf{r})).
\label{eq:gamma_def}	
\end{align}
Note that we have kept only the terms of order $1/\delta_\xi$ and $1$. Let us use the orthogonality of the eigenfunctions to determine the coefficients $a_{1,2}$ and $b_{1,2}$. Apply the Helmholtz operator $\mathcal{L}- k^2$ to (\ref{eq:rewrite-Vl-part3}) to get rid of $F$:
\begin{align}
- \delta (\tmmathbf{r}_s) & =  \delta_{\lambda}^{(\ell)} (a_1
V_{\ell}^{(1)} (\tmmathbf{r}) + a_2 V_{\ell}^{(2)} (\tmmathbf{r}))  \nonumber \\&+ \sum_{j
\neq \ell} b_j \delta_{\lambda}^{(j)} V_j (\tmmathbf{r};
\tmmathbf{\xi}^{\star}) +\mathcal{L} [\Upsilon] - k^2 \Upsilon
(\tmmathbf{r})  \label{eq:after-op-L-applied}
\end{align}
where we defined
\begin{align}
\delta_{\lambda}^{(\ell)} \assign \lambda_{\ell} (\tmmathbf{\xi}^{\star}) -
k^2 &= (k^{\star})^2 - k^2 \text{ and } \delta_{\lambda}^{(j)} \assign \lambda_j
(\tmmathbf{\xi}^{\star}) - k^2,  
\end{align}
for $j \neq \ell$. We assume that $\delta_{\lambda}^{(\ell)}$ is small in the sense that $\delta_{\lambda}^{(\ell)}=\mathcal{O}(\delta_\xi)$.
%\[
%\delta_{\lambda}^{(\ell)} \sim \delta\xi.
%\]
Then,  taking the product ($\left. \left( {\rm \ref{eq:after-op-L-applied}} \right. \right),
V_{\ell}^{(\mathscr{m})}$) where $\mathscr{m} \in \{ 1, 2 \}$, using orthogonality and keeping terms of order $1$ leads to the following system of equations for $a_{1,2}$:
\begin{align}
- (V_{\ell}^{(\mathscr{m})} (\tmmathbf{r}_s))^{\dagger} & \approx 
\delta_{\lambda}^{(\ell)} a_{\mathscr{m}} + (\mathcal{L} [\Upsilon],
V_{\ell}^{(\mathscr{m})}) - (\Upsilon, \mathcal{L}V_{\ell}^{(\mathscr{m})}).
\label{eq:res1-raph}
\end{align}

The idea is now to use Green's theorem in order to rewrite the last two terms
as an integral over $\partial S_0$ (see Figure~\ref{fig:def-A-B-cell}). Since $\mathcal{L}$ is not self-adjoint in $\tmop{dom}_{\tmmathbf{\xi}}$, a boundary term emerges in Green's formula.  Namely, we have
%\begin{widetext}
\begin{eqnarray}
(\mathcal{L} [\Upsilon], V_{\ell}^{(\mathscr{m})}) - (\Upsilon,
\mathcal{L}V_{\ell}^{(\mathscr{m})}) & = & \int_{\partial S_0} [\Upsilon (\tmmathbf{n} \cdot \nabla
V_{\ell}^{(\mathscr{m})})^{\conj} - (\tmmathbf{n} \cdot \nabla \Upsilon)
(V_{\ell}^{(\mathscr{m})})^{\conj}] \mathd s,  \label{eq:weird-int-term}
\end{eqnarray}
where we took into account that $\Upsilon$ is zero on $\ptl \Omega$.
Then, after some algebra and using the periodicity of $V_{\ell}^{(\mathscr{m})}$ and (\ref{eq:gamma_def}) we get: 
\begin{eqnarray*}
(\mathcal{L} [\Upsilon], V_{\ell}^{(\mathscr{m})}) - (\Upsilon,
\mathcal{L}V_{\ell}^{(\mathscr{m})}) & = & - \delta_{\xi_1} (a_1
Y_{\mathscr{m}, 1}^{(L)} + a_2 Y_{\mathscr{m}, 2}^{(L)}) - \delta_{\xi_2}
(a_1 Y_{\mathscr{m}, 1}^{(B)} + a_2 Y_{\mathscr{m}, 2}^{(B)})
\end{eqnarray*}
where we defined $\delta_{\xi_1} = (\xi_1-\xi^\star_1)$, $\delta_{\xi_2} = (\xi_2-\xi^\star_2)$ and
\begin{eqnarray*}
Y_{\mathscr{m}, \mathscr{k}}^{(L)} & = & i \int_{\partial S_0^L}
V_{\ell}^{(\mathscr{k})} (\tmmathbf{r}) (\tmmathbf{n} \cdot \nabla
V_{\ell}^{(\mathscr{m})})^{\conj} - (\tmmathbf{n} \cdot \nabla
V_{\ell}^{(\mathscr{k})} (\tmmathbf{r})) (V_{\ell}^{(\mathscr{m})})^{\conj}
\mathd s, \\
Y_{\mathscr{m}, \mathscr{k}}^{(B)} & = & i \int_{\partial S_0^B}
V_{\ell}^{(\mathscr{k})} (\tmmathbf{r}) (\tmmathbf{n} \cdot \nabla
V_{\ell}^{(\mathscr{m})})^{\conj} - (\tmmathbf{n} \cdot \nabla
V_{\ell}^{(\mathscr{k})} (\tmmathbf{r})) (V_{\ell}^{(\mathscr{m})})^{\conj}
\mathd s.
\end{eqnarray*}

Finally, we obtain the following expression for the coefficients $\bda=(a_1,a_2)^\transpose$:
\begin{eqnarray*}
\tmmathbf{a}& = & (\mathcal{Y}- \delta_{\lambda}^{(\ell)}{\rm I}_2)^{-1} \left(
\begin{array}{c}
	(V_{\ell}^{(1)} (\tmmathbf{r}_s))^{\conj}\\
	(V_{\ell}^{(2)} (\tmmathbf{r}_s))^{\conj}
\end{array} \right),
\end{eqnarray*}
where ${\rm I}_2$ is the $2\times2$ identity matrix and
\begin{equation}
\mathcal{Y}= \delta_{\xi_1} \mathcal{Y}^{(L)} + \delta_{\xi_2}
\mathcal{Y}^{(B)}, \quad \mathcal{Y}^{(L)} = \left( \begin{array}{cc}
	Y_{1, 1}^{(L)} & Y_{1, 2}^{(L)}\\
	Y_{2, 1}^{(L)} & Y_{2, 2}^{(L)}
\end{array} \right), \quad  \mathcal{Y}^{(B)} = \left( \begin{array}{cc}
	Y_{1, 1}^{(B)} & Y_{1, 2}^{(B)}\\
	Y_{2, 1}^{(B)} & Y_{2, 2}^{(B)}
\end{array} \right), 
\label{eq:matrixYformula}
\end{equation}
%\end{widetext}
therefore recovering the formula (\ref{eq:asandbs}).

Using a similar approach, one can obtain the expressions for the coefficients $b_j$:
\begin{eqnarray*}
b_j & \approx & - \frac{V_j (\tmmathbf{r}_s ; \tmmathbf{\xi}^{\star})} 
{\delta_{\lambda}^{(j)}} \cdot
\end{eqnarray*}

Note that for small $\delta^{(\ell)}_\lambda$ and $\delta_\xi$ the coefficients $a_\mathscr{m}$ dominate over the coefficients $b_j$, which is consistent with the assumption (\ref{eq:pert_dep}).

%
%\subsection{Exploiting orthogonality to approximate the dispersion diagram}
%
%Everything done so far in this appendix holds if we start from
%(\ref{eq:double-eigenproj-Vl}) instead of starting from
%(\ref{eq:double-eigenproj-F}). As a result, we get the following equation for
%the coefficient $\tilde{\tmmathbf{a}} = (\tilde{a}_1, \widetilde{a_2})$:
%\begin{eqnarray*}
%	\mathcal{Y} \tilde{\tmmathbf{a}} - \delta_{\lambda}^{(\ell)}
%	\tilde{\tmmathbf{a}} & = & 0
%\end{eqnarray*}
%This equation only admits non-zero solutions if $\det (\mathcal{Y}-
%\delta_{\lambda}^{\ell} I_2) = 0$. Solving such equation for $(\tmmathbf{\xi},
%k)$ therefore provides a l to the dispersion diagram.

\section{FEM model of a phononic crystal}
\label{FEM_construction}
\subsection{Structure of the stiffness and mass matrices for quasi-periodic boundary conditions}

%Let us consider discrete analogue of the problem formulated in (\ref{sec:continuous-example}). First, consider one unit cell $S_0/\Omega$. As it is shown in Figure \ref{fig:FEM_unit_cell}, $S_0/\Omega$ could be easily covered by a mesh of equilateral triangles (here we are using linear triangular finite elements). 

Following a classical FEM procedure, one can obtain the stiffness matrix $\rm{K}_0$ and the mass matrix $\rm{M}_0$ for a cell, such that the infinite matrices $\rm{K}$ and $\rm{M}$ are made up of infinitely many such $N_{\text{cell}}\times N_{\text{cell}}$ block matrices $\rm{K}_0$ and $\rm{M}_0$ respectively. Note that by definition of $\rm{K}_0$ and $\rm{M}_0$, we can define
\begin{equation*}
E_{\text{kin}}^{\bdm} = \frac{1}{2} \dot{\rm U}_{\bdm}^{\dagger} \textrm{K}_0 \dot{\rm U}_{\bdm}, \text{ and }
E_{\text{pot}}^{\bdm} = \frac{1}{2} {\rm U}_{\bdm}^{\dagger} \textrm{M}_0 {\rm U}_{\bdm},
\end{equation*}
corresponding to the kinetic and potential energy of the $\bdm^\text{th}$ meshed cell respectively, where ${\rm U_{\bdm}}$ is a vector of field values for the $\bdm^\text{th}$ cell, and $\dot{\{\}}$ stands for time derivative, and where, whenever used on a vector, ${(\cdot)}^{\dagger}$ stands for the conjugate transpose. 
%All elements of matrix $\textrm{K}^0$ corresponding to non-neighbouring nodes should be equal to zero. 
The $N_{\text{cell}}\times1$ vector ${\rm F(\bdxi)}$ defined by Floquet-Bloch transformation in Example \ref{ex:FEM-stuff} can be written in the block form
%Let ${\rm U}$ be represented in the block form
\begin{equation}
{\rm F} = (
{\rm F}_i^\transpose, 
{\rm F}_l^\transpose,
{\rm F}_b^\transpose,
{\rm F}_r^\transpose,
{\rm F}_u^\transpose,
{\rm F}_{c_1}^\transpose,
{\rm F}_{c_2}^\transpose,
{\rm F}_{c_3}^\transpose,
{\rm F}_{c_4}^\transpose)^\transpose,
\end{equation}
%\begin{equation}
%    {\rm U} = \begin{pmatrix}
%    {\rm U}_i \\ 
%    {\rm U}_l \\
%    {\rm U}_b \\
%    {\rm U}_r \\
%    {\rm U}_u \\
%    {\rm U}_{c_1} \\
%    {\rm U}_{c_2} \\
%    {\rm U}_{c_3} \\
%    {\rm U}_{c_4} \\
%\end{pmatrix},
%\end{equation}
\noindent where $\{i, l, b, r,u\}$ correspond to the \{inner, left, bottom, right, upper\} nodes, and $\{c_1,c_2,c_3,c_4\}$ indicate the \{bottom left, upper left, upper right, bottom right\} corner nodes (see Figure \ref{fig:FEM_unit_cell}). Due to the $\bdxi$ quasi-periodicity of ${\rm F}(\bdxi)$, the following conditions should be satisfied:
\begin{alignat}{2}
\label{eq:app_quasi}
{\rm F}_r &= e^{-i \xi_1} {\rm F}_l,
\quad
{\rm F}_u &&= e^{-i \xi_2} {\rm F}_b,
\\
{\rm F}_{c_2} &= e^{-i \xi_2} {\rm F}_{c_1},
\quad
{\rm F}_{c_3} &&= e^{-i (\xi_1 + \xi_2)} {\rm F}_{c_1},
\quad
{\rm F}_{c_4} = e^{-i \xi_1} {\rm F}_{c_1}. \nonumber
\end{alignat}
%where $\xi_1, \xi_2 \in [-\pi, \pi]$. Define vector 
%\begin{equation}
%    \boldsymbol{\xi} = \begin{pmatrix}
%    \xi_1 \\ 
%    \xi_2 \\
%    \end{pmatrix}.
%\end{equation}
Similarly, write the $N_{\text{cell}}\times N_{\text{cell}}$ matrix $\rm K_0$ in the block form:
\begin{equation}
\textrm{K}_0 = \begin{pmatrix}\label{appendix_block_matrix}
	\textrm{K}_{ii} & \textrm{K}_{il} & \textrm{K}_{ib} & \textrm{K}_{ir} & \textrm{K}_{iu} & \textrm{K}_{ic_1} & \dots & \textrm{K}_{ic_4}\\
	\textrm{K}_{li} & \textrm{K}_{ll} & \textrm{K}_{lb} & \textrm{K}_{lr} & \textrm{K}_{lu} & \textrm{K}_{lc_1} & \dots & \textrm{K}_{lc_4}\\
	\textrm{K}_{bi} & \textrm{K}_{bl} & \textrm{K}_{bb} & \textrm{K}_{br} & \textrm{K}_{bu} & \textrm{K}_{bc_1} & \dots & \textrm{K}_{bc_4}\\
	\textrm{K}_{ri} & \textrm{K}_{rl} & \textrm{K}_{rb} & \textrm{K}_{rr} & \textrm{K}_{ru} & \textrm{K}_{rc_1} & \dots & \textrm{K}_{rc_4}\\
	\textrm{K}_{ui} & \textrm{K}_{ul} & \textrm{K}_{ub} & \textrm{K}_{ur} & \textrm{K}_{uu} & \textrm{K}_{uc_1} & \dots & \textrm{K}_{uc_4}\\
	\textrm{K}_{c_1i} & \textrm{K}_{c_1l} & \textrm{K}_{c_1b} & \textrm{K}_{c_1r} & \textrm{K}_{c_1u} & \textrm{K}_{c_1c_1} & \dots & \textrm{K}_{c_1c_4}\\
	\vdots & \vdots & \vdots & \vdots & \vdots & \vdots & \ddots & \vdots \\
	\textrm{K}_{c_4i} & \textrm{K}_{c_4l} & \textrm{K}_{c_4b} & \textrm{K}_{c_4r} & \textrm{K}_{c_4u} & \textrm{K}_{c_4c_1} & \dots & \textrm{K}_{c_4c_4}\\
	
\end{pmatrix}.
\end{equation}
Following the discussion in Example \ref{ex:FEM-stuff}, we can introduce the reduced $\hat{N}\times 1$ vector $\hat{\rm F}$ that contains all the necessary information and is defined by
\begin{equation}
\label{eq:hat_F}
{\rm \hat{F}} = (
{\rm F}_i^\transpose,
{\rm F}_l^\transpose,
{\rm F}_b^\transpose,
{\rm F}_{c_1}^\transpose)^\transpose,
\end{equation}
so that $\hat{N}=N_i+N_l+N_b+1=N_{\text{cell}}-N_r-N_u-3$, where $\{N_i,N_l,N_b,N_r,N_u\}$ correspond to the number of $\{\text{inner},\text{left},\text{bottom},\text{right},\text{upper}\}$ nodes in a cell.
Then, it follows from (\ref{eq:app_quasi})--(\ref{appendix_block_matrix}) that the associated reduced matrix $\hat{\textrm{K}}(\boldsymbol{\xi})$ introduced in Example \ref{ex:FEM-stuff} has the following structure: 
%\begin{widetext}
\begin{equation}
	\hat{\textrm{K}}(\boldsymbol{\xi}) \equiv \begin{pmatrix}\label{appendix_block_matrix_r}
		\textrm{K}_{ii} & \textrm{K}_{il} + \textrm{K}_{ir}e^{-i\xi_1} & \textrm{K}_{ib} + \textrm{K}_{iu} e^{-i \xi_2}  & \hat{\textrm{K}}_{ic_1}(\boldsymbol{\xi})\\
		\textrm{K}_{li} + \textrm{K}_{ri} e^{i\xi_1} & \textrm{K}_{ll} + \textrm{K}_{rr} & \hat{\textrm{K}}_{lb}(\boldsymbol{\xi}) & \hat{\textrm{K}}_{lc_1}(\boldsymbol{\xi})\\
		\textrm{K}_{bi} + \textrm{K}_{ri} e^{i \xi_2} & \hat{\textrm{K}}_{bl}(\boldsymbol{\xi}) & \textrm{K}_{bb} + \textrm{K}_{uu} & \hat{\textrm{K}}_{bc_1}(\boldsymbol{\xi})\\
		\hat{\textrm{K}}_{c_1i}(\boldsymbol{\xi}) & \hat{\textrm{K}}_{c_1l}(\boldsymbol{\xi}) &  \hat{\textrm{K}}_{c_1b}(\boldsymbol{\xi}) & \hat{\textrm{K}}_{c_1c_1}
	\end{pmatrix},
\end{equation}
%\end{widetext}
where 
\begin{align*}
\hat{\textrm{K}}_{lb}(\boldsymbol{\xi}) &= \textrm{K}_{lb} + \textrm{K}_{rb} e^{i\xi_1} + \textrm{K}_{ru} e^{i(\xi_1 - \xi_2)} + \textrm{K}_{lu} e^{-i \xi_2},\\
\hat{\textrm{K}}_{bl}(\boldsymbol{\xi}) &=  \textrm{K}_{bl} + \textrm{K}_{ul} e^{i\xi_2} + \textrm{K}_{br} e^{-i\xi_1} + \textrm{K}_{ur} e^{i(\xi_2 - \xi_1)},\\
\hat{\textrm{K}}_{c_1c_1} &=  \textrm{K}_{c_1c_1} + \textrm{K}_{c_2c_2} + \textrm{K}_{c_3c_3} + \textrm{K}_{c_4c_4},\\
\hat{\textrm{K}}_{ic_1}(\boldsymbol{\xi}) &=  \textrm{K}_{ic_1} + \textrm{K}_{ic_2} e^{-i \xi_2} + \textrm{K}_{ic_3} e^{-i (\xi_1 + \xi_2)} + \textrm{K}_{ic_4} e^{-i \xi_1},\\
\hat{\textrm{K}}_{c_1i}(\boldsymbol{\xi}) &=  \textrm{K}_{c_1i} + \textrm{K}_{c_2i} e^{i \xi_2} + \textrm{K}_{c_3i} e^{i (\xi_1 + \xi_2)} + \textrm{K}_{c_4i} e^{i \xi_1},\\
\hat{\textrm{K}}_{lc_1}(\boldsymbol{\xi}) &=  \textrm{K}_{lc_1} + \left(\textrm{K}_{lc_2} + \textrm{K}_{lc_3} \right)  e^{-i \xi_2} + \textrm{K}_{lc_4},\\
\hat{\textrm{K}}_{c_1l}(\boldsymbol{\xi}) &=  \textrm{K}_{c_1l} + \left(\textrm{K}_{c_2l} + \textrm{K}_{c_3l} \right)  e^{i \xi_2} + \textrm{K}_{c_4l},\\
\hat{\textrm{K}}_{bc_1}(\boldsymbol{\xi}) &=  \textrm{K}_{bc_1} + \textrm{K}_{bc_2} + \left(\textrm{K}_{bc_3} + \textrm{K}_{bc_4} \right)  e^{i \xi_1}, \\
\hat{\textrm{K}}_{c_1b}(\boldsymbol{\xi}) &=  \textrm{K}_{c_1b} + \textrm{K}_{c_2b} +\left(\textrm{K}_{c_3b} + \textrm{K}_{c_4b} \right)  e^{-i \xi_1}.
\end{align*}
The matrix $\rm \hat{M}(\bdxi)$ can be built in a similar way and has the same structure as $\rm \hat{K}(\bdxi)$.

\subsection{Perturbation  analysis near a double eigenvalue}
Let us for simplicity consider the case of a lumped mass matrix, i.e. ${\hat{\rm M}}$ being diagonal, and thus not depending on $\tmmathbf{\xi}$. Let ${\hat{\rm V}}_{1, 2}(\boldsymbol{\xi}^\star)$ be the orthonormal eigenvectors of the problem (\ref{eq:eigenequation}) corresponding to the double eigenvalue \RED{$\hat{\lambda}_1(\bdxi^\star)=\hat{\lambda}_2(\bdxi^\star)=(k^\star)^2$}, so that
\begin{equation}
	\label{eq:V_norm}
	\hat{\rm V}_{i}^{\dagger}(\boldsymbol{\xi}^\star) \hat{\rm M} \hat{\rm V}_{j}(\boldsymbol{\xi}^\star) = \delta_{i j},
\end{equation}
where $\delta_{ij}$ is the Kronecker delta. \RED{In the same spirit as in the continuous case, we assume that when $\bdxi$ is close to $\bdxi^\star$, we can write}
\begin{align}
	\RED{F(\bdxi)\approx \sum_{j}^{\hat{N}} a_j(\bdxi) \hat{V}_j(\bdxi^\star)}.
	\label{eq:new-numerics-F-xi-star}
\end{align}		
Let us find the leading terms in the approximation (\ref{eq:new-numerics-F-xi-star}). %at a double eigenvalue point $\bdxi^\star$.
To do this, expand the matrix $\hat{\textrm{K}}(\tmmathbf{\xi})$ in a Taylor series near $\tmmathbf{\xi}^\star$ to get
\begin{equation}\label{eq:diff K}
	\hat{\textrm{K}}(\boldsymbol{\xi}) \approx \hat{\textrm{K}} (\boldsymbol{\xi}^\star) + \delta_{\xi_1} \textrm{D}_1 + \delta_{\xi_2} \textrm{D}_2, %\text{where } \textrm{D}_i = \frac{\partial \hat{\textrm{K}}}{\partial \xi_i}(\tmmathbf{\xi}^\star)
\end{equation}
where 
\begin{equation}
	\textrm{D}_i = \frac{\partial \hat{\textrm{K}}}{\partial \xi_i}(\tmmathbf{\xi}^\star),
\end{equation} 
for $i = {1, 2}$, and $\delta_{\xi_i}$ are defined in the same way as in Appendix \ref{app:A}.

Considering (\ref{eq:V_norm})--(\ref{eq:diff K}), the equation (\ref{eq:discreteFloqueteq}) becomes: %\RED{In the eq below think $\hat{\lambda}_j^2$ should be $\hat{\lambda}_j$} 
%\begin{widetext}
	\begin{align}\label{eq:eq of motion near Dirac point}
		[\delta_{\xi_1} \textrm{D}_1 +  \delta_{\xi_2} \textrm{D}_2  +  \delta_\lambda\hat{\textrm{M}}] (a_1 \hat{\rm V}_1(\boldsymbol{\xi}^\star) + a_2\hat{\rm V}_2(\boldsymbol{\xi}^\star)) + \sum_{j=3}^{\hat{N}} a_j\left[(\RED{\hat{\lambda}_j(\bdxi^\star)} - (k^{\star})^2)\hat{\textrm{M}}  \hat{\rm V}_j(\boldsymbol{\xi}^\star) + [\delta_{\xi_1} \textrm{D}_1 +  \delta_{\xi_2} \textrm{D}_2 + \delta_\lambda\hat{\textrm{M}}] \hat{\rm V}_j(\boldsymbol{\xi}^\star) \right] &= {\rm \hat E},
	\end{align} 
%\end{widetext}
where $\delta_\lambda = (k^\star)^2 - k^2$. Upon defining $\delta_\xi$ as in Appendix \ref{app:A}, let us further assume that
$$
a_1, a_2 = \mathcal{O}(1/\delta_\xi), \quad
a_\RED{j} =\mathcal{O}(1), \quad \RED{j} = 3, 4, ..., \hat{N}.
$$
After multiplying (\ref{eq:eq of motion near Dirac point}) successively by $\hat{\rm V}_1^\dagger(\boldsymbol{\xi})$, and $\hat{\rm V}_2^\dagger(\boldsymbol{\xi})$ and collecting the $\mathcal{O}(1)$ terms, we obtain the system of equations:
\begin{equation}\label{eq:eq of motions_eq1}
	({\mathcal{Y}} - \delta_\lambda {\rm I}_2)
	\begin{pmatrix}
		a_1\\
		a_2
	\end{pmatrix} = 
	%\begin{pmatrix}
	%\hat{\rm V}_1^\dagger(\bdxi^\star){\rm \hat E}\\
	%\hat{\rm V}_2^\dagger(\bdxi^\star){\rm \hat E}
	%\end{pmatrix}=
	\begin{pmatrix}
		\tilde{V}_1\\
		\tilde{V}_2
	\end{pmatrix}, 
\end{equation}
where $\mathcal{Y} = \delta_{\xi_1}{\rm Y}^{(1)} + \delta_{\xi_2}{\rm Y}^{(2)}$, and for $i\in\{1,2\}$ we define
\begin{align*}
	{\rm Y}^{(i)} &= - \begin{pmatrix}
		\hat{\rm V}_1^\dagger(\boldsymbol{\xi}^\star)  \textrm{D}_i \hat{\rm V}_1(\boldsymbol{\xi}^\star) & \hat{\rm V}_1^\dagger(\boldsymbol{\xi}^\star)  \textrm{D}_i \hat{\rm V}_2(\boldsymbol{\xi}^\star) \\
		\hat{\rm V}_2^\dagger(\boldsymbol{\xi}^\star)  \textrm{D}_i \hat{\rm V}_1(\boldsymbol{\xi}^\star) & \hat{\rm V}_2^\dagger(\boldsymbol{\xi}^\star)  \textrm{D}_i \hat{\rm V}_2(\boldsymbol{\xi}^\star) \\
	\end{pmatrix}, \text{ and } \\ \tilde {\rm V}_i &= -{\rm \hat V}_i^\dagger(\bdxi^\star){\rm \hat E}.
\end{align*}
Solving (\ref{eq:eq of motions_eq1}) by direct inversion, we obtain an analogue of (\ref{eq:a_coef_cont}): 
\begin{equation}
	\begin{pmatrix}
		a_1\\
		a_2
	\end{pmatrix} = 
	\frac{1}{D(\tmmathbf{\xi})}
	\begin{pmatrix}
		\mathcal{Y}_{22} - \delta_\lambda & - \mathcal{Y}_{12}\\
		- \mathcal{Y}_{21} & \mathcal{Y}_{11} - \delta_\lambda
	\end{pmatrix}
	\begin{pmatrix}
		\tilde{\rm V}_1\\
		\tilde{\rm V}_2
	\end{pmatrix}
\end{equation}
where
\begin{equation}\label{D}
	{D}(\boldsymbol{\xi}) = \det({\mathcal{Y}} - \delta_\lambda {\rm I}_2)= (\delta_\lambda)^2-\delta_\lambda \trace(\mathcal{Y})+\det(\mathcal{Y})
\end{equation}
%\begin{equation}\label{D}
%   {D}(\boldsymbol{\xi}) = \det({\mathcal{Y}} {\color{red}-} \delta_\lambda {\rm I}_2)= -\left(\delta_\lambda \right)^2 +  C_1\delta_\lambda\delta_{\xi_1} + C_2\delta_\lambda\delta_{\xi_2}+C_3 \left(\delta_{\xi_1} \right)^2 + C_4  \left(\delta_{\xi_2} \right)^2 + C_5 \delta_{\xi_1} \delta_{\xi_2} 
%\end{equation}
is the local dispersion function. 
%The coefficients $C_1, ..., C_5$ are: 
%\[
%   C_1 = {Y}^{(1)}_{11} + {Y}^{(1)}_{22},\quad
%   C_2 = {Y}^{(2)}_{11} + {Y}^{(2)}_{22}, \quad
%   C_3 = -\det \textrm{Y}^{(1)},\quad
%   C_4 = -\det \textrm{Y}^{(2)},
%\]
%\[
%   C_5 =  {Y}^{(1)}_{21} {Y}^{(2)}_{12} + {Y}^{(1)}_{12} {Y}^{(2)}_{21}  - {Y}^{(1)}_{11} {Y}^{(2)}_{22} -{Y}^{(1)}_{22} {Y}^{(2)}_{11}.
%\]
For the Example \ref{ex:FEM-stuff}, direct computations show that
$\trace(\mathcal{Y})=0$, so that $D$ reduces to:
\[
{D}(\boldsymbol{\xi}) = \left(\delta_\lambda \right)^2 - C_1 \left(\delta_{\xi_1} \right)^2 - C_2 \left(\delta_{\xi_2} \right)^2 - C_3 \delta_{\xi_1} \delta_{\xi_2},
\] 
where the coefficients $C_1, C_2, C_3$ are given by 
%\begin{widetext}
\begin{align*}
	%  C_1 = {Y}^{(1)}_{11} + {Y}^{(1)}_{22},\quad
	% C_2 = {Y}^{(2)}_{11} + {Y}^{(2)}_{22}, \quad
	C_1 &= -\det(\textrm{Y}^{(1)}),\quad
	C_2 = -\det(\textrm{Y}^{(2)}), \\
	C_3 &=  -{Y}^{(1)}_{11} {Y}^{(2)}_{22}-{Y}^{(1)}_{22} {Y}^{(2)}_{11}+{Y}^{(1)}_{21} {Y}^{(2)}_{12} + {Y}^{(1)}_{12} {Y}^{(2)}_{21}.
\end{align*}
%\end{widetext}
Upon introducing the change of variables $\bdxi \leftrightarrow \bdzeta$, so that
$ \boldsymbol{\xi}  - \boldsymbol{\xi}^\star = {\Psi}^\star \boldsymbol{\zeta}$, where
\begin{eqnarray*}
	\boldsymbol{\zeta} = \begin{pmatrix}
		\zeta_1 \\
		\zeta_2
	\end{pmatrix},
	\Psi^\star = \begin{pmatrix}
		2 \left(4 C_1  - C_3^2 / C_2\right)^{-1/2}  & 0  \\
		-C_3/C_2 \left(4 C_1  - C_3^2 / C_2\right)^{-1/2} & C_2^{-1/2}
	\end{pmatrix},
\end{eqnarray*}
%{\color{red}(\textbf{AIK}: please adapt this last matrix to current notations, I have not done it due to being lazy)}{\color{green} done} 
we can obtain a canonical form for ${D}(\boldsymbol{\xi})$:
\begin{equation}
	D = \left(\delta_\lambda\right) ^2 - \zeta_1^2 - \zeta_2^2.
\end{equation}
One can then rewrite $\mathcal{Y}$ as $\mathcal{Y} = \zeta_1 { \rm \tilde Y}^{(1)} + \zeta_2 { \rm \tilde Y}^{(2)}$, and check directly that ${\tilde Y}^{(1,2)}$ have the same structure as in (\ref{Yanzats}). %\RED{Should we introduce $\rho$ here?}
%
%The resulting field near Dirac point could be written as a double integral over $\boldsymbol{\zeta}$ by applying (\ref{eq:change of varibales xi zetta}) to (\ref{eq:inverse-BF-transform}):
%\begin{equation}\label{eq:discrete U result}
%   U(\boldsymbol{m}) = e^{i \boldsymbol{\xi}^* \boldsymbol{m}} \iint\limits_{-\pi}^{\pi}  \left( a_1(\boldsymbol{\zeta}) \hat{U}_1(\boldsymbol{\xi}^*) + a_2(\boldsymbol{\zeta}) \hat{U}_2(\boldsymbol{\xi}^*) \right) e^{i \boldsymbol{\zeta} \boldsymbol{\alpha}} d\boldsymbol{\zeta},
%\end{equation}
%where $\boldsymbol{\alpha} = \boldsymbol{\Psi}^T \boldsymbol{m}$ and $d\boldsymbol{\zeta} = d\zeta_1 d\zeta_2$.

%{\color{red}(\textbf{AIK}: please adapt this next red section to current notations, lots of things are weird with (C.14), missing dot in exponential, missing mionus sign in the exponential, what is $V^{(j)}$, problematic transpose sign, what exactly is the vector ${\rm U}_{loc}$ (I think it is an approximation to the vector ${\rm U}_{\bdm}$ I introduced earlier?), etc)}{\color{green} done}

Finally, estimating (\ref{eq:discret_inverse_floquet}) in the far field, we obtain
\begin{equation}
	\label{eq:U_loc_dirac}
	{\rm \hat  U}_{\text{loc}}(\tmmathbf{m}) \approx \frac{e^{-i \boldsymbol{m}\cdot\boldsymbol{\xi}^\star } \det {\Psi}^\star}{4 \pi^2} \times \sum_{i,j=1}^{2}{\rm \hat V}_i{\mathcal I}_{i,j}{\rm\tilde  V}_j,
\end{equation}
where
\begin{align*}
	\mathcal{I}_{11} &= -\tilde{Y}^{(1)}_{11} \mathcal{J}_1(\tmmathbf{\alpha},\delta_{\lambda}) -\tilde{Y}^{(2)}_{11}  \mathcal{J}_2(\tmmathbf{\alpha},\delta_{\lambda}) -\delta_{\lambda}\mathcal{J}_0(\tmmathbf{\alpha},\delta_{\lambda}),\\ 
	\mathcal{I}_{12} &=
	-(\tilde{Y}^{(1)}_{21})^\conj\mathcal{J}_1(\tmmathbf{\alpha},\delta_{\lambda}) - (\tilde{Y}^{(2)}_{21})^\conj\mathcal{J}_2(\tmmathbf{\alpha},\delta_{\lambda}),\\
	\mathcal{I}_{21} &= -\tilde{Y}^{(1)}_{21}\mathcal{J}_1(\tmmathbf{\alpha},\delta_{\lambda}) - {(\tilde{Y}^{(2)}_{21})}\mathcal{J}_2(\tmmathbf{\alpha},\delta_{\lambda}),\\
	\mathcal{I}_{22} &= \tilde{Y}^{(1)}_{11} \mathcal{J}_1(\tmmathbf{\alpha},\delta_{\lambda}) +\tilde{Y}^{(2)}_{11}  \mathcal{J}_2 -\delta_{\lambda}\mathcal{J}_0(\tmmathbf{\alpha},\delta_{\lambda}),
\end{align*}
%			\[
%			\mathcal{I}_{12} =
%			-(\tilde{Y}^{(1)}_{21})^\conj\mathcal{J}_1(\tmmathbf{\alpha},\delta_{\lambda}) - (\tilde{Y}^{(2)}_{21})^\conj\mathcal{J}_2(\tmmathbf{\alpha},\delta_{\lambda}),\quad \mathcal{I}_{21} = -\tilde{Y}^{(1)}_{21}\mathcal{J}_1(\tmmathbf{\alpha},\delta_{\lambda}) - {(\tilde{Y}^{(2)}_{21})}\mathcal{J}_2(\tmmathbf{\alpha},\delta_{\lambda}),
%			\]
%			\[
%			\quad \mathcal{I}_{22} = \tilde{Y}^{(1)}_{11} \mathcal{J}_1(\tmmathbf{\alpha},\delta_{\lambda}) +\tilde{Y}^{(2)}_{11}  \mathcal{J}_2 -\delta_{\lambda}\mathcal{J}_0(\tmmathbf{\alpha},\delta_{\lambda}),
%			\]
where
$
\tmmathbf{\alpha}=
(\Psi^{\star})^{\transpose} \tmmathbf{m}, 
$
and the $\mathcal{J}_i$ are the canonical integrals defined in  (\ref{eq:canon_int_Dirac}). Here ${\rm \hat  U}_{\text{loc}}(\tmmathbf{m}) $  is a vector of size $\hat N$ organized in a way similar to (\ref{eq:hat_F}) for a given $\tmmathbf{m}$, with entries  tending to  $u(\tmmathbf{m},j)$  as $|\tmmathbf{m}|\to \infty$.

\begin{figure*}[h]
		\centering{
			\includegraphics[width=0.45\textwidth]{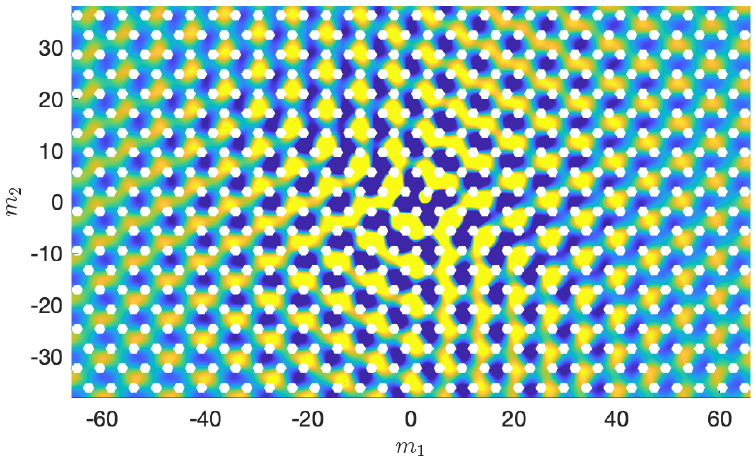}
			\includegraphics[width=0.45\textwidth]{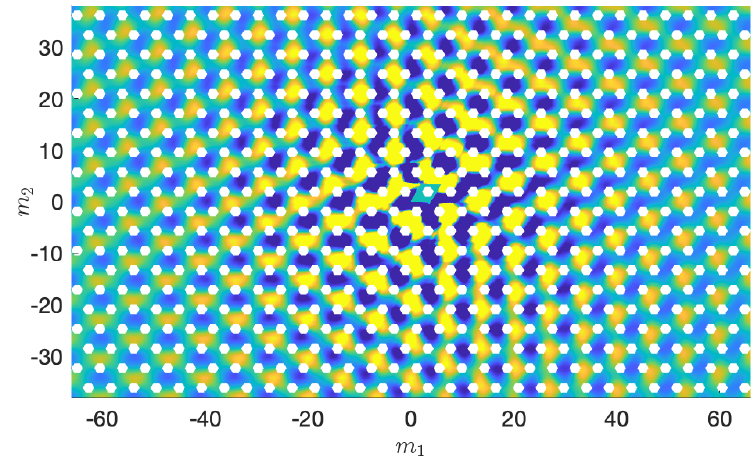}
			\caption{Direct modelling (left) and asymptotic estimation using (\ref{eq:U_loc_dirac}) (right)}
			\label{fig:FEM_numerics}
		}	
	\end{figure*}

Let us illustrate (\ref{eq:U_loc_dirac}) numerically. One way to evaluate ${\rm \hat U}_{\text{loc}}$ directly is to take a large enough domain of cells and consider a time domain problem with a harmonic point source:
\begin{equation}
	{\rm K}{\rm U}(t) + \frac{1}{c^2}{\rm M}\frac{d^2}{dt^2}{\rm U}(t) = e^{\RED{-}i(\omega^\star +\delta\omega) t}E,
\end{equation}  
where ${\rm K}$ and ${\rm M}$ are the standard mass and stiffness matrices, $\omega^\star$ is the frequency corresponding to the Dirac point $k^\star = \omega^\star/c$, $\delta\omega$ is a small perturbation, and $c$ is the wave velocity. The latter is a system of ordinary differential equations, and can be solved numerically using some time integration scheme such as the 4$^\text{th}$ order Runge-Kutta method. Then,
$
{\rm U}(t)e^{\RED{i}\omega^\star t}
$
provides an approximation to ${\rm \hat U}_{\text{loc}}$ at the nodes where the wave process has stabilized, that are far enough from the domain boundary and not affected by the reflected field. The result (real part of the field) of the time domain modelling along with the asymptotic estimation (\ref{eq:U_loc_dirac})  is presented in Figure~\ref{fig:FEM_numerics}. The geometry of one cell is presented in Figure~\ref{fig:FEM_unit_cell}, left. The crystal was excited gradually by a point source in  the $\tmmathbf{0}^{\rm th}$ cell. If we introduce a numbering in the cell by a pair of indices going along the sides of the cell, and let the left bottom node be indexed by $(0,0)$, then the node $(6,4)$ was excited. The  excitation has been turned on gradually in order not to excite higher modes of the cell. Time domain modelling has been conducted for a mesh of $10^4$ cells, with $3400$ steps in time, and time step and mesh step equal to $0.1$ and $0.55$ respectively. The other parameters were as follows:
$c = 1$, $\delta_\lambda=10^{-3}$ and $\omega^{\star} = 0.8156$.

%\onecolumngrid			
%\begin{center}
%\end{center}
%\twocolumngrid

%%% Bibliography %%%
{\small
	\bibliography{biblio}
	\bibliographystyle{unsrt}
}

\end{document}